\documentclass[a4paper,reqno,12pt,oneside]{amsart}
\usepackage{amsmath,amsrefs,geometry}
\usepackage{color}

\numberwithin{equation}{section}

\DeclareMathOperator{\divergence}{div}
\DeclareMathOperator{\Scal}{Scal}
\DeclareMathOperator{\Ric}{Ric}
\DeclareMathOperator{\Weyl}{Weyl}
\DeclareMathOperator{\Eucl}{Eucl}
\DeclareMathOperator{\Vol}{Vol}
\DeclareMathOperator{\inward}{in}
\DeclareMathOperator{\outward}{out}

\newcommand{\<}{\left<}
\renewcommand{\>}{\right>}
\renewcommand{\(}{\left(}
\renewcommand{\)}{\right)}
\renewcommand{\[}{\left[}
\renewcommand{\]}{\right]}

\providecommand{\bgl}[1]{\mbox{\boldmath$#1$}}

\usepackage{ulem}

\newtheorem{theorem}{Theorem}[section]
\newtheorem{proposition}[theorem]{Proposition}
\newtheorem{lemma}[theorem]{Lemma}

\begin{document}

\date{\today}

\title[The effect of linear perturbations on the Yamabe problem]{The effect of linear perturbations on\\the Yamabe problem} 

\author{Pierpaolo Esposito}
\address{Pierpaolo Esposito, Dipartimento di Matematica, Universit\`a degli Studi ``Roma Tre'',\,Largo S. Leonardo Murialdo 1, 00146 Roma, Italy}
\email{esposito@mat.uniroma3.it}

\author{Angela Pistoia}
\address{Angela Pistoia, Dipartimento SBAI, Universit\`a di Roma ``La Sapienza'', via Antonio Scarpa 16, 00161 Roma, Italy}
\email{pistoia@dmmm.uniroma1.it}

\author{J\'er\^ome V\'etois}
\address{J\'er\^ome V\'etois, Universit\'e de Nice - Sophia Antipolis, Laboratoire J.-A. Dieudonn\'e, UMR CNRS-UNS 7351, Parc Valrose, 06108 Nice Cedex 2, France}
\email{vetois@unice.fr}

\thanks{The first author was partially supported by the Prin project ``Critical Point Theory and Perturbative Methods for Nonlinear Differential Equations" and the Firb-Ideas project ``Analysis and Beyond". The third author was partially supported by the ANR grant ANR-08-BLAN-0335-01.}

\begin{abstract}
In conformal geometry, the Compactness Conjecture asserts that the set of Yamabe metrics on a smooth, compact, aspherical  Riemannian manifold $\(M,g\)$  is compact. Established in the locally conformally flat case by Schoen~\cites{Sch2,Sch3} and for $n\leq 24$ by Khuri--Marques--Schoen~\cites{KhuMarSch}, it has revealed to be generally false for $n\geq 25$ as shown by Brendle~\cite{Bre} and Brendle--Marques~\cite{BreMar}. A stronger version of it, the compactness under perturbations of the Yamabe equation, is addressed here with respect to the linear geometric potential $\frac{n-2}{4(n-1)} \Scal_g$, $\Scal_g$ being the Scalar curvature of $\(M,g\)$. We show that a-priori $L^\infty$--bounds fail for linear perturbations  on all manifolds with $n\geq 4$ as well as a-priori gradient $L^2$--bounds fail for non-locally conformally flat manifolds with $n\ge6$ and for locally conformally flat manifolds with $n\ge7$. In several situations, the results are optimal.   Our proof combines a finite dimensional reduction and the construction of a suitable ansatz for the solutions  generated by a family of varying metrics in the conformal class of $g$.
\end{abstract}

\maketitle

\section{Introduction}\label{Sec1}

Letting $\(M,g\)$ be a smooth, compact Riemannian $n$--manifold, $n\geq 3$, we are concerned with the so-called Yamabe equation
\begin{equation}\label{Yam}
\Delta_g u+ \alpha_n \Scal_g u=c u^{2^*-1},\quad u>0\quad\text{in }M\,,
\end{equation}
where $\Delta_g:=-\divergence_g\nabla$ is the Laplace--Beltrami operator, $\alpha_n:=\frac{n-2}{4(n-1)}$, $\Scal_g$ is the Scalar curvature of the manifold, $2^*=\frac{2n}{n-2}$ is the critical Sobolev exponent, and $c \in \mathbb{R}$. The geometric problem of finding a metric $\widetilde g$ in the conformal class $[g]=\{\phi g:\, \phi \in C^\infty(M),\, \phi>0 \}$ of $g$ with constant Scalar curvature is equivalent to solving \eqref{Yam} through the setting $\widetilde g=u^{4/\(n-2\)} g$. The constant $c$ can be restricted to the values $-1/1$ or $0$ depending on whether the {\it Yamabe invariant}
$$\mu_g(M)=\inf_{\widetilde g\in\[g\]}\(\Vol_{\widetilde g}\(M\)^{\frac{2-n}{n}}\int_M\Scal_{\widetilde g}dv_{\widetilde g}\)$$
of $\(M,g\)$ has negative/positive sign or vanishes, respectively, where $\Vol_{\widetilde g}\(M\)$ is the volume of the manifold $\(M,\widetilde g\)$.

\medskip
The Yamabe problem, raised by H.~Yamabe~\cite{Yam} in '60, was firstly solved by Trudinger~\cite{Tru} when $\mu_g(M)\leq 0$. In this case, the solution is unique (up to a normalization when $\mu_g(M)=0$). Aubin~\cite{Aub2} then solved the Yamabe problem in the non-locally conformally flat (non-l.c.f.~for short) case with $n\geq 6$, and Schoen~\cite{Sch1} solved it in the remaining cases of low dimensions $3\le n\le5$ and locally conformally flat (l.c.f.~for short) manifolds. In this paper, we restrict our attention to the case where $\(M,g\)$ has {\it positive Yamabe invariant} $\mu_g(M)>0$, and we set $c=1$ in \eqref{Yam}.

\medskip
After the complete resolution of the Yamabe problem, one can attempt to describe the solution set of \eqref{Yam}. A well-known conjecture claims the compactness of Yamabe metrics for manifolds $\(M,g\)$ which are not conformally equivalent to $(\mathbb{S}^n,g_0)$ ($\(M,g\)\not=(\mathbb{S}^n,g_0)$ for short), namely the convergence, up to a subsequence, in $C^2\(M\)$ of any sequence of solutions of equation \eqref{Yam}. Referred to in literature as the Compactness Conjecture, by elliptic regularity theory, it amounts to prove a-priori $L^\infty$--bounds on the set of solutions to \eqref{Yam}. In the basic example of the round sphere $(\mathbb{S}^n,g_0)$, by the works of Lelong--Ferrand~\cite{Lel} and Obata~\cite{Oba}, the set of solutions of \eqref{Yam} is explicit and non-compact. The Compactness Conjecture arose after the work of Schoen~\cites{Sch2,Sch3} who first proved it  in the l.c.f.~case, and also proposed a strategy, based on the Pohozaev identity, for proving it in the non-l.c.f.~case. The Compactness Conjecture has been then proved in the low-dimensional case by Li--Zhu~\cite{LiZhu} for $n = 3$, by  Druet~\cite{Dru2} for $n\leq 5$, by Marques~\cite{Mar} for $n \leq 7$, by Li--Zhang~\cites{LiZha2,LiZha3} for $n \leq 11$, and finally, by Khuri--Marques--Schoen~\cite{KhuMarSch} for $n\le24$.  Unexpectedly, the dimension $n=24$, which arises in~\cite{KhuMarSch} as the maximal dimension for a suitable quadratic form to be positive definite, has revealed to be optimal by the counter-examples constructed in dimensions $n\ge25$ by Brendle~\cite{Bre} and Brendle--Marques~\cite{BreMar}. We also refer to the constructions by Ambrosetti--Malchiodi~\cite{AmbMal} and Berti--Malchiodi~\cite{BerMal}  in case of background metrics which have a finite regularity. All these constructions are  made on $(\mathbb{S}^n,g)$ with a metric $g$ close to the round one $g_0$.

\medskip
Replacing the geometric potential $\alpha_n \Scal_g$ in \eqref{Yam} with a general potential   $\kappa\in C^1\(M\)$ such that $\kappa\(\xi\)\not=\alpha_n \Scal_g\(\xi\)$ for all $\xi\in M$, the Compactness Conjecture is esentially still true for solutions with bounded energy of
\begin{equation}\label{Yam1}
\Delta_g u+\kappa u=u^{2^*-1},\quad u>0\quad\text{in }M\,,
\end{equation}
provided that $n\ge4$, as shown by Druet~\cite{Dru1}*{Section~4}.

\medskip
Towards a deeper understanding of the Yamabe equation, one can address a stronger version of the Compactness Conjecture. One can ask whether or not the compactness property is preserved under perturbations of the equation under exam, which is equivalent to have or not uniform a-priori $L^\infty$--bounds for solutions of the perturbed problem. This question has been introduced and investigated in Druet~\cites{Dru1,Dru2}, Druet--Hebey~\cites{DruHeb1,DruHeb2}, Druet--Hebey--Robert~\cite{DruHebRob}, and, under the name of stability, in Druet--Hebey~\cite{DruHeb3} and Druet--Hebey--Vetois~\cite{DruHebVet}. The aim of our paper is to establish non-compactness properties for linear perturbations of the geometric potential $\alpha_n \Scal_g$ in \eqref{Yam}. In case $3\leq n\leq 5$, Druet~\cite{Dru2} obtained uniform $L^\infty$--bounds for solutions of \eqref{Yam1} along potentials $\kappa_\varepsilon\leq \alpha_n \Scal_g$, $\kappa_\varepsilon \to \kappa$  in $C^2\(M\)$ as $\varepsilon \to 0$, with $\(M,g\)\not=(\mathbb{S}^n,g_0)$ in case $\kappa=\alpha_n \Scal_g$. The same result is  strongly expected to be true in the l.c.f.~case and generally for $n \leq 24$, and it is still true, as shown by Druet~\cite{Dru2}, when $n\geq 6$ and $\kappa<\alpha_n \Scal_g$. In dimension $n=3$,  Li--Zhu~\cite{LiZhu}*{Theorem~0.3} have obtained uniform $L^\infty$--bounds in case $\kappa_\varepsilon\leq \alpha_n\Scal_g+\eta_0$, for some $\eta_0>0$  when $\(M,g\)\not=(\mathbb{S}^n,g_0)$.

\medskip
Let us briefly review the previous results of non-compactness for equations of type \eqref{Yam1}. Apart from the trivial case of the Yamabe equation on $(\mathbb{S}^n,g_0)$, the first non-compactness result is due to Hebey--Vaugon~\cite{HebVau} who proved that in the l.c.f.~case with $n\ge4$, there always exists $\widetilde{g}\in\[g\]$ such that the equation $\Delta_{\widetilde{g}}u+\alpha_n\max_M( \Scal_{\widetilde{g}})u=u^{2^*-1}$ in $M$ is not compact. Families of non-compact solutions have then been explicitly constructed on $(\mathbb{S}^n,g_0)$ by Druet~\cite{Dru1} and Druet--Hebey~\cite{DruHeb1} for linear perturbations of the potential $\kappa=\alpha_n\Scal_g$ when $n\ge6$ and in case $\kappa>\alpha_n\Scal_g$ when $n=6$ (see also the survey paper by Druet--Hebey~\cite{DruHeb2} for the case of quotients of $(\mathbb{S}^n,g_0)$). When $(\kappa-\alpha_n\Scal_g)$ is a positive constant, on $\(M,g\)=(\mathbb{S}^n,g_0)$, Chen--Wei--Yan~\cite{ChenWeiYan} have constructed infinitely many solutions with unbounded energy when $n\geq 5$, and Hebey--Wei~\cite{HebWei} have constructed non-compact solutions with bounded energy for an infinite number of constant $\kappa$ in case $n=3$. Concerning the potential $\kappa=\alpha_n\Scal_g$, so far, the only available examples of non-compact solutions for \eqref{Yam} and its linear perturbations are in the case of $(\mathbb{S}^n,g)$, with $g$ close to $g_0$ and $\kappa=\alpha_n \Scal_g$ by~Ambrosetti--Malchiodi~\cite{AmbMal}, Berti--Malchiodi~\cite{BerMal}, Brendle~\cite{Bre}, and Brendle--Marques~\cite{BreMar}, or with $g=g_0$ and $\kappa$ close to $\alpha_n \Scal_{g_0}$ by Druet~\cite{Dru1} and Druet--Hebey~\cite{DruHeb1} (see also Druet--Hebey~\cite{DruHeb2} for the case of quotients of $(\mathbb{S}^n,g_0)$).

\medskip
In this paper, for $n\geq 4$, we exhibit   the general failure of compactness properties for
\begin{equation}\label{Eq1}
\Delta_gu+(\alpha_n\Scal_g+\varepsilon h) u=u^{2^*-1},\quad u>0\quad\text{in }M\,,
\end{equation}
where $h$ is a  $C^1$ or $C^{0,\alpha}$--function in $M$, $\alpha\in\(0,1\)$,  with $\max_Mh>0$ and $\varepsilon>0$ is a small parameter. As a by-product, we obtain that the Compactness Conjecture completely fails down under the effect of linear perturbations (with the correct sign) of the Yamabe equation \eqref{Yam} on every manifold $\(M,g\)$ with $n\geq 4$ (but it is still true for $n=3$ by Li--Zhu~\cite{LiZhu}). Our results, together with those by Druet~\cites{Dru1,Dru2}, give a sharp picture of the situation. Even more than the failure of a-priori $L^\infty$--bounds, we show that a-priori gradient $L^2$--bounds fail for non-locally conformally flat manifolds with $n\ge6$ and for locally conformally flat manifolds with $n\ge7$.

\medskip
To be more precise, we say that a family $\(u_\varepsilon\)_\varepsilon$ of solutions to equation \eqref{Eq1} {\it blows up} at some point $\xi_0\in M$ if there holds $\sup_Uu_\varepsilon\to+\infty$ as $\varepsilon\to0$, for all neighborhoods $U$ of $\xi_0$ in $M$. Let $E:M \to (-\infty,+\infty]$ be defined as
\begin{equation}\label{ThEq1}
E(\xi)= \left\{\begin{array}{ll}
h\(\xi\) A_\xi^{-\frac{2}{n-2}} &\text{if }n=4,5\text{ or }\(M,g\)\text{ l.c.f.}\\
h\(\xi\) \left|\Weyl_g\(\xi\)\right|^{-1}_g &\text{if }n\ge6\text{ and }\(M,g\)\text{ non-l.c.f.}
\end{array}\right.
\end{equation}
with the convention that $1/0=+\infty$. Here, $A_\xi$ is defined  in \eqref{Eq3} and $\Weyl_g$ is the Weyl curvature tensor of the manifold. 
In dimensions $n=3,4,5$ or if the manifold is l.c.f., up to a conformal change of metric (depending on $\xi \in M$), the Green's function $G_g(\cdot, \xi)$ has an asymptotic expansion of the form
\begin{equation}\label{Eq3}
G_g\(\exp_\xi y,\xi\)=\beta_n^{-1}\left|y\right|^{2-n}+A_\xi+\operatorname{O}\(\left|y\right|\)
\end{equation}
as $y\to0$, where $\beta_n:=\(n-2\)\omega_{n-1}$, $\omega_{n-1}$ is the volume of the unit $\(n-1\)$--sphere, and $A_\xi \in \mathbb{R}$, see Lee--Parker~\cite{LeePar}. The geometric quantity $A_\xi$ depends only on the manifold $\(M,g\)$ and the point $\xi$, is smooth with respect to $\xi$, and can be identified with the {\it mass} of a stereographic projection of the manifold with respect to $\xi$. We refer to Lee--Parker~\cite{LeePar} for the definition of the mass and a discussion about its role in general relativity. In particular, for manifolds $\(M,g\) \not=\(\mathbb{S}^n,g_0 \)$, in case $n=4,5$ and in the l.c.f.~case with $n\ge6$, we have that $A_\xi>0$ by the result of Schoen--Yau~\cites{SchYau1,SchYau2}, and thus $E(\xi)<+\infty$ for all $\xi\in M$. In the non-l.c.f.~case with $n\ge6$, we have that $\Weyl_g\not\equiv0$, and thus $E(\xi)\not\equiv+\infty$. Our first result concerns the existence of solutions blowing-up at one point and reads as:

\begin{theorem}\label{ThTh}{\bf (existence of solutions with a single blow-up point in case $n\ge4$)}
Let $\(M,g\) \not=\(\mathbb{S}^n,g_0 \)$ be a smooth compact Riemannian manifold with $n\ge4$ and $\mu_g(M)>0$, and $h\in C^{0,\alpha}\(M\)$, $\alpha\in\(0,1\)$, { be} so that $\max_Mh>0$. In the non-l.c.f.~case with $n \geq 6$, assume in addition that {$\min\{|\Weyl_g(\xi)|_g\,:\,h(\xi)>0\}>0$}. Then for $\varepsilon>0$ small, equation \eqref{Eq1} has a solution $u_\varepsilon \in C^{2,\alpha}\(M\)$ such that the family $\(u_\varepsilon\)_\varepsilon$ blows up, up to a sub-sequence, as $\varepsilon \to 0$ at some point $\xi_0$ so that $E(\xi_0)=\max_ME$.
\end{theorem}

Let us mention that the Compactness Conjecture does hold for (\ref{Yam}) when $n\geq 6$ as soon as $|\Weyl_g(\xi)|>0$ for all $\xi \in M$, as it follows by Li-Zhang~\cites{LiZha2} and Marques ~\cite{Mar}. The following result concerns multiplicity of solutions with a single blow-up point. Isolated critical points of $E$ with non-trivial local degree include non-degenerate critical points of $E$. The result reads as:

\begin{theorem}\label{ThThTh}{\bf (multiplicity of solutions with a single blow-up point in case $n\ge4$)}
Let $\(M,g\)\not=\(\mathbb{S}^n,g_0\)$ be a smooth compact Riemannian manifold  with $n\ge4$ and $\mu_g(M)>0$, and $h\in C^1\(M\)$. { For any} isolated critical point $\xi_0$ of $E$ with non-trivial local degree and $h(\xi_0)>0$, for $\varepsilon>0$ small, equation \eqref{Eq1} has a solution $u_\varepsilon\in C^{2,\alpha}\(M\)$, $\alpha\in(0,1)$, such that the family $\(u_\varepsilon\)_\varepsilon$ blows up, up to a sub-sequence, at $\xi_0$ as $\varepsilon \to 0$.
\end{theorem}

As already said, by the result of Li--Zhu~\cite{LiZhu}, such blowing-up solutions as in Theorems~\ref{ThTh} and~\ref{ThThTh} do not exist in dimension $n=3$. The last results, Theorems~\ref{ThThThThNonlcf},~\ref{ThThThTh}, and~\ref{ThThThThTh} below, claim the existence of solutions which blow up at more than one point. The first result concerns the non-l.c.f.~case with $n\geq 6$ and reads as:

\begin{theorem}\label{ThThThThNonlcf}{\bf (existence of solutions with more than one blow-up point in the non-l.c.f.~case with $n\ge6$)}
Let $\(M,g\)$ be a smooth compact non-l.c.f.~Riemannian manifold with $n\ge6$ and $\mu_g(M)>0$. Let $k\ge2$ be an integer and $h_k\in C^{0,\alpha}\(M\)$, $\alpha\in\(0,1\)$, be so that the set $\{\xi\in M\,:\,h_k\(\xi\)>0\}$ has $k$ connected components $\mathcal{C}_1,\dotsc,\mathcal{C}_k$ and $\min\{|\Weyl_g(\xi)|_g\,:\,h(\xi)\geq 0\}>0$. Then, for $\varepsilon>0$ small, equation \eqref{Eq1} has a solution $u_{k,\varepsilon}\in C^{2,\alpha}\(M\)$ such that the family $\(u_{k,\varepsilon}\)_\varepsilon$ blows up, up to a sub-sequence, as $\varepsilon \to 0$ at $k$ distinct points $\(\xi_0\)_1,\dots,\(\xi_0\)_k$ so that 
$$\frac{h_k\((\xi_0)_j\)}{\left|\Weyl_g \((\xi_0)_j\) \right|_g}  =\max_{\,\xi \in \mathcal{C}_j} \frac{h_k \(\xi\) }{\left|\Weyl_g \(\xi\)\right|_g}$$ 
for all $j=1,\dotsc,k$. Moreover, there holds $\lim_{k\to+\infty}\limsup_{\varepsilon\to0}\left\|\nabla u_{k,\varepsilon}\right\|_{L^2\(M\)}=+\infty$.
\end{theorem}

Each blow-up point $\(\xi_0\)_j$, $j=1,\dotsc,k$, in Theorem~\ref{ThThThThNonlcf} maximizes the same function $E$ as in Theorem~\ref{ThTh} for single blow-up points. On the contrary, in the remaining cases, there is a strong interaction between the blow-up points, and these are not anymore related to maximum points of the function $E$ defined in \eqref{ThEq1}. Concerning the l.c.f.~case with $n\geq 7$, we prove the following result:

\begin{theorem}\label{ThThThTh}{\bf (existence of solutions with more than one blow-up point in the l.c.f.~case with $n\ge7$)}
Let $\(M,g\) \not=\(\mathbb{S}^n,g_0 \)$ be a smooth compact l.c.f.~Riemannian manifold with $n\geq 7$ and $\mu_g(M)>0$, and $h\in C^{0,\alpha}\(M\)$, $\alpha\in\(0,1\)$, be so that $\max_Mh>0$. Then for any integer {$k\ge2$}, for $\varepsilon>0$ small, equation \eqref{Eq1} has a solution $u_{k,\varepsilon}\in C^{2,\alpha}\(M\)$ such that the family $\(u_{k,\varepsilon}\)_\varepsilon$ blows up, up to a sub-sequence, at $k$ distinct points $\(\xi_0\)_1,\dots,\(\xi_0\)_k$ in $M$  as $\varepsilon \to 0$. Moreover, there holds $\lim_{k\to+\infty}\limsup_{\varepsilon\to0}\left\|\nabla u_{k,\varepsilon}\right\|_{L^2\(M\)}=+\infty$.
\end{theorem}

The location of $\xi_0=(\(\xi_0\)_1,\dots,\(\xi_0\)_k)$ in Theorem~\ref{ThThThTh} is related to maximum points of a ``reduced energy'' given in \eqref{eridk2}, and the assumption $n\geq 7$ guarantees that such ``reduced energy'' achieves its maximum value.

\medskip
Our last result, Theorem~\ref{ThThThThTh} below, concerns the l.c.f.~case for $n=6$. This case reveals to be even more intricate than the case of higher dimensions. For any integer $k\ge2$, define $\Delta_k:=\left\{(\xi_1,\dots,\xi_k)\in M^k \ :\ \xi_i=\xi_j\ \hbox{for}\ i\not=j\right\}$. For any $\bgl\xi:=(\xi_1,\dots,\xi_k)\in M^k\backslash\Delta_k$, let $A_{k,\bgl\xi}$ be the symmetric $k\times k$ matrix with entries
\begin{equation}\label{matrix1}
\big(A_{k,\bgl\xi}\big)_{ij}:=\left\{\begin{aligned}&A_{\xi_i}&&\text{if }i=j\\&G_g\(\xi_i,\xi_j\)&&\text{if }i\ne j,\end{aligned}\right.
\end{equation}
where $A_{\xi_i}$ is as in \eqref{Eq3}. When $A_{k,\bgl\xi}$ is invertible, let $E_k:M^k\backslash\Delta_k\to\mathbb{R}$ be defined as
\begin{equation}\label{matrix2}
E_k\(\bgl\xi\):=\big<H\big(\bgl\xi\big),A_{k,\bgl\xi}^{-1}.H\big(\bgl\xi\big)\big>,
\end{equation}
where $H\(\bgl\xi\):=\(h\(\xi_1\),\dotsc,h\(\xi_k\)\)$ and $\<\cdot,\cdot\>$ is the Euclidean scalar product. { Here, contrary to the situation with $k=1$, the definition of $E_k\(\bgl\xi\)$ allows to consider the case $\(M,g\)=\(\mathbb{S}^n,g_0\)$. In this case, $A_{k,\bgl\xi}$\vspace{-3pt} has all entries equal to zero on the diagonal. In particular, when $k=2$ or 3, observe that $A_{k,\bgl\xi}$ is invertible, and thus $E_k\(\bgl\xi\)$ is well-defined, for all $\bgl\xi\in\(\mathbb{S}^n\)^k\backslash\Delta_k$.} { Our result in the l.c.f.~case with $n=6$ reads as:}

\begin{theorem}\label{ThThThThTh}{\bf (existence of solutions with more than one blow-up point in the l.c.f.~case with $n=6$)}
Let $\(M,g\)$ be a smooth compact l.c.f.~Riemannian manifold with $n=6$ and $\mu_g(M)>0$,  and $h \in C^1\(M\)$. Let $k\ge2$ be an integer, and assume that $E_k$ has an isolated critical point $\bgl\xi_0:=\(\(\xi_0\)_1,\dots,\(\xi_0\)_k\)$ with non-trivial local degree and $A_{k,\bgl\xi_0}^{-1}.H\big(\bgl\xi_0\big)$ has positive coordinates. Then for $\varepsilon>0$ small, equation \eqref{Eq1} has a solution {$u_{k,\varepsilon}\in C^{2,\alpha}\(M\)$}, $\alpha \in (0,1)$, such that the family {$\(u_{k,\varepsilon}\)_\varepsilon$} blows up, up to a sub-sequence, at $\(\xi_0\)_1,\dots,\(\xi_0\)_k$ as $\varepsilon \to 0$.   
\end{theorem}

Contrary to the assumptions in the previous theorems, here it seems unclear in general when the function $E_k$ admits an isolated critical point with non-trivial local degree. An easy situation where we can construct $h$ and $\bgl\xi_0$ satisfying the assumptions in Theorem~\ref{ThThThThTh} is the case $k=2$ on the round sphere $\(\mathbb{S}^n,g_0\)$. Indeed, in this case, we find that $E_2\(\bgl\xi\)=2h\(\xi_1\)h\(\xi_2\)G_{g_0}\(\xi_1,\xi_2\)^{-1}$ for all $\bgl\xi=\(\xi_1,\xi_2\)$, $\xi_1\ne\xi_2$, which has always a maximum point in $(\mathbb{S}^n \cap \{ h\geq 0\})^2 \setminus  \Delta_2$ provided that $\max_M h >0$. It is clear that the maximum point is non-degenerate for several $h'$s (in a generic sense).

\medskip
Let us finally  compare problem \eqref{Eq1} with its Euclidean counter-part on a smooth bounded domain $\Omega \subset \mathbb{R}^n$, $n\ge4$, with homogeneous Dirichlet boundary condition
\begin{equation}\label{bn}
\Delta_{\Eucl} u+\lambda u=u^{2^*-1}  \ \hbox{in}\ \Omega,\quad u>0 \ \hbox{in}\ \Omega,\quad u=0\ \hbox{on}\ \partial\Omega.
\end{equation}
For $\lambda \geq 0$, a direct minimization method (for the corresponding Rayleigh quotient) never gives rise to any solution of \eqref{bn}, and moreover, no solution exists at all if $\Omega$ is star-shaped as shown by Poho{\v{z}}aev~\cite{poho}. Moreover, following the arguments developed by Ben~Ayed--El~Mehdi--Grossi--Rey~\cite{begr}, problem  \eqref{bn} has no solutions with a single blow-up point as $\lambda \to 0^+$. The effect of the geometry, which is crucial to provide a solution for the Yamabe problem (corresponding to $\lambda=0$ in \eqref{bn}) by minimization, is also relevant to producing solutions of \eqref{Eq1} (corresponding to $\lambda \to 0^+$ in \eqref{bn}) with a single blow-up point as stated in Theorems~\ref{ThTh} and~\ref{ThThTh}. On the contrary, Theorem~\ref{ThThThTh} has a partial counter-part on domains with nontrivial topology, see Musso--Pistoia~\cite{mupi1} and Pistoia--Rey~\cite{pire}. When $\lambda<0$, solutions of \eqref{bn} can be found by direct minimization as shown by Brezis--Nirenberg~\cite{bn}, and exhibit a single blow-up point as $\lambda \to 0^-$ as shown by Han~\cite{h}, in contrast with the compactness property proved by Druet~\cites{Dru1}. Solutions of \eqref{bn} with a single blow-up point, see Rey~\cites{Rey,rey2}, and with multiple blow-up points, see Bahri--Li--Rey~\cites{blr} and Musso--Pistoia~\cite{mupi4}, as $\lambda \to 0^-$ have been constructed in a very general way. Since the manifold with boundary $(\Omega,dx)$ is l.c.f., notice that the Green's function $G(\cdot,\xi)$ still has an asymptotic expansion of the form \eqref{Eq3}, but the constant $A_\xi$ is always negative in this case. The different sign of $A_\xi$ is the analytical reason of the completely different picture we have for equations \eqref{Eq1} and \eqref{bn}.

\medskip
The paper is organized in the following way. In Section~\ref{Sec2} we describe the perturbative method we use to attack existence issues of blowing-up solutions. We describe the main steps of such an approach, leading to the general result Theorem~\ref{Th}, and we deduce from it our main results concerning solutions with a single blow-up point.   A crucial  point is to produce a suitable ansatz for the solutions. Inspired by the approach of Lee--Parker~\cite{LeePar}, which unifies the previous constructions of Aubin~\cite{Aub2} and Schoen~\cite{Sch1} in the resolution of the Yamabe problem, we build up general approximating solutions $W_{\varepsilon,t,\xi}$ for \eqref{Eq1} which approximation rates are evaluated in Section~\ref{Sec3}. An   important  point here is  that we  allow the metric $g$ to vary in the conformal class so to gain flatness at each point $\xi \in M$. An alternative, less geometric approach can be devised in the non-l.c.f.~case, see Esposito--Pistoia--V\'etois~\cite{EPV}, by keeping $g$ fixed and slightly correcting the basic ansatz via linearization so to account for the local geometry. Thanks to the solvability theory of the linearized operator for \eqref{Eq1} at $W_{\varepsilon,t,\xi}$, we are led to study critical points of a finite-dimensional functional $\mathcal{J}_\varepsilon\(t,\xi\)$. A key step is to obtain an asymptotic expansion of $\mathcal{J}_\varepsilon\(t,\xi\)$ and to identify a ``reduced energy'' as the main order term. This step is performed in Section~\ref{Sec4} in $C^0$--sense and completes the proof of  Theorem~\ref{ThTh}. In Section~\ref{Sec5}, we investigate the existence of solutions with $k$ blow-up points, yielding to the proofs of Theorems~\ref{ThThThThNonlcf},~\ref{ThThThTh}, and~\ref{ThThThThTh}. The $C^1$--expansion of $\mathcal{J}_\varepsilon\(t,\xi\)$ is addressed in Section~\ref{Sec6}, completing the proof of Theorem~\ref{ThThTh}. The appendix is devoted to some technical issues.

\medskip\noindent 
{\bf Acknowledgments:} this work has been initiated and partially carried out during the visits of the third author to the University of ``Roma La Sapienza'' in November~2009 and to the University of ``Roma Tre'' in November~2010. The third author gratefully acknowledges the hospitality and the financial support of these two institutions.

\section{Scheme of the proof}\label{Sec2}

In this section, we aim to give the scheme of proof for Theorem~\ref{Th} below. First, let us set some notations. We denote the conformal Laplacian of the manifold by
\begin{equation}\label{Eq4}
L_g:=\Delta_g+\alpha_n\Scal_g\,,
\end{equation}
where $\alpha_n:=\frac{n-2}{4(n-1)}$. The conformal covariance of $L_g$ expresses as
\begin{equation}\label{confinv}
L_{\hat g}(\phi)= u^{-(2^*-1)} L_g (u \phi)\qquad \forall \,\phi \in C^2(M),
\end{equation}
for all $\hat g=u^{2^*-2} g$ in the conformal class $[g]$ of $g$. We assume that the manifold has {\it positive Yamabe invariant}, which is equivalent to assuming the coercivity of $L_g$. We let $H^2_1\(M\)$ be the Sobolev space of all functions in $L^2\(M\)$ with gradient in $L^2\(M\)$ equipped with the scalar product
\begin{equation}\label{Eq5}
\<u,v\>_{L_g}:=\int_M\<\nabla u,\nabla v\>_gdv_g+\alpha_n\int_M\Scal_guvdv_g\,,
\end{equation}
where $dv_g$ is the volume element of the manifold. We let $\left\|\cdot\right\|_{L_g}$ be the norm induced by $\<\cdot,\cdot\>_{L_g}$. For any $u\in L^q\(M\)$, we denote the $L^q$--norm of $u$ by $\left\|u\right\|_q:=\(\int_M|u|^qdv_g\)^{1/q}$. We define $\left\|u\right\|_{1,2}:=\(\left\|\nabla u\right\|^2_2+\left\|u\right\|^2_2\)^{1/2}$. By coercivity of $L_g$, we get that the norms $\left\|\cdot\right\|_{L_g}$ and $\left\|\cdot\right\|_{1,2}$ are equivalent.

\medskip
We let $i_g$ be the injectivity radius of the manifold $\(M,g\)$. By compactness of $M$, we get that there exists a positive real number $r_0$ such that $r_0<i_g$. In case $\(M,g\)$ is locally conformally flat, there exists a family $\(g_\xi\)_{\xi\in M}$ of smooth conformal metrics to $g$ such that $g_\xi$ is flat in the geodesic ball $B_\xi\(r_0\)$. In case $\(M,g\)$ is not locally conformally flat, we fix $N>n$, and we provide ourselves with a family $\(g_\xi\)_{\xi\in M}$ of smooth conformal metrics to $g$ such that
\begin{equation}\label{Eq6}
\left|\exp_\xi^*g_\xi\right|\(y\)=1+\operatorname{O}\big(\left|y\right|^N\big)
\end{equation}
uniformly with respect to $\xi\in M$ and $y\in T_\xi M$, $\left|y\right|\ll1$, where $\left|\exp_\xi^*g_\xi\right|$ is the determinant of $g_\xi$ in the geodesic normal coordinates of $g_\xi$ at $\xi$. Such coordinates are said to be {\it conformal normal coordinates} of order $N$ on the manifold. Here and in the sequel, the exponential map $\exp_\xi$ is always intended with respect to the metric $g_\xi$. We refer to Lee--Parker~\cite{LeePar} for a proof of the existence of conformal normal coordinates of any finite order, see also the later proofs by Cao~\cite{Cao} and G{\"u}nther~\cite{Gun} of the existence of conformal normal coordinates which are volume preserving near a given point (with no remainder term in \eqref{Eq6}). For any $\xi\in M$, we let $\varLambda_\xi$ be the smooth positive function in $M$ such that $g_\xi=\varLambda_\xi^{2^*-2}g$. In both cases (locally conformally flat or not), the metric $g_\xi$ can be chosen smooth with respect to $\xi$ and such that $\varLambda_\xi\(\xi\)=1$ and $\nabla\varLambda_\xi\(\xi\)=0$. We let $G_g$ and $G_{g_\xi}$ be the respective Green's functions of $L_g$ and $L_{g_\xi}$. Using the fact that $\varLambda_\xi\(\xi\)=1$, by \eqref{confinv}, we find that
\begin{equation}\label{Eq7}
G_g\(\cdot,\xi\)=\varLambda_\xi (\cdot)\ G_{g_\xi}\(\cdot,\xi\).
\end{equation}

\medskip
By compactness of $M$ and since $g_\xi$ is smooth with respect to $\xi$, decreasing if necessary the real number $r_0$, we may assume that $r_0<i_{g_\xi}$ for all $\xi\in M$, where $i_{g_\xi}$ is the injectivity radius of the manifold $\(M,g_\xi\)$. For $\varepsilon>0$ small and for any $t>0$, we define
\begin{equation}\label{Eq8}
\delta_\varepsilon\(t\):=\left\{\begin{aligned}
&e^{-\frac{t}{\varepsilon}}&&\text{if }n=4\\
&t\varepsilon^{\frac{1}{n-4}}&&\text{if }n=5\text{ or (}n\ge6\text{ and }\(M,g\)\text{ l.c.f.)}\\
&t\ell^{-1}\(\varepsilon\)&&\text{if }n=6\text{ and }\(M,g\)\text{ non-l.c.f.}\\
&t\sqrt\varepsilon&&\text{if }n\ge7\text{ and }\(M,g\)\text{ non-l.c.f.},
\end{aligned}\right.
\end{equation}
where $\ell:\(0,e^{-1/2}\)\to\(0,e^{-1}/2\)$ is given by $\ell\(\delta\):=-\delta^2\ln\delta$, and $\delta_\varepsilon:=\delta_\varepsilon\(1\)$. For $\varepsilon>0$ small and for any $t>0$, $\xi\in M$, inspired by the approach of Lee--Parker~\cite{LeePar}, we define $W_{\varepsilon,t,\xi}$ in $M$ by
\begin{equation}\label{Eq9}
W_{\varepsilon,t,\xi}\(x\)=G_g\(x,\xi\)\widehat{W}_{\varepsilon,t,\xi}\(x\),
\end{equation}
with
\begin{equation}\label{Eq10}
\widehat{W}_{\varepsilon,t,\xi}\(x\):=\left\{\begin{aligned}
&\beta_n\delta_\varepsilon\(t\)^{\frac{2-n}{2}}d_{g_\xi}\(x,\xi\)^{n-2}U\(\delta_\varepsilon\(t\)^{-1} d_{g_\xi}(x,\xi) \)&&\text{if }d_{g_\xi}\(x,\xi\)\le r_0\\
&\beta_n\delta_\varepsilon\(t\)^{\frac{2-n}{2}}r_0^{n-2}U\(\delta_\varepsilon\(t\)^{-1}r_0\)&&\text{if }d_{g_\xi}\(x,\xi\)>r_0,
\end{aligned}\right.
\end{equation}
where $\beta_n=(n-2)\omega_{n-1}$, $\omega_{n-1}$ is the volume of the unit $\(n-1\)$--sphere, $\delta_\varepsilon\(t\)$ is as in \eqref{Eq8}, and
\begin{equation}\label{Eq11}
U\(r\):=\(\frac{\sqrt{n\(n-2\)}}{1+r^2}\)^{\frac{n-2}{2}}.
\end{equation}
The function $U$ generates a family $U_\delta(r)=\delta^{\frac{2-n}{2}} U(\delta^{-1}r)$, $\delta>0$, of radial solutions of the equation $\Delta_{\Eucl}U =U^{2^*-1}$ in $\mathbb{R}^n$, where $\Delta_{\Eucl}:=-\divergence_{\Eucl}\nabla$ is the Laplace operator with respect to the Euclidean metric. With these definitions, by \eqref{Eq7}, we observe that $W_{\varepsilon,t,\xi}$ rewrites as
$$W_{\varepsilon,t,\xi}=\varLambda_\xi [\beta_n G_{g_\xi}(x,\xi) d_{g_\xi}(x,\xi)^{n-2}]U_{\delta_\varepsilon(t)}(d_{g_\xi}(x,\xi))$$
for all $x\in M$ such that $d_{g_\xi}(x,\xi)\leq r_0$.

\medskip
Let us spend few words to comment on the choice of $W_{\varepsilon,t,\xi}$. Since, by Lemma~\ref{Lem4}, the function $\beta_n  G_{g_\xi}(x,\xi) d_{g_\xi}(x,\xi)^{n-2}$ is very close to $1$ as $x \to \xi$, we have that $W_{\varepsilon,t,\xi}$ is a small correction of $\varLambda_\xi U_{\delta_\varepsilon(t)}(d_{g_\xi}(x,\xi))$ near $\xi$. Since, by \eqref{confinv}, we have that
$$L_g (\varLambda_\xi U_{\delta_\varepsilon(t)}(d_{g_\xi}(x,\xi)))=\varLambda_\xi^{2^*-1}L_{g_\xi}(U_{\delta_\varepsilon(t)}(d_{g_\xi}(x,\xi)))$$
in view of the flatness of $g_\xi$ at $\xi$ (see \eqref{Eq6}) and $\Delta_{\Eucl}U_\delta =U_\delta^{2^*-1}$ in $\mathbb{R}^n$ it is natural to expect that $W_{\varepsilon,t,\xi}$ is a very good approximating solution to \eqref{Eq1} near $\xi$. Away from $\xi$, the function $W_{\varepsilon,t,\xi}$ behaves like $\beta_n [n(n-2)]^{\frac{n-2}{4}} \delta_\varepsilon(t)^{\frac{n-2}{2}} G_g(x,\xi)$, which is still a good approximating solution to \eqref{Eq1} in that region.

\medskip
We define $V_0,\dotsc,V_n:\mathbb{R}^n\to\mathbb{R}$ by
\begin{equation}\label{Eq12}
V_0\(y\):=\frac{\left|y\right|^2-1}
{\(1+\left|y\right|^2\)^{\frac{n}{2}}}\quad\text{and}\quad V_i\(y\):=\frac{y_i}
{\(1+\left|y\right|^2\)^{\frac{n}{2}}}
\end{equation}
for all $y \in \mathbb{R}^n$ and $i=1,\dotsc,n$. By Bianchi--Egnell~\cite{BiaEgn}, any solution $v\in D^{1,2}\(\mathbb{R}^n\)$ to the equation $\Delta_{\Eucl}v=\(2^*-1\)U^{2^*-2}v$ is a linear combination of the functions $V_0,\dotsc,V_n$. We let $\chi$ be a smooth cutoff function in $\mathbb{R}_+$ such that $0\le\chi\le1$ in $\mathbb{R}_+$, $\chi=1$ in $\[0,r_0/2\]$, and $\chi=0$ in $\[r_0,\infty\)$. For $\varepsilon>0$ small and for any $t>0$, $\xi\in M$, and $\omega\in T_\xi M$, we define $Z_{\varepsilon,t,\xi},Z_{\varepsilon,t,\xi,\omega}:M\to\mathbb{R}$ by
\begin{equation}\label{Eq13}
Z_{\varepsilon,t,\xi}\(x\):=G_g\(x,\xi\)\widehat{Z}_{\varepsilon,t,\xi}\(x\)\quad\text{and}\quad Z_{\varepsilon,t,\xi,\omega}\(x\):=G_g\(x,\xi\)\widehat{Z}_{\varepsilon,t,\xi,\omega}\(x\),
\end{equation}
where
\begin{align}
\widehat{Z}_{\varepsilon,t,\xi}\(x\)&:=\beta_n\chi\(d_{g_\xi}\(x,\xi\)\)\delta_\varepsilon\(t\)^{\frac{2-n}{2}}d_{g_\xi}\(x,\xi\)^{n-2}V_0\(\delta_\varepsilon\(t\)^{-1}\exp_\xi^{-1}x\),\label{Eq14}\\
\widehat{Z}_{\varepsilon,t,\xi,\omega}\(x\)&:=\beta_n\chi\(d_{g_\xi}\(x,\xi\)\)\delta_\varepsilon\(t\)^{\frac{2-n}{2}}d_{g_\xi}\(x,\xi\)^{n-2}\<V\(\delta_\varepsilon\(t\)^{-1}\exp_\xi^{-1}x\),\omega\>_g\label{Eq15}
\end{align}
with $V\(y\)=\(V_1\(y\),\dotsc,V_n\(y\)\)$ for all $y\in T_\xi M\cong\mathbb{R}^n$. We define
\begin{align}
K_{\varepsilon,t,\xi}&:=\left\{\lambda Z_{\varepsilon,t,\xi}+Z_{\varepsilon,t,\xi,\omega}\,:\quad\lambda\in\mathbb{R}\quad\text{and}\quad\omega\in T_\xi M\right\},\label{Eq16}\\
K^\perp_{\varepsilon,t,\xi}&:=\left\{\phi\in H^2_1\(M\):\,\,\,\,\<\phi,Z_{\varepsilon,t,\xi}\>_{L_g}=0\,\,\,\,\text{and}\,\,\,\,\<\phi,Z_{\varepsilon,t,\xi,\omega}\>_{L_g}=0\,\,\,\,\text{for all }\omega\in T_\xi M\right\},\label{Eq17}
\end{align}
where the scalar product $\<\cdot,\cdot\>_{L_g}$ is as in \eqref{Eq5}. We let $\varPi_{\varepsilon,t,\xi}$ and $\varPi^\perp_{\varepsilon,t,\xi}$ be the respective projections of $H^2_1\(M\)$ onto $K_{\varepsilon,t,\xi}$ and $K^\perp_{\varepsilon,t,\xi}$.

\medskip
We intend to construct solutions to equation \eqref{Eq1} of the form
$$u_\varepsilon:=W_{\varepsilon,t,\xi}+\phi_\varepsilon\,,$$
where $t>0$, $\xi \in M$, $\phi_\varepsilon\in K^\perp_{\varepsilon,t,\xi}$, and $W_{\varepsilon,t,\xi}$ is as in \eqref{Eq9}. We re-write equation \eqref{Eq1} as the couple of equations
\begin{align}
&\varPi_{\varepsilon,t,\xi}\(W_{\varepsilon,t,\xi}+\phi-L_g^{-1}\(f_\varepsilon\(W_{\varepsilon,t,\xi}+\phi\)\)\)=0\,,\label{Eq19}\\
&\varPi^\perp_{\varepsilon,t,\xi}\(W_{\varepsilon,t,\xi}+\phi-L_g^{-1}\(f_\varepsilon\(W_{\varepsilon,t,\xi}+\phi\)\)\)=0\,,\label{Eq20}
\end{align}
where $L_g$ is as in \eqref{Eq4} and
\begin{equation}\label{Eq21}
f_\varepsilon\(u\):=u_+^{2^*-1}-\varepsilon hu\,,
\end{equation}
with $u_+=\max\(u,0\)$. We begin with solving equation \eqref{Eq20} in Proposition~\ref{Pr1}, a rather standard result in this context (see for instance Musso--Pistoia~\cite{mupi4}) which proof is skipped for shortness.

\begin{proposition}\label{Pr1}
Given two positive real numbers $a<b$, there exists a positive constant $C=C\(a,b,n,M,g,h\)$ such that for $\varepsilon>0$ small, for any $t\in\[a,b\]$ and $\xi\in M$, there exists a unique function $\phi_{\varepsilon,t,\xi}\in K^\perp_{\varepsilon,t,\xi}$ which solves equation \eqref{Eq20} and satisfies
\begin{equation}\label{Pr1Eq1}
\left\|\phi_{\varepsilon,t,\xi}\right\|_{1,2}\le C\left\|R_{\varepsilon,t,\xi}\right\|_{1,2}\,,
\end{equation}
where $R_{\varepsilon,t,\xi}:=W_{\varepsilon,t,\xi}-L_g^{-1}\(f_\varepsilon\(W_{\varepsilon,t,\xi}\)\)$. Moreover, $\phi_{\varepsilon,t,\xi}$ is continuously differentiable with respect to $t$ and $\xi$.
\end{proposition}

In Proposition~\ref{Pr2} below, we give a crucial estimate for $\left\|R_{\varepsilon,t,\xi}\right\|_{1,2}$. The proof of Proposition~\ref{Pr2} is presented in Section~\ref{Sec3}.

\begin{proposition}\label{Pr2}
Given two positive real numbers $a<b$, there exists a positive constant $C=C\(a,b,n,M,g,h\)$ such that for $\varepsilon>0$ small, for any $t\in\[a,b\]$ and $\xi\in M$, there holds
\begin{equation}\label{Pr2Eq1}
\left\|R_{\varepsilon,t,\xi}\right\|_{1,2}\le C\left\{\begin{aligned}&\varepsilon e^{-\frac{t}{\varepsilon}}&&\text{if }n=4\\&\varepsilon^{\frac{5}{2}}&&\text{if }n=5\\
&\varepsilon^2 \left|\ln \varepsilon\right|^{\frac{2}{3}}&&\text{if }n=6\text{ and }\(M,g\)\text{ l.c.f.}\\
&\varepsilon^{\frac{n+2}{2(n-4)}}&&\text{if }n\ge7\text{ and }\(M,g\)\text{ l.c.f.}\\
&\varepsilon^2\left|\ln \varepsilon\right|^{-\frac{1}{3}}&&\text{if }n=6\text{ and }\(M,g\)\text{ non-l.c.f.}\\
&\varepsilon^2&&\text{if }n\ge7\text{ and }\(M,g\)\text{ non-l.c.f.}
\end{aligned}\right.
\end{equation}
where $R_{\varepsilon,t,\xi}$ is as in Proposition~\ref{Pr1}.
\end{proposition}

For $\varepsilon>0$ small, we define $J_\varepsilon:H^2_1\(M\)\to\mathbb{R}$ by
\begin{equation}\label{Eq22}
J_\varepsilon\(u\):=\frac{1}{2}\int_M\left|\nabla u\right|^2_gdv_g+\frac{1}{2}\int_M\(\alpha_n\Scal_g+\varepsilon h\)u^2dv_g-\frac{1}{2^*}\int_M\left|u\right|^{2^*}dv_g.
\end{equation}
For any $t>0$ and $\xi\in M$ we define
\begin{equation}\label{Eq23}
\mathcal{J}_\varepsilon\(t,\xi\):=J_\varepsilon\(W_{\varepsilon,t,\xi}+\phi_{\varepsilon,t,\xi}\),
\end{equation}
where $\phi_{\varepsilon,t,\xi}$ is given by Proposition~\ref{Pr1}. We can solve equation \eqref{Eq19} by searching critical points of $\mathcal{J}_\varepsilon$, as it follows from \eqref{Eq27}--\eqref{Eq28} and \eqref{Pr3Eq9}--\eqref{Pr3Eq10}. To this aim, it becomes crucial to have the asymptotic expansion of $\mathcal{J}_\varepsilon$ given by Proposition~\ref{Pr3} below. The proof of Proposition~\ref{Pr3} strongly relies on Propositions~\ref{Pr1} and~\ref{Pr2}, and is presented in section~\ref{Sec4}.

\medskip
We define the ``reduced energy'' $\widetilde E:(0,\infty)\times M \to \mathbb{R}$ as follows:
\begin{equation} \label{tildeE} \widetilde E(t,\xi)=\left\{ \begin{array}{ll}e^{-\frac{2t}{\varepsilon}}(c_2th(\xi)-c_3A_\xi)&\hbox{if }n=4\\
c_2 t^2 h(\xi)-c_3 t^{n-2} A_\xi &\hbox{if }n=5 \hbox{ or }(n\geq 6 \text{ and }\(M,g\)\text{ l.c.f.})\\
c_2 t^2 h(\xi)-c_3 t^4 \left|\Weyl_g\(\xi\)\right|_g^2 &\hbox{if }n\geq 6\text{ and }\(M,g\)\text{ non-l.c.f.},
\end{array} \right.
\end{equation}
where $c_2,c_3>0$, $\Weyl_g$ is the Weyl curvature tensor of the manifold and $A_\xi$ is as in \eqref{Eq3}.

\begin{proposition}\label{Pr3}
Let $p\in\left\{0,1\right\}$ and assume that $h\in C^{0,\alpha}\(M\)$, $\alpha\in\(0,1\)$, in case $p=0$ and $h\in C^1\(M\)$ in case $p=1$. Then there holds
$$\mathcal{J}_\varepsilon\(t,\xi\)=c_1+\left\{\begin{array}{ll}
\widetilde E(t,\xi)+\operatorname{o}\big(e^{-\frac{2\varepsilon}{t}}\big) &\text{if }n=4\\
\varepsilon^{\frac{n-2}{n-4}}\widetilde E(t,\xi)+\operatorname{o}\big(\varepsilon^{\frac{n-2}{n-4}}\big) &\text{if }n=5\,\text{or\,(}n\ge6\,\text{and}\,\(M,g\)\text{ l.c.f.)}\\
\varepsilon^2 \big(\ln \frac{1}{\varepsilon} \big)^{-1}  \widetilde E(t,\xi)+\operatorname{o} \big(\varepsilon^2 \(\ln \frac{1}{\varepsilon} \)^{-1} \big) &\text{if }n=6\text{ and }\(M,g\)\text{non-l.c.f.}
\\
\varepsilon^2 \widetilde E(t,\xi)+\operatorname{o}\(\varepsilon^2\) &\text{if }n\ge7 \text{ and }\(M,g\)\text{non-l.c.f.}
\end{array}\right.$$
as $\varepsilon \to 0$, $C^p$--uniformly with respect to $\xi\in M$ and $t$ in compact subsets of $(0,\infty)$, where $\widetilde E$ is given by \eqref{tildeE} and $c_1,c_2,c_3>0$ depend only on $n$.
\end{proposition}

We are now ready to state the following general result.

\begin{theorem}\label{Th}
Let $\(M,g\)\not=\(\mathbb{S}^n,g_0 \)$ be a smooth compact Riemannian manifold with $n\ge4$ and $\mu_g(M)>0$. Let $p\in\left\{0,1\right\}$ and assume that $h\in C^{0,\alpha}\(M\)$, $\alpha\in\(0,1\)$, in case $p=0$ and $h\in C^1\(M\)$ in case $p=1$. Assume that there exists a $C^p$--{\it stable critical set} $\mathcal{\widetilde D} \subset (0,\infty) \times M$, $p=0,1$, of the function $\widetilde E$. Then for $\varepsilon>0$ small, equation \eqref{Eq1} has a solution $u_\varepsilon \in C^{2,\alpha}\(M\)$, $\forall\: \alpha \in (0,1)$ if $p=1$, such that the family $\(u_\varepsilon\)_\varepsilon$ blows up, up to a sub-sequence, at some $\xi_0 \in \pi(\mathcal{\widetilde D})$ as $\varepsilon\to+\infty$. Here, $\pi: (0,\infty)\times M \to M$ is  the projection operator onto the second component. Moreover, when $p=1$, $\xi_0$ is a critical point of $E$ with $h(\xi_0)>0$, where $E$ is defined in \eqref{ThEq1}.
\end{theorem}

The notion of stability we are using is essentially taken from Li~\cite{Li0}. We say that a compact set $\mathcal{\widetilde D} \subset (0,\infty)\times M$ is a $C^p$--{\it stable critical set}, $p\in \mathbb{N}$, if for any compact neighborhood $\widetilde U$ of $\mathcal{\widetilde D}$ in $(0,\infty) \times M$, there exists $\delta>0$ such that, if $\widetilde J\in C^1(\widetilde U)$ and $\|\mathcal{\widetilde J}-\widetilde E\|_{C^p(\widetilde U)}\leq \delta$, then $\mathcal{\widetilde J}$ has at least one critical point in $\widetilde U$. Since $\widetilde E$ depends on $\varepsilon$ when $n=4$, the above assumption $\|\mathcal{\widetilde J}-\widetilde E\|_{C^p(\widetilde U)}\leq \delta$ needs to be interpreted in this case as: $|\mathcal{\widetilde J}(t,\xi)-\widetilde E(t,\xi)| \leq \delta e^{-\frac{2t}{\varepsilon}}$ for all $(t,\xi) \in \widetilde U$, when $p=0$;
$|\mathcal{\widetilde J}(t,\xi)-\widetilde E(t,\xi)|+\varepsilon |\partial_t(\mathcal{\widetilde J}(t,\xi)-\widetilde E(t,\xi))| +|\nabla_\xi (\mathcal{\widetilde J}(t,\xi)-\widetilde E(t,\xi))|_g \leq \delta e^{-\frac{2t}{\varepsilon}}$ for all $(t,\xi) \in \widetilde U$, when $p=1$. Observe also that we do not require the set $\mathcal{\widetilde D}$ to be composed of critical points of $\widetilde E$, as it would be intuitively reasonable. Indeed, we want   to include the case where $\mathcal{\widetilde D}$ is given by almost critical points of $\widetilde E$, as it arises for $n=4$ in the proof of Theorem~\ref{ThTh}. However, since $\widetilde U$ can shrink onto $\mathcal{\widetilde D}$, by compactness of $\mathcal{\widetilde D}$, we have that $\mathcal{\widetilde D}$ contains at least one critical point of $\widetilde E$.

\smallskip
\proof[Proof of Theorem~\ref{Th}]
Let $\widetilde U$ be a compact neighborhood of the $C^p$--{\it stable critical set} $\mathcal{\widetilde D}$ in $(0,\infty) \times M$, $p=0,1$. Given any $\delta>0$, by Proposition~\ref{Pr3}, we have that $\mathcal{\widetilde J}_\varepsilon:=\mu_\varepsilon^{-1} (\mathcal{J}_\varepsilon-c_1)$ satisfies $\|\mathcal{\widetilde J}_\varepsilon-\widetilde E\|_{C^p(\widetilde U)}\leq \delta$ for $\varepsilon$ sufficiently small, where $\mu_\varepsilon=1$ if $n=4$, $\mu_\varepsilon:=\varepsilon^{\frac{n-2}{n-4}}$ if $n=5$ or $n\ge6$ and $\(M,g\)$ l.c.f., $\mu_\varepsilon:=\varepsilon^2 (\ln \frac{1}{\varepsilon})^{-1}$ if $n=6$ and $\(M,g\)$ non-l.c.f., $\mu_\varepsilon:=\varepsilon^2$ if $n\ge7$ and $\(M,g\)$ non-l.c.f. By definition of a $C^p$--{\it stable critical set}, $p=0,1$, it follows that $\mathcal{J}_\varepsilon$ has a critical point $\(t_\varepsilon,\xi_\varepsilon\)\in \widetilde U$ for $\varepsilon$ small. Up to a subsequence and taking $\widetilde U$ smaller and smaller, we can assume that $\(t_\varepsilon,\xi_\varepsilon\) \to \(t_0,\xi_0\)$ as $\varepsilon \to 0$ with $\xi_0 \in \pi(\mathcal{\widetilde D})$. As already observed, we get that $u_\varepsilon=W_{\varepsilon,t_\varepsilon,\xi_\varepsilon}+\phi_{\varepsilon,t_\varepsilon,\xi_\varepsilon}$ is a critical point of $J_\varepsilon$, and thus, by elliptic regularity, a classical solution of \eqref{Eq1}. Since $\left\|\phi_{\varepsilon,t_\varepsilon,\xi_\varepsilon} \right\|_{1,2}\to0$, by definition of $W_{\varepsilon,t_\varepsilon,\xi_\varepsilon}$, it is easily seen that $u_\varepsilon>0$ and $u_\varepsilon^{2^*} \rightharpoonup K_n^{-n}\delta_{\xi_0}$  in the measures sense as $\varepsilon \to 0$ (see for instance Rey~\cite{Rey}), where $K_n$ is given by \eqref{Eq24} and $\delta_{\xi_0}$ denotes the Dirac mass measure at $\xi_0$. From very basic facts concerning the asymptotic analysis of solutions of Yamabe-type equations (see for instance Druet--Hebey~\cites{DruHeb2} and Druet--Hebey--Robert~\cite{DruHebRob}) we get that the family $(u_\varepsilon)_\varepsilon$ of solutions to \eqref{Eq1} blows up at the point $\xi_0$ as $\varepsilon \to 0$. Finally, when $p=1$, we can pass to the limit into the equations $\partial_t \mathcal{J}_\varepsilon \(t_\varepsilon,\xi_\varepsilon\)=0$ and $\nabla_\xi\mathcal{J}_\varepsilon \(t_\varepsilon,\xi_\varepsilon\)=0$ as $\varepsilon \to 0$ to get that $h(\xi_0)>0$ and $\nabla_\xi E(\xi_0)=0$ in view of $t_0>0$, where $E$ is given by \eqref{ThEq1}. This ends the proof of Theorem~\ref{Th}.
\endproof

We now apply Theorem~\ref{Th} to deduce Theorems~\ref{ThTh} and~\ref{ThThTh}. To this aim, given $\xi \in M$ with $h(\xi)>0$, define $t(\xi)$ as
\begin{equation}\label{txi}
t(\xi):=\left\{ \begin{array}{ll} \frac{c_3 A_\xi}{c_2 h(\xi)}+\frac{\varepsilon}{2} &\hbox{if }n=4\\
\left(\frac{2c_2 h(\xi)}{(n-2)c_3 A_\xi}\right)^{\frac{1}{n-4}}&\hbox{if }n=5 \hbox{ or }(n\geq 6 \text{ and }\(M,g\)\text{ l.c.f.})\\
\left(\frac{c_2 h(\xi)}{2 c_3 \left|\Weyl_g\(\xi\)\right|_g^2}\right)^{\frac{1}{2}}    &\hbox{if }n\geq 6\text{ and }\(M,g\)\text{ non-l.c.f.}
\end{array} \right.
\end{equation}
with the convention that $1/0=+\infty$. One easily checks that every $t\(\xi\)<+\infty$ is a global maximum point of $\widetilde E$ in $t$. In the proofs of Theorems~\ref{ThTh} and~\ref{ThThTh} below, we show that the $C^p$--{\it stable critical set} $\mathcal{\widetilde D}$ in Theorem~\ref{Th} can be constructed as $\mathcal{\widetilde D}:=\{(t(\xi),\xi):\,\xi \in \mathcal{D} \},$ where $\mathcal{D}$ is a $C^p$--{\it stable critical set} of $\xi\mapsto\widetilde E(t(\xi),\xi)$ with $h>0$ in $\mathcal{D}$. Since we have that
\begin{equation} \label{Etxixi}
\widetilde E(t(\xi),\xi)=\left\{ \begin{array}{ll} \frac{\varepsilon}{2} c_2 e^{-1} h(\xi) e^{-\frac{2 c_3}{\varepsilon c_2}\cdot\frac{A_\xi}{h(\xi)} } &\hbox{if }n=4\\
\frac{2^{\frac{2}{n-4}} (n-4) c_2^{\frac{n-2}{n-4}} }{(n-2)^{\frac{n-2}{n-4}} c_3^{\frac{2}{n-4}}}\cdot\frac{h(\xi)^{\frac{n-2}{n-4}}}{A_\xi^{\frac{2}{n-4}}}&\hbox{if }n=5 \hbox{ or }(n\geq 6 \text{ and }\(M,g\)\text{ l.c.f.})\\
\frac{c_2^2}{4c_3}\cdot\frac{h(\xi)^2}{\left|\Weyl_g\(\xi\)\right|_g^2}  &\hbox{if }n\geq 6\text{ and }\(M,g\)\text{ non-l.c.f.},
\end{array} \right.
\end{equation}
the role played by the map $E$ defined in \eqref{ThEq1} becomes clear. We prove Theorem~\ref{ThTh} as follows.

\smallskip
\proof[Proof of Theorem~\ref{ThTh}]
Notice that $\sup_M E>0$ in view of $\max_M h>0$. Letting $\(\xi_k\)_{k\in\mathbb{N}}$ be a maximizing sequence for $E$, by compactness of $M$, we can assume that $\xi_k \to \xi$ as $k \to +\infty$ with $h(\xi_k)>0$ for all $k$. Since {$\min\{|\Weyl_g(\xi)|_g\,:\,h(\xi)>0\}>0$} in the non-l.c.f.~case with $n \geq 6$, we have that $\(E(\xi_k)\)_k$ is bounded. Since then $\sup_M E<+\infty$, it follows that $\sup_M E$ is achieved, and it makes sense to define $\mathcal{D}:=\{ \xi \in M:\ E(\xi)=\sup_M E\}$ (possibly coinciding with the whole $M$ if $E$ is a constant function) as the maximal set of $E$ in $M$. Correspondingly, define $\mathcal{\widetilde D}:=\{(t(\xi),\xi):\ \xi \in \mathcal{D}\}$, where $t(\xi)$ is given by \eqref{txi} and is well defined for all $\xi \in \mathcal{D}$ in view of $h>0$ in $\mathcal{D}$. Moreover, $t(\xi)$ is clearly bounded away from zero on $\mathcal{D}$, and thus the set $\mathcal{\widetilde D}$ is a compact set in $(0,\infty) \times M$.
To show that $\mathcal{\widetilde D}$ is a $C^0$--{\it stable critical set}, let $\widetilde U$ be a compact neighborhood of $\mathcal{\widetilde D}$ in $(0,\infty) \times M$. Taking $\widetilde U$ smaller if necessary, we can assume that $\widetilde U=\{(t,\xi): \
t \in [t(\xi)-\eta,t(\xi)+\eta],\ \xi \in U \}$, where $\eta>0$ is small and $U$ is a closed neighborhood of $\mathcal{D}$ in $M$ so that $h>0$ in $U$. There hold
\begin{itemize}
\item for $n\geq 5$, by the definition of $t(\xi)$ and the simple relation between $\widetilde E(t(\xi),\xi)$ and $E(\xi)$ (see \eqref{ThEq1} and \eqref{Etxixi}), we have that $\widetilde E(t(\xi)\pm \eta,\xi)<\widetilde E(t(\xi),\xi)\leq \sup_{\mathcal{\widetilde D}} \widetilde E$ if $\xi \in U$, and $\widetilde E(t,\xi)\leq \widetilde E(t(\xi),\xi)<\sup_{\mathcal{\widetilde D}} \widetilde E$ if $t \in [t(\xi)-\eta,t(\xi)+\eta]$ and $\xi \in \partial U$;
\item for $n=4$, by the definition of $\widetilde E$ and $t(\xi)$, we have that for any $\xi \in U$, there holds
$$\widetilde E(t(\xi)\pm \eta,\xi)=\widetilde E(t(\xi),\xi) e^{\mp \frac{2\eta }{\varepsilon}}(1 \pm \frac{2\eta}{\varepsilon}) < \sup_{\mathcal{\widetilde D}} \widetilde E$$
when $\varepsilon$ is small, in view of $e^{-\frac{2\eta }{\varepsilon}}(1 + \frac{2\eta}{\varepsilon})\to 0$ and $e^{\frac{2\eta }{\varepsilon}}(1 -\frac{2\eta}{\varepsilon}) \to -\infty$ as $\varepsilon \to 0$, and for any $t \in [t(\xi)-\eta,t(\xi)+\eta]$ and $\xi \in \partial U$, there holds
$$\widetilde E(t,\xi)\leq \widetilde E(t(\xi),\xi)< \sup_{\mathcal{\widetilde D}} \widetilde E$$
when $\varepsilon$ is small, in view of $\displaystyle\sup_{\partial U}\, e^{-\frac{2 c_3}{\varepsilon c_2} \frac{A_\xi}{h(\xi)}}=\operatorname{o}\Big(e^{-\frac{2 c_3}{\varepsilon c_2 \max_M E}}\Big)$ as $\varepsilon \to 0$.
\end{itemize}
In conclusion, by compactness of
$$\partial \widetilde{U}=\{(t,\xi):\ t\in\{t(\xi)-\eta,t(\xi)+\eta\},\ \xi \in U\}\cup\{(t,\xi):\ t \in [t(\xi)-\eta,t(\xi)+\eta],\ \xi \in \partial U\},$$
we get that $\sup_{\partial \widetilde U}\widetilde E<\sup_{\mathcal{\widetilde D}} \widetilde E$. It follows that if $\|\mathcal{\widetilde J}-\widetilde E\|_{C^0(\widetilde U)}\leq \delta$, with $\delta <\frac{1}{2}[\sup_{\mathcal{\widetilde D}} \widetilde E-\sup_{\partial \widetilde U}\widetilde E]$, then we get that
$$\sup_{\partial \widetilde U}\mathcal{\widetilde J}<\sup_{\mathcal{\widetilde D}}\mathcal{\widetilde J}.$$
Then $\mathcal{\widetilde J}$ achieves its maximum value in $\widetilde U$ at some interior point of $\widetilde U$, which is a critical point of $\mathcal{\widetilde J}$. It follows that $\mathcal{\widetilde D}$ is a $C^0$--{\it stable critical set} of $\widetilde E$ as desired. By Theorem~\ref{Th}, we then get that for $\varepsilon>0$ small, equation \eqref{Eq1} has a solution $u_\varepsilon \in C^{2,\alpha}\(M\)$ such that the family $\(u_\varepsilon\)_\varepsilon$ blows up, up to a sub-sequence, at some $\xi_0 \in \pi(\mathcal{\widetilde D})$ as $\varepsilon\to+\infty$. Moreover, by definition of $\mathcal{\widetilde D}$ , we get that $\widetilde E\(\xi_0\)=\max_M\widetilde E$. This ends the proof of Theorem~\ref{ThTh}.
\endproof

Since $t(\xi)$ is a maximum point of $\widetilde E$ in $t$, minimum points or saddle points of $E$   provide critical points of  $\widetilde E$ which in general are not $C^0$--{\it stable critical points} of $\widetilde E$. To cover these cases, we need to use Theorem~\ref{Th} with $p=1$. We assume that $h\in C^1\(M\)$. We can still define $C^1$--{\it stable critical sets} of $E$ as in the case of $(0,\infty)\times M$, but in general they don't give rise to $C^1$--{\it stable critical sets} of $\widetilde E$. In $M$, we restrict the notion of $C^1$--stability to isolated critical points of $E$ with non-trivial local degree, which is still sufficiently general to include non-degenerate critical points of $E$. We prove Theorem~\ref{ThThTh} as follows.

\smallskip
\proof[Proof of Theorem~\ref{ThThTh}]
We only need to show that the set $\mathcal{\widetilde D}:=\{ (t(\xi_0),\xi_0) \}$ is a $C^1$--{\it stable critical set} of $\widetilde E$, where $t(\xi_0)>0$ is well defined in view of $h(\xi_0)>0$ and $\Weyl_g\(\xi_0\)\not=0$ when $n\geq 6$ and $\(M,g\)$ non-l.c.f. To this aim, let $\widetilde U$ be any compact neighborhood of $\mathcal{\widetilde D}$ in $(0,\infty) \times M$ and $\delta>0$ be any given small number. Taking $\widetilde U$ smaller if necessary, we can assume that $\widetilde U=I \times B_{\xi_0}\(r_0\)$, where $I$ is a closed interval in $(0,\infty)$ containing $t(\xi_0)$ in its interior and $B_{\xi_0}\(r_0\)$ is the geodesic ball of center $\xi_0$ and radius $r_0$ with respect to $g_{\xi_0}$, with $r_0<i_{g_{\xi_0}}$, the injectivity radius of the manifold $\(M,g_{\xi_0}\)$. For $n\geq 5$, the assumption $\|\mathcal{\widetilde J}-\widetilde E\|_{C^1(\widetilde U)}\leq \delta$ gives that
\begin{equation}
\big|\nabla \mathcal{\widetilde J}\(t,\exp_{\xi_0}\eta\)  -\nabla \widetilde E\(t,\exp_{\xi_0}\eta\)\big|\leq \delta
\label{ThEq2}
\end{equation}
as $\varepsilon\to0$, uniformly with respect to $\eta\in B_0\(r_0\)$ and $t \in I$, where $\nabla:=\(d/dt,\nabla_{\hspace{-2pt}\eta}\)$. Since $\xi_0$ is a $C^1$--{\it stable critical point} of $E$ with $h\(\xi_0\)>0$, where $E$ is given by \eqref{ThEq1}, we find that $\nabla \widetilde E\(t,\exp_{\xi_0}\eta\)$ has an isolated zero at $\(t(\xi_0),0\)$. Since $t(\xi_0)$ is a non-degenerate critical point of $t \mapsto \widetilde E\(t,\exp_{\xi_0}\eta\)$, taking $\widetilde U$ smaller if necessary, we also have that $\deg\big(\nabla \widetilde E, \widetilde U,0\big)\ne0$  (see for instance~\cite{FonGan}). It follows from \eqref{ThEq2} that if $\delta$ is small enough, then $\mathcal{\widetilde J}$ has at least one critical point $\(t,\xi\)\in \widetilde U$. When $n=4$, the assumption $\|\mathcal{\widetilde J}-\widetilde E\|_{C^1(\widetilde U)}\leq \delta$ gives (by definition) that
\begin{equation}
\varepsilon e^{\frac{2t}{\varepsilon}}\big|\partial_t \mathcal{\widetilde J} \(t,\exp_{\xi_0}\eta\)- \partial_t \widetilde E\(t,\exp_{\xi_0}\eta\)\big|+ e^{\frac{2t}{\varepsilon}}\big|\nabla_\eta \mathcal{\widetilde J} \(t,\exp_{\xi_0}\eta\)- \nabla_\eta \widetilde E\(t,\exp_{\xi_0}\eta\)\big| \leq C_0 \delta \label{ThEq3}
\end{equation}
as $\varepsilon\to0$, uniformly with respect to $\eta\in B_0\(r_0\)$ and $t \in I$, for some $C_0>0$. Letting
$$\Psi(t,\eta)=\left(-2 c_2 t h(\exp_{\xi_0}\eta)+2c_3 A_{\exp_{\xi_0}\eta},c_2 t \nabla_\eta h(\exp_{\xi_0}\eta)-c_3 \nabla_\eta A_{\exp_{\xi_0}\eta}\right)$$
and $\widetilde \Psi(t,\xi)=\Psi(t,\exp_{\xi_0}^{-1}\xi)$, by \eqref{ThEq3} we deduce that
\begin{equation}
\Big|\left(\varepsilon e^{\frac{2t}{\varepsilon}} \partial_t \mathcal{\widetilde J} \(t,\exp_{\xi_0}\eta\), e^{\frac{2t}{\varepsilon}} \nabla_\eta \mathcal{\widetilde J} \(t,\exp_{\xi_0}\eta\)\right)-\Psi (t,\eta)\Big| \leq C_0 \delta +\operatorname{O}(\varepsilon)
\label{ThEq4}
\end{equation}
as $\varepsilon\to0$, uniformly with respect to $\eta\in B_0\(r_0\)$ and $t \in I$. Arguing as above, the map $\Psi$ has an isolated zero at $\big(\frac{c_3 A_{\xi_0}}{c_2 h(\xi_0)},0\big)$ with $\deg\big(\widetilde \Psi, \widetilde U,0\big)\ne0$, and then by  \eqref{ThEq4}, it follows that if $\delta$ is small enough, then $\mathcal{\widetilde J}$ has at least one critical point $\(t,\xi\)\in \widetilde U$.
\endproof

\section{The error estimate}\label{Sec3}

This section is devoted to the estimate of $R_{\varepsilon,t,\xi}$.

\smallskip
\proof[Proof of Proposition~\ref{Pr2}]
All our estimates in this proof are uniform with respect to $t\in\[a,b\]$, $\xi\in M$ and $\varepsilon\in\(0,\varepsilon_0\)$, for some fixed $\varepsilon_0>0$. We let $r_0$ be as in Section~\ref{Sec2}. For any $\phi\in H^2_1\(M\)$, an integration by parts gives that
\begin{multline}\label{Pr2Eq2}
\<L_g^{-1}\(f_\varepsilon\(W_{\varepsilon,t,\xi}\)\)-W_{\varepsilon,t,\xi},\phi\>_{L_g}=\int_M\(f_\varepsilon\(W_{\varepsilon,t,\xi}\)-L_gW_{\varepsilon,t,\xi}\)\phi dv_g\\
-\int_{\partial B_\xi\(r_0\)}\(\partial_{\nu_{\inward}}W_{\varepsilon,t,\xi}+\partial_{\nu_{\outward}}W_{\varepsilon,t,\xi}\)\phi d\sigma_g\,,
\end{multline}
where $\partial B_\xi\(r_0\)$ is the boundary of the geodesic ball with respect to $g_\xi$ of center $\xi$ and radius $r_0$, $\partial_{\nu_{\inward}}$ and $\partial_{\nu_{\outward}}$ are the derivatives with respect to the respective inward and outward, unit, normal vectors to $\partial B_\xi\(r_0\)$, and $d\sigma_g$ is the volume element on $\partial B_\xi\(r_0\)$. By Sobolev's and trace's embeddings, it follows from \eqref{Pr2Eq2} that
\begin{multline}\label{Pr2Eq3}
\left\|L_g^{-1}\(f_\varepsilon\(W_{\varepsilon,t,\xi}\)\)-W_{\varepsilon,t,\xi}\right\|_{1,2}=\operatorname{O}\Big(\left\|f_\varepsilon\(W_{\varepsilon,t,\xi}\)-L_gW_{\varepsilon,t,\xi}\right\|_{L^{\frac{2n}{n+2}}(M) }\\
+\left\|\partial_{\nu_{\inward}}W_{\varepsilon,t,\xi}+\partial_{\nu_{\outward}}W_{\varepsilon,t,\xi} \right\|_{L^{\frac{2\(n-1\)}{n}}\(\partial B_\xi\(r_0\)\)}\Big).
\end{multline}
Regarding the second term in the right hand side of \eqref{Pr2Eq3}, on $\partial B_\xi\(r_0\)$, we find that
\begin{align}
\partial_{\nu_{\inward}}W_{\varepsilon,t,\xi}+\partial_{\nu_{\outward}}W_{\varepsilon,t,\xi}&=\beta_n \delta_\varepsilon\(t\)^{\frac{2-n}{2}} G_g\(\cdot,\xi\) \frac{d}{dr}\(r^{n-2}U(\delta_\varepsilon(t)^{-1}r)\)\Big|_{r=r_0}\nonumber\\
&=n^{\frac{n-2}{4}}\(n-2\)^{\frac{n+2}{4}}\beta_nG_g\(\cdot,\xi\)\frac{\delta_\varepsilon\(t\)^{\frac{n+2}{2}}r_0^{n-3}}{\(\delta_\varepsilon\(t\)^2+r_0^2\)^{\frac{n}{2}}}=\operatorname{O}\(\delta_\varepsilon(t)^{\frac{n+2}{2}}\).\label{Pr2Eq4}
\end{align}
Regarding the first term in the right hand side of \eqref{Pr2Eq3}, we observe that in $M \setminus B_\xi\(r_0\)$, there holds
\begin{equation}\label{Pr2Eq5bis}
f_\varepsilon\(W_{\varepsilon,t,\xi}\)-L_gW_{\varepsilon,t,\xi}=f_\varepsilon\(W_{\varepsilon,t,\xi}\)=\operatorname{O}\(\delta_\varepsilon(t)^{\frac{n+2}{2}}+\varepsilon \delta_\varepsilon(t)^{\frac{n-2}{2}}\)
\end{equation}
in view of
$$\widehat{W}_{\varepsilon,t,\xi}=\beta_n \delta_\varepsilon(t)^{\frac{2-n}{2}} r_0^{n-2}U\(\delta_\varepsilon(t)^{-1}r_0\)=\operatorname{O}(\delta_\varepsilon(t)^{\frac{n-2}{2}})$$
and $L_g G_g(\cdot,\xi)=0$. By conformal covariance \eqref{confinv} of $L_g$ and by \eqref{Eq7}, in $ B_\xi\(r_0\)$, we can write that
\begin{equation}\label{Pr2Eq5}
f_\varepsilon\(W_{\varepsilon,t,\xi}\)-L_gW_{\varepsilon,t,\xi}=\varLambda_\xi^{2^*-1}\big[\big(G_{g_\xi}\(\cdot,\xi\)\widehat{W}_{\varepsilon,t,\xi}\big)^{2^*-1}-L_{g_\xi}\big(G_{g_\xi}\(\cdot,\xi\)\widehat{W}_{\varepsilon,t,\xi}\big)\big]-\varepsilon hW_{\varepsilon,t,\xi}.
\end{equation}
Since $\widehat{W}_{\varepsilon,t,\xi}\(\xi\)=0$ and $L_{g_\xi}G_{g_\xi}\(\cdot,\xi\)=\delta_\xi$, we get that
\begin{equation}\label{Pr2Eq6}
L_{g_\xi}\big(G_{g_\xi}\(\cdot,\xi\)\widehat{W}_{\varepsilon,t,\xi}\big)=G_{g_\xi}\(\cdot,\xi\)\Delta_{g_\xi}\widehat{W}_{\varepsilon,t,\xi}-2\big<\nabla G_{g_\xi}\(\cdot,\xi\),\nabla\widehat{W}_{\varepsilon,t,\xi}\big>_{g_\xi}\,.
\end{equation}
Since $\widehat{W}_{\varepsilon,t,\xi}\circ\exp_\xi$ is radially symmetrical in $B_0\(r_0\)$, writing $\Delta_{g_\xi}\widehat{W}_{\varepsilon,t,\xi}\(\exp_\xi y\)$ in polar coordinates, by \eqref{Eq6}, we find that
\begin{align}
\Delta_{g_\xi}\widehat{W}_{\varepsilon,t,\xi}\(\exp_\xi y\)&=\Delta_{\Eucl}\big(\widehat{W}_{\varepsilon,t,\xi}\circ\exp_\xi\big)\(y\)+\operatorname{O}\(\left|y\right|^{N-1}\big|\nabla\big(\widehat{W}_{\varepsilon,t,\xi}\circ\exp_\xi\big)\(y\)\big|\)\nonumber\\
&=\beta_n\delta_\varepsilon\(t\)^{-\frac{n+2}{2}}\left|y\right|^{n-2}U\( \frac{y}{\delta_\varepsilon(t)}\)^{2^*-1}-2n^{\frac{n-2}{4}}\(n-2\)^{\frac{n+6}{4}}\beta_n\frac{\delta_\varepsilon\(t\)^{\frac{n+2}{2}}\left|y\right|^{n-4}}{\(\delta_\varepsilon\(t\)^2+\left|y\right|^2\)^{\frac{n}{2}}}\nonumber\\
&\quad+\operatorname{O}\(\frac{\delta_\varepsilon(t)^{\frac{n+2}{2}}\left|y\right|^{N+n-4}}{\(\delta_\varepsilon(t)^2+\left|y\right|^2\)^{\frac{n}{2}}}\)\label{Pr2Eq7}
\end{align}
uniformly with respect to $y\in B_0\(r_0\)$, in view of $\Delta_{\Eucl}U=U^{2^*-1}$ in $\mathbb{R}^n$. Moreover, since $\widehat{W}_{\varepsilon,t,\xi}\circ\exp_\xi$ is radially symmetrical, we get that
\begin{multline} \label{added}
\big<\nabla G_{g_\xi}\(\exp_\xi y,\xi\),\nabla\widehat{W}_{\varepsilon,t,\xi}\(\exp_\xi y\)\big>_{g_\xi}=
\partial_r  \[G_{g_\xi}\(\exp_\xi y,\xi\) \] \partial_r \big[\widehat{W}_{\varepsilon,t,\xi} \circ \exp_\xi(y)\big]\\
=  \partial_r  \[G_{g_\xi}\(\exp_\xi y,\xi\) \]  n^{\frac{n-2}{4}}\(n-2\)^{\frac{n+2}{4}}\beta_n\frac{\delta_\varepsilon\(t\)^{\frac{n+2}{2}}\left|y\right|^{n-3}}{\(\delta_\varepsilon\(t\)^2+\left|y\right|^2\)^{\frac{n}{2}}}
\end{multline}
for all $y\in B_0\(r_0\)$. Inserting \eqref{Pr2Eq7} and \eqref{added} into \eqref{Pr2Eq6}, we get that
\begin{multline} \label{added1}
L_{g_\xi}\big(G_{g_\xi}\(\cdot,\xi\)\widehat{W}_{\varepsilon,t,\xi}\big)(\exp_\xi y)=
 \beta_n G_{g_\xi}\(\exp_\xi y,\xi\) \left|y\right|^{n-2} \delta_\varepsilon\(t\)^{-\frac{n+2}{2}} U\(\delta_\varepsilon(t)^{-1}y\)^{2^*-1}\\
-2 n^{\frac{n-2}{4}}\(n-2\)^{\frac{n+2}{4}} \beta_n \((n-2) G_{g_\xi}\(\exp_\xi y,\xi\)
+|y| \partial_r\[G_{g_\xi}\(\exp_\xi y,\xi\)\]\) \frac{\delta_\varepsilon\(t\)^{\frac{n+2}{2}}\left|y\right|^{n-4}}{\(\delta_\varepsilon\(t\)^2+\left|y\right|^2\)^{\frac{n}{2}}}\\
+\operatorname{O}\(\frac{\delta_\varepsilon(t)^{\frac{n+2}{2}}G_{g_\xi}\(\exp_\xi y,\xi\) \left|y\right|^{N+n-4}}{\(\delta_\varepsilon(t)^2+\left|y\right|^2\)^{\frac{n}{2}}}\)
\end{multline}
in $B_0\(r_0\)$. Using Lemma~\ref{Lem4}, by \eqref{added1} we find that
\begin{multline}\label{Pr2Eq8}
\big(G_{g_\xi}\(\exp_\xi y,\xi\)\widehat{W}_{\varepsilon,t,\xi}\(\exp_\xi y\)\big)^{2^*-1}-L_{g_\xi}\big(G_{g_\xi}\(\cdot,\xi\)\widehat{W}_{\varepsilon,t,\xi}\big)\(\exp_\xi y\)\\
=\frac{\delta_\varepsilon(t)^{\frac{n+2}{2}}}{\(\delta_\varepsilon(t)^2+\left|y\right|^2\)^{\frac{n}{2}}}\times\left\{\begin{aligned}
&\operatorname{O}\(|y|^{n-4}\)&&\text{if }n=4,5\text{ or }\(M,g\)\text{ l.c.f.}\\
&\operatorname{O}\(|y|^2\ln|y|\) &&\text{if }n=6\text{ and }\(M,g\)\text{ non-l.c.f.}\\
&\operatorname{O}\(|y|^2\)&&\text{if }n\ge7\text{ and }\(M,g\)\text{ non-l.c.f.}
\end{aligned}\right.
\end{multline}
in $B_0\(r_0\)$. It follows from \eqref{Pr2Eq8} that
\begin{align}
&\int_{B_\xi\(r_0\)} \varLambda_\xi^{\frac{2n}{n-2}} \left|\big(G_{g_\xi}\(\cdot,\xi\)\widehat{W}_{\varepsilon,t,\xi}\big)^{2^*-1}-L_{g_\xi}\big(G_{g_\xi}\(\cdot,\xi\)\widehat{W}_{\varepsilon,t,\xi}\big)\right|^{\frac{2n}{n+2}}dv_g\nonumber\\
&\qquad=\int_0^{r_0}\frac{\delta_\varepsilon(t)^n dr}{(\delta_\varepsilon(t)^2+r^2)^{\frac{n^2}{n+2}}}
\times \left\{\begin{aligned}
&\operatorname{O}\(r^{\frac{2n^2}{n+2}-1+\frac{n(n-6)}{n+2}}\)&&\text{if }n=4,5\text{ or }\(M,g\)\text{ l.c.f.}\\
&\operatorname{O}\(r^8\left|\ln r\right|\)&&\text{if }n=6\text{ and }\(M,g\)\text{ non-l.c.f.}\\
&\operatorname{O}\(r^{\frac{2n^2}{n+2}-1-\frac{n(n-6)}{n+2}}\)&&\text{if }n\ge7\text{ and }\(M,g\)\text{ non-l.c.f.}
\end{aligned}\right.\nonumber\allowdisplaybreaks\\
&\qquad=\left\{\begin{aligned}
&\operatorname{O}\Big(\delta_\varepsilon(t)^{\frac{2n\(n-2\)}{n+2}}\Big)&&\text{if }n=4,5\\
&\operatorname{O}\(\delta_\varepsilon(t)^6\left|\ln\delta_\varepsilon(t)\right|\)&&\text{if }n=6\text{ and }\(M,g\)\text{ l.c.f.}\\
&\operatorname{O}\(\delta_\varepsilon(t)^n\)&&\text{if }n\ge7\text{ and }\(M,g\)\text{ l.c.f.}\\
&\operatorname{O}\Big(\delta_\varepsilon(t)^6\left|\ln\delta_\varepsilon(t) \right|^2 \Big)&&\text{if }n=6\text{ and }\(M,g\)\text{ non-l.c.f.}\\
&\operatorname{O}\Big(\delta_\varepsilon(t)^{\frac{8n}{n+2}}\Big)&&\text{if }n\ge7\text{ and }\(M,g\)\text{ non-l.c.f.}
\end{aligned}\right.\label{Pr2Eq9}
\end{align}
Moreover, by Lemma~\ref{Lem4}, we find that
\begin{equation}\label{Pr2Eq10}
\int_{B_\xi\(r_0\)}\left|hW_{\varepsilon,t,\xi}\right|^{\frac{2n}{n+2}}dv_g=\operatorname{O}\(\int_0^{r_0}\frac{\delta_\varepsilon(t)^{\frac{n\(n-2\)}{n+2}}r^{n-1} dr}{\(\delta_\varepsilon(t)^2+r^2\)^{\frac{n\(n-2\)}{n+2}}}\)=\left\{\begin{aligned}
&\operatorname{O}\(\delta_\varepsilon(t)^{\frac{n\(n-2\)}{n+2}}\)&&\text{if }n=4,5\\
&\operatorname{O}\(\delta_\varepsilon(t)^3\left|\ln\delta_\varepsilon(t) \right|\)&&\text{if }n=6\\
&\operatorname{O}\Big(\delta_\varepsilon(t)^{\frac{4n}{n+2}}\Big)&&\text{if }n\ge7.
\end{aligned}\right.
\end{equation}
By \eqref{Pr2Eq5bis}, \eqref{Pr2Eq5}, \eqref{Pr2Eq9}, and \eqref{Pr2Eq10}, we get that
\begin{align}
&\left\|f_\varepsilon\(W_{\varepsilon,t,\xi}\)-L_gW_{\varepsilon,t,\xi}\right\|_{\frac{2n}{n+2}} \nonumber\\
&\qquad=\left\{\begin{aligned}
&\operatorname{O}\Big(\delta_\varepsilon(t)^{n-2}+\varepsilon \delta_\varepsilon(t)^{\frac{n-2}{2}} \Big)&&\text{if }n=4,5\\
&\operatorname{O}\(\delta_\varepsilon(t)^4 \left|\ln\delta_\varepsilon(t)\right|^{\frac{2}{3}}+\varepsilon \delta_\varepsilon(t)^2 \left|\ln\delta_\varepsilon(t)\right|^{\frac{2}{3}} \)&&\text{if }n=6\text{ and }\(M,g\)\text{ l.c.f.}\\
&\operatorname{O}\(\delta_\varepsilon(t)^{\frac{n+2}{2}}+\varepsilon \delta_\varepsilon(t)^2 \)&&\text{if }n\ge7\text{ and }\(M,g\)\text{ l.c.f.}\\
&\operatorname{O}\Big(\delta_\varepsilon(t)^4 \left|\ln\delta_\varepsilon(t) \right|^{\frac{4}{3}} +\varepsilon \delta_\varepsilon(t)^2 \left|\ln\delta_\varepsilon(t)\right|^{\frac{2}{3}} \Big)&&\text{if }n=6\text{ and }\(M,g\)\text{ non-l.c.f.}\\
&\operatorname{O}\Big(\delta_\varepsilon(t)^4 +\varepsilon \delta_\varepsilon(t)^2 \Big)&&\text{if }n\ge7\text{ and }\(M,g\)\text{ non-l.c.f.}
\end{aligned}\right.\label{Pr2Eq10bis}
\end{align}
Finally, \eqref{Pr2Eq1} follows from \eqref{Eq8}, \eqref{Pr2Eq4}, and \eqref{Pr2Eq10bis} in view of \eqref{Pr2Eq3}. This ends the proof of Proposition~\ref{Pr2}.
\endproof

\section{The reduced energy}\label{Sec4}

In  Lemma~\ref{Lem1} below, we give an asymptotic expansion of $J_\varepsilon\(W_{\varepsilon,t,\xi}\)$ as $\varepsilon\to0$, where $J_\varepsilon$ is given by \eqref{Eq22}. We let $K_n$ be the sharp constant for the embedding of $D^{1,2}\(\mathbb{R}^n\)$ into $L^{2^*}\(\mathbb{R}^n\)$. It has been proved independently by Rodemich~\cite{Rod}, Aubin~\cite{Aub1}, and Talenti~\cite{Tal} that
\begin{equation}\label{Eq24}
K_n=\sqrt{\frac{4}{n\(n-2\)\omega_n^{2/n}}}\,,
\end{equation}
where $\omega_n$ is the volume of the unit $n$--sphere.

\begin{lemma}\label{Lem1}
We let $K_n$ be as in \eqref{Eq24}, $\Weyl_g$ be the Weyl curvature tensor of the manifold, and for any $\xi\in M$, we let $A_\xi$ be as in \eqref{Eq3}. As $\varepsilon\to0$, the following expansions do hold:
\begin{enumerate}
\renewcommand{\labelenumi}{(\roman{enumi})}
\item when $n=4$,
\begin{equation}\label{Lem1Eq1}
J_\varepsilon\(W_{\varepsilon,t,\xi}\)=\frac{1}{4}K_4^{-4}+4 \omega_3  h(\xi) \varepsilon \delta_\varepsilon(t)^2\ln \frac{1}{\delta_\varepsilon(t)}-4  \beta_4^2 A_\xi \delta_\varepsilon(t)^2+\hbox{h.o.t.,}\end{equation}
\item when $n=5$ or ($n\ge6$ and $\(M,g\)$ is l.c.f.),
\begin{equation}\label{Lem1Eq2}
J_\varepsilon\(W_{\varepsilon,t,\xi}\)=\frac{1}{n}K_n^{-n}+[n(n-2)]^{\frac{n-2}{2}}\( \frac{(n-1)\omega_n}{2^{n-1}(n-4)}h(\xi) \varepsilon\delta_\varepsilon(t)^2-   \frac{\beta_n^2}{2} A_\xi\delta_\varepsilon(t)^{n-2}\)+\hbox{h.o.t.,}
\end{equation}
\item when $n=6$ and $\(M,g\)$ is non-l.c.f.,
\begin{equation}\label{Lem1Eq3}
J_\varepsilon\(W_{\varepsilon,t,\xi}\)=\frac{1}{6}K_6^{-6}+45 \omega_6 h\(\xi\) \varepsilon\delta_\varepsilon(t)^2-\frac{3}{4} \omega_6 \left|\Weyl_g\(\xi\)\right|_g^2 \delta_\varepsilon(t)^4 \ln \frac{1}{\delta_\varepsilon(t)}+\hbox{h.o.t.,}
\end{equation}
\item when $n\ge7$ and $\(M,g\)$ is non-l.c.f.,
\begin{multline}\label{Lem1Eq4}
J_\varepsilon\(W_{\varepsilon,t,\xi}\)=\frac{1}{n}K_n^{-n}+\frac{[n(n-2)]^{\frac{n-2}{2}}}{2^{n-1}(n-4)} \omega_n \( (n-1)h(\xi) \varepsilon\delta_\varepsilon(t)^2- \frac{n-2}{48 (n-6) } \left|\Weyl_g\(\xi\)\right|_g^2 \delta_\varepsilon(t)^4\)\\
+\hbox{h.o.t.}
\end{multline}
\end{enumerate}
uniformly with respect to $\xi\in M$ and $t$ in compact subsets of $(0,\infty)$, where $\hbox{h.o.t.}$ stands for a term which is asymptotically smaller than one of the previous terms in the expansion as $\varepsilon \to 0$.
\end{lemma}

\proof
All our estimates in this proof are uniform with respect to $\xi\in M$, $t$ in compact subsets of $(0,\infty)$, and $\varepsilon\in\(0,\varepsilon_0\)$ for some fixed $\varepsilon_0>0$. Since $\widehat{W}_{\varepsilon,t,\xi}$ is a constant and $L_{g_\xi} G_{g_\xi}(\cdot,\xi)=0$ in $M\setminus B_\xi(r_0)$, by \eqref{confinv} and \eqref{Eq7}, we get that
\begin{eqnarray*}
&&\int_M\big|\nabla W_{\varepsilon,t,\xi}\big|^2_g dv_g+\alpha_n \int_M \Scal_g W_{\varepsilon,t,\xi}^2dv_g \nonumber \\
&&\qquad=\int_M\(L_g W_{\varepsilon,t,\xi}\)W_{\varepsilon,t,\xi} dv_g +\int_{\partial B_\xi\(r_0\)}\(\partial_{\nu_{\inward}}W_{\varepsilon,t,\xi}+\partial_{\nu_{\outward}}W_{\varepsilon,t,\xi}\) W_{\varepsilon,t,\xi} d\sigma_g \nonumber\allowdisplaybreaks\\
&&\qquad=\int_M L_{g_\xi} \big(G_{g_\xi}(\cdot,\xi) \widehat{W}_{\varepsilon,t,\xi}\big) G_{g_\xi}(\cdot,\xi) \widehat{W}_{\varepsilon,t,\xi} dv_{g_\xi}+\operatorname{O}\(\delta_\varepsilon(t)^n \) \nonumber\\
&&\qquad=\int_{B_\xi(r_0)} L_{g_\xi} \big(G_{g_\xi}(\cdot,\xi) \widehat{W}_{\varepsilon,t,\xi}\big) G_{g_\xi}(\cdot,\xi) \widehat{W}_{\varepsilon,t,\xi} dv_{g_\xi}+\operatorname{O}\(\delta_\varepsilon(t)^n \)
\end{eqnarray*}
in view of \eqref{Pr2Eq4} and $W_{\varepsilon,t,\xi}=\operatorname{O}\big(\delta_\varepsilon(t)^{\frac{n-2}{2}} \big)$ on $\partial B_\xi\(r_0\)$. In the estimates below, we make use of \eqref{added1} along with \eqref{Eq6}, Lemma~\ref{Lem4}, and
$$G_{g_\xi}(x,\xi) \widehat{W}_{\varepsilon,t,\xi}(x)=\beta_n G_{g_\xi}(x,\xi)d_{g_\xi}(x,\xi)^{n-2} [n(n-2)]^{\frac{n-2}{4}}\frac{\delta_\varepsilon(t)^{\frac{n-2}{2}}}{(\delta_\varepsilon(t)^2+d_{g_\xi}(x,\xi)^2)^{\frac{n-2}{2}}}$$
for all $x\in B_\xi\(r_0\)$. When $n=4,5$ or $\(M,g\)$ is l.c.f., we can deduce that
\begin{eqnarray}\label{Lem1Eq6}
&& \int_M\big|\nabla W_{\varepsilon,t,\xi}\big|^2_g dv_g+\alpha_n \int_M \Scal_g W_{\varepsilon,t,\xi}^2dv_g =
\frac{1}{2} [n(n-2)]^{\frac{n}{2}} \omega_{n-1} I_n^{\frac{n-2}{2}} \nonumber \\
&&\qquad+n^{\frac{n-2}{2}}\(n-2\)^{\frac{n}{2}} \beta_n   \omega_{n-1} I_n^{n-2} A_\xi \delta_\varepsilon\(t\)^{n-2} +\operatorname{O}\(\delta_\varepsilon(t)^{n-1}\),
\end{eqnarray}
where we denote $I^q_p:=\int_0^{+\infty}\(1+r\)^{-p}r^qdr$ for all $p,\, q$ such that $p-q>1$, and we use that $I^q_p=\frac{q}{p-q-1}I^{q-1}_p=\frac{p}{p-q-1}I^q_{p+1}$. Concerning the remaining cases, we find that
\begin{align}\label{Lem1Eq8}
\int_M\big|\nabla W_{\varepsilon,t,\xi}\big|^2_g dv_g+\alpha_n \int_M \Scal_g W_{\varepsilon,t,\xi}^2dv_g&=6912 \omega_5 I_6^2\\
&\quad-\frac{8}{5}\omega_5 I_6^4 \big|\Weyl_g\(\xi\)\big|_g^2\delta_\varepsilon\(t\)^4\ln\delta_\varepsilon\(t\)+\operatorname{O}\(\delta_\varepsilon(t)^4 \)\nonumber
\end{align}
when $n=6$ and $\(M,g\)$ is non-l.c.f., and
\begin{eqnarray} \label{Lem1Eq10}
&&\int_M\big|\nabla W_{\varepsilon,t,\xi}\big|^2_g dv_g+\alpha_n \int_M \Scal_g W_{\varepsilon,t,\xi}^2dv_g =\frac{1}{2}[n (n-2)]^{\frac{n}{2}}\omega_{n-1}I_n^{\frac{n-2}{2}}\\
&&+\frac{n^{\frac{n-2}{2}}(n-2)^{\frac{n+4}{2}}}{48(n+2)(n-1)}\omega_{n-1}I_n^{\frac{n+2}{2}}\delta_\varepsilon(t)^4 \Bigg(\frac{\big|\Weyl_g(\xi)\big|_g^2}{12(n-6)}+\frac{1}{n}\Delta_{\Eucl}(\Scal_{g_\xi} \circ \exp_\xi)(0)\Bigg)+\operatorname{O}(\delta_\varepsilon(t)^5) \nonumber
\end{eqnarray}
when $n\ge7$ and $\(M,g\)$ is non-l.c.f., in view of symmetry properties. Now, we estimate the term
\begin{align*}
\int_M W_{\varepsilon,t,\xi}^{2^*}dv_g&=\int_M \( G_{g_\xi}(\cdot,\xi) \widehat{W}_{\varepsilon,t,\xi}\)^{2^*} dv_{g_\xi}\\
&=[n(n-2)]^{\frac{n}{2}}\int_{B_0(r_0)} \[\beta_n G_{g_\xi}(\exp_\xi y,\xi) |y|^{n-2}\]^{2^*} \frac{\delta_\varepsilon(t)^n}{(\delta_\varepsilon(t)^2+|y|^2)^n}dy
+\operatorname{O}\(\delta_\varepsilon(t)^n \)
\end{align*}
in view of \eqref{Eq6}. By Lemma~\ref{Lem4}, we get that
\begin{eqnarray} \label{Lem1Eq7-9-11}
&&\int_M W_{\varepsilon,t,\xi}^{2^*}dv_g=\frac{1}{2}[n(n-2)]^{\frac{n}{2}}\omega_{n-1}I_n^{\frac{n-2}{2}}\\
&&+\left\{ \begin{aligned}
&n^{\frac{n+2}{2}}\(n-2\)^{\frac{n-2}{2}}\beta_n\omega_{n-1}I_n^{n-2}A_\xi\delta_\varepsilon\(t\)^{n-2}+\operatorname{O}\(\delta_\varepsilon(t)^{n-1}\) && \hbox{if }n=4,5 \hbox{ or }\(M,g\) \hbox{ l.c.f.}\\
&-\frac{72}{5}\omega_5I_6^4\big|\Weyl_g\(\xi\)\big|_g^2\delta_\varepsilon\(t\)^4\ln\delta_\varepsilon\(t\)+\operatorname{O}\(\delta_\varepsilon(t)^4 \) && \hbox{if }n=6 \hbox{ and }\(M,g\) \hbox{ non-l.c.f.}\\
&\frac{n^{\frac{n+2}{2}}(n-2)^{\frac{n}{2}}}{48(n-4)(n-1)}\omega_{n-1}I_n^{\frac{n+2}{2}}\delta_\varepsilon\(t\)^4
\Bigg(\frac{\big|\Weyl_g\(\xi\)\big|_g^2 }{12\(n-6\)} &&\\
&\qquad+\frac{1}{n}\Delta_{\Eucl}(\Scal_{g_\xi}\circ \exp_\xi)(0)\Bigg)+\operatorname{O}\(\delta_\varepsilon(t)^5 \), && \hbox{if }n\geq 7 \hbox{ and }\(M,g\) \hbox{ non-l.c.f.}\end{aligned} \right. \nonumber
\end{eqnarray}
by using symmetry properties when $n\geq 7$ and $\(M,g\)$ is non-l.c.f. Moreover, we find that
\begin{equation}\label{Lem1Eq12}
\int_MhW_{\varepsilon,t,\xi}^2dv_g=\left\{\begin{aligned}
&-8\omega_3h\(\xi\)\delta_\varepsilon\(t\)^2\ln\delta_\varepsilon\(t\)+\operatorname{O}\(\delta_\varepsilon(t)^2\)&&\text{if }n=4\\
&\frac{1}{2}[n(n-2)]^{\frac{n-2}{2}}\omega_{n-1}I_{n-2}^{\frac{n-2}{2}}h\(\xi\)\delta_\varepsilon\(t\)^2+\operatorname{O}\big(\delta_\varepsilon(t)^{\frac{5}{2}}\big)&&\text{if }n\ge5
\end{aligned}\right.
\end{equation}
as $\varepsilon\to0$, in view of $\varLambda_\xi(\xi)=1$, \eqref{Eq6}, and Lemma~\ref{Lem4}. Successive integrations by parts give that
\begin{equation}\label{Lem1Eq13}
I_n^{\frac{n-2}{2}}=\frac{\omega_n}{2^{n-1}\omega_{n-1}}\,,\,\ I_n^{n-2}=\frac{1}{n-1}\,,
\end{equation}
and (if $n\ge5$)
\begin{equation}\label{Lem1Eq14}
I_{n-2}^{\frac{n-2}{2}}=\frac{\(n-1\)\omega_n}{2^{n-3}\(n-4\)\omega_{n-1}}\,,\,\ I_n^{\frac{n+2}{2}}=\frac{n\(n+2\)\omega_n}{2^{n-1}\(n-2\)\(n-4\)\omega_{n-1}}\,.
\end{equation}
Moreover, for any $\xi\in M$, since $g_\xi$ defines conformal normal coordinates of order $N\ge5$, see Lee--Parker~\cite{LeePar}*{Theorem~5.1}, and since $\varLambda_\xi\(\xi\)=1$, we get that
\begin{equation}\label{Lem1Eq15}
\Delta_{\Eucl}(\Scal_{g_\xi}\circ \exp_\xi)(0)=\frac{1}{6}\big|\Weyl_{g_\xi}\(\xi\)\big|^2_{g_\xi}=\frac{1}{6}\big|\Weyl_g\(\xi\)\big|^2_g\,.
\end{equation}
Finally, \eqref{Lem1Eq1}--\eqref{Lem1Eq4} follow from \eqref{Lem1Eq6}--\eqref{Lem1Eq15} in view of \eqref{Eq24}.
\endproof

We end this section by proving the validity of the expansion for $\mathcal{J}_\varepsilon$ in Proposition~\ref{Pr3} uniformly with respect to $\xi \in M$, $t$ in compact subsets of $(0,\infty)$ and $\varepsilon\in\(0,\varepsilon_0\)$. To this aim, it suffices to prove the (uniform) expansion
\begin{equation}\label{Eq26}
\mathcal{J}_\varepsilon\(t,\xi\)=J_\varepsilon\(W_{\varepsilon,t,\xi}\)+\left\{\begin{aligned}&\operatorname{o}\(\delta_\varepsilon(t)^2\)&&\text{if }n=4\\&\operatorname{o}\(\varepsilon\delta_\varepsilon(t)^2\)&&\text{if }n\ge5\end{aligned}\right.
\end{equation}
as $\varepsilon\to0$. Indeed, the expansion of $\mathcal{J}_\varepsilon$ in Proposition~\ref{Pr3} follows from combining Lemma~\ref{Lem1} with \eqref{Eq26} and letting $\delta_\varepsilon(t)$ be as in \eqref{Eq8}.  Since
\begin{multline}\label{Pr3Step1Eq1}
\mathcal{J}_\varepsilon\(t,\xi\)-J_\varepsilon\(W_{\varepsilon,t,\xi}\)=\<W_{\varepsilon,t,\xi}-L_g^{-1}\(f_\varepsilon\(W_{\varepsilon,t,\xi}\)\),\phi_{\varepsilon,t,\xi}\>_{L_g}+\frac{1}{2}\(\left\|\phi_{\varepsilon,t,\xi}\right\|^2_{L_g}+\varepsilon \int_M h \phi_{\varepsilon,t,\xi}^2 dv_g \)\\
-\frac{1}{2^*}\int_M\(\(W_{\varepsilon,t,\xi}+\phi_{\varepsilon,t,\xi}\)_+^{2^*}-W_{\varepsilon,t,\xi}^{2^*}-2^*W_{\varepsilon,t,\xi}^{2^*-1}\phi_{\varepsilon,t,\xi}\)dv_g\,,
\end{multline}
by Cauchy--Schwarz inequality, we get that
\begin{equation}\label{Pr3Step1Eq2}
\left|\<W_{\varepsilon,t,\xi}-L_g^{-1}\(f_\varepsilon\(W_{\varepsilon,t,\xi}\)\),\phi_{\varepsilon,t,\xi}\>_{L_g}\right|\le\left\|W_{\varepsilon,t,\xi}-L_g^{-1}\(f_\varepsilon\(W_{\varepsilon,t,\xi}\)\)\right\|_{L_g}\left\|\phi_{\varepsilon,t,\xi}\right\|_{L_g},\end{equation}
and by the Mean Value Theorem and the H\"older's inequality, we get that
\begin{multline}\label{Pr3Step1Eq3}
\int_M\(\(W_{\varepsilon,t,\xi}+\phi_{\varepsilon,t,\xi}\)_+^{2^*}-W_{\varepsilon,t,\xi}^{2^*}-2^*W_{\varepsilon,t,\xi}^{2^*-1}\phi_{\varepsilon,t,\xi}\)dv_g\\
=\operatorname{O}\(\left\|\phi_{\varepsilon,t,\xi}\right\|^2_{2^*}\left\|W_{\varepsilon,t,\xi}\right\|^{2^*-2}_{2^*}+\left\|\phi_{\varepsilon,t,\xi}\right\|^{2^*}_{2^*}\)=\operatorname{O}\(\left\|\phi_{\varepsilon,t,\xi}\right\|^2_{2^*}+\left\|\phi_{\varepsilon,t,\xi}\right\|^{2^*}_{2^*}\).
\end{multline}
By Proposition~\ref{Pr1} and \eqref{Pr3Step1Eq1}--\eqref{Pr3Step1Eq3}, it follows that
$$\mathcal{J}_\varepsilon\(t,\xi\)-J_\varepsilon\(W_{\varepsilon,t,\xi}\)=\operatorname{O}\(\|R_{\varepsilon,t,\xi}\|_{1,2}^2\)$$
in view of Sobolev's embedding $L^{2^*}(M) \hookrightarrow H_1^2(M)$ and the equivalence between the two norms $\|\cdot\|_{L_g}$ and $\|\cdot \|_{1,2}$. Proposition~\ref{Pr2} now yields that the estimate \eqref{Eq26} does hold $C^0$--uniformly.

\section{Existence of $k$--bubbles}\label{Sec5}

The previous analysis can be extended to solutions which have $k$ distinct blow-up points, $k \geq 2$. In the non-l.c.f.~case with $n\geq 6$, the ``reduced energy'' $\widetilde E_k:(0,\infty)^k \times (M^k \setminus \Delta_k) \to \mathbb{R}$, $\Delta_k:=\left\{(\xi_1,\dots,\xi_k)\in M^k \ :\ \xi_i=\xi_j\ \hbox{for}\ i\not=j\right\}$, which governs the location of these blow-up points, is just a super-position of the one for each single point:
$$\widetilde E_k( \bgl{t}, \bgl{\xi}):=\sum_{i=1}^k \widetilde E(t_i,\xi_i),\qquad \bgl{t}:=(t_1,\dots,t_k)\in (0,\infty)^k,\quad \bgl{\xi}:=(\xi_1,\dots,\xi_k)\in M^k \setminus \Delta_k\,.$$
Theorem~\ref{Th} works as well in this context in the following way: as soon as we find $k$ distinct $C^p$--{\it stable}, $p=0,1$, {\it critical sets} $\mathcal{\widetilde D}_1, \dots, \mathcal{\widetilde D}_k$ of $\widetilde E(t,\xi)$, we can construct a  family {$\(u_{k,\varepsilon}\)_\varepsilon$} of solutions to \eqref{Eq1} which blows up, up to a sub-sequence, at points $\(\xi_0\)_1,\dots, \(\xi_0\)_k$ with $\bgl{\xi}_0 \in  \mathcal{\widetilde D}_1 \times \dots \times \mathcal{\widetilde D}_k$ as $\varepsilon \to 0$. { From this result, we deduce Theorem~\ref{ThThThThNonlcf} exactly as in the case $k=1$.}

\medskip
{ Theorem~\ref{ThThTh} has its counter-part too: in the non-l.c.f.~case with $n\geq 6$, solutions with $k$ blow-up points do exist provided that $k$ is at most the number of isolated critical points of $E(\xi)=h\(\xi\) \left|\Weyl_g\(\xi\)\right|^{-1}_g$ with non-trivial local degree and $h>0$.}

\medskip
When $n=4$, the energy for the approximating function \eqref{Eq2m} is not suitable due to the dependence in $t$ of the smallness rate of $\delta_\varepsilon\(t\)$ in \eqref{Eq8}.

\medskip When $n=5$ or $\(M,g\)$ is l.c.f., the picture is completely different. There is an effective interaction between different blow-up points as expressed by the following ``reduced energy'' $\widetilde E_k :(0,\infty)^k \times (M^k \setminus \Delta_k) \to \mathbb{R}$:
\begin{equation}\label{eridk2}
\widetilde E_k(\bgl{t}, \bgl{\xi})=
c_2 \displaystyle \sum_{i=1}^k t_i^2 h(\xi_i)-c_3 \sum_{i=1}^k t_i^{n-2} A_{\xi_i} -c_3 \sum_{i,j=1 \atop i\not=j}^k t_i^{\frac{n-2}{2}} t_j^{\frac{n-2}{2}} G_g (\xi_i,\xi_j),
\end{equation}
where $c_2,c_3>0$ and $A_\xi$ is as in \eqref{Eq3}. Theorem~\ref{Th} is still valid in this context by simply replacing $\widetilde E(t,\xi)$ with $\widetilde E_k(\bgl{t}, \bgl{\xi})$. It is no longer possible in general to relate critical sets of $\widetilde E_k(\bgl{t}, \bgl{\xi})$ with that of an explicit $E_k(\bgl{\xi})$ as when $k=1$, with the exception of the case $n=6$ for which we prove Theorem~\ref{ThThThThTh} with $E_k$ defined as in \eqref{matrix2}. In case $n\geq 7$ (in such a way that $\frac{n-2}{2}>2$), we can produce a $C^0$--{\it stable critical set} of $\widetilde E_k$ through its maximal set, yielding to Theorem~\ref{ThThThTh}. In this section, we first sketch the proof of Theorem~\ref{Th} in case $k \geq 2$, with $\widetilde E(t,\xi)$ replaced by $\widetilde E_k(\bgl{t}, \bgl{\xi})$, and we then derive from it Theorems~\ref{ThThThTh} and~\ref{ThThThThTh}.

\medskip
For $\bgl{t}:=(t_1,\dots,t_k)\in (0,\infty)^k$ and $\bgl{\xi}:=(\xi_1,\dots,\xi_k)\in M^k\setminus \Delta_k$, define $\delta_\varepsilon(t_i)$, $i=1,\dots,k$, as in \eqref{Eq8}. The $k$--bubbles approximating function is given by
\begin{equation}\label{Eq2m}
W_{\varepsilon,\bgl{t} ,\bgl{\xi} }:=\sum\limits_{i=1 }^k W_{\varepsilon,t_i,\xi_i }\,,
\end{equation}
where $W_{\varepsilon,t_i,\xi_i }$ is defined in \eqref{Eq9}--\eqref{Eq10}. Since $\varLambda_\xi$ are positive functions depending smoothly in $\xi \in M$, there exists $C_0>0$ so that
$C_0^{-1}\leq \Lambda_\xi(x)\leq C_0$ for all $x,\xi \in M$, and then
\begin{equation} \label{comparedist}
C_0^{-\frac{2}{n-2}}d_g(x,y) \leq d_{g_\xi}(x,y)\leq C_0^{\frac{2}{n-2}} d_g(x,y)
\end{equation}
for all $x,y,\xi\in M$. The number $r_0$ in \eqref{Eq10} is also assumed to satisfy
 $r_0< C_0^{-\frac{2}{n-2}} {d_g(\xi_i,\xi_j)\over 2}$ for all $i\not=j$ in such a way that $\{d_{g_{\xi_i}}(x,\xi_i)\leq r_0\} \cap \{d_{g_{\xi_j}}(x,\xi_j)\leq r_0\} =\emptyset$ in view of \eqref{comparedist}. We look for a solution of  \eqref{Eq1} in the form
$$ u_{k,\varepsilon}:=W_{\varepsilon,\bgl{t},\bgl{\xi}}+\phi_\varepsilon\,,$$
where $\bgl{t} \in(0,\infty)^k$, $\bgl\xi \in M^k\setminus \Delta_k$, and $\phi_\varepsilon\in K^\perp_{\varepsilon,\bgl t,\bgl \xi}$ with
\begin{align*}
K_{\varepsilon,\bgl t,\bgl \xi}:=\bigcup_{i=1}^k K_{\varepsilon, t_i,  \xi_i}\quad\text{and}\quad K^\perp_{\varepsilon,\bgl t,\bgl\xi} :=\bigcap_{i=1}^k K^\perp_{\varepsilon, t_i,  \xi_i} \,,
\end{align*}
where $K_{\varepsilon, t_i,  \xi_i}$ and $K^\perp_{\varepsilon, t_i,  \xi_i}$ are as in \eqref{Eq16}--\eqref{Eq17}. Letting  $\varPi_{\varepsilon,\bgl{t},\bgl{\xi}}$ and $\varPi^\perp_{\varepsilon,\bgl{t},\bgl{\xi}}$ be the respective projections of $H^2_1\(M\)$ onto $K_{\varepsilon,\bgl{t},\bgl{\xi}}$ and $K^\perp_{\varepsilon,\bgl{t},\bgl{\xi}}$, we rewrite equation \eqref{Eq1} as the system \eqref{Eq19}--\eqref{Eq20}, with $W_{\varepsilon,t,\xi}$, $\varPi_{\varepsilon,t,\xi}$, $\varPi^\perp_{\varepsilon,t,\xi}$ replaced by $W_{\varepsilon,\bgl{t},\bgl{\xi}}$, $\varPi_{\varepsilon,\bgl{t},\bgl{\xi}}$, $\varPi^\perp_{\varepsilon,\bgl{t},\bgl{\xi}}$, respectively. We begin with solving equation \eqref{Eq20} in Proposition~\ref{Pr1m} below, which is, as already observed, a well known result in this context (see for instance Musso--Pistoia~\cite{mupi4}).

\begin{proposition}\label{Pr1m}
Given  positive real numbers $a<b$ and $\eta$, there exists a positive constant $C=C\(a,b,\eta,k,n,M,g,h\)$ such that for $\varepsilon>0$ small, for any $\bgl t\in\[a,b\]^k$ and $\bgl \xi\in M^k$ with $d_g(\xi_i,\xi_j)\ge\eta$ for all $i\not=j$, there exists a unique function $\phi_{\varepsilon,\bgl t,\bgl \xi}\in K^\perp_{\varepsilon,\bgl t,\bgl \xi}$ which solves equation \eqref{Eq20} and satisfies
$$\big\|\phi_{\varepsilon,\bgl t,\bgl \xi}\big\|_{1,2}\le C\big\|R_{\varepsilon,\bgl t,\bgl \xi}\big\|_{1,2}\,,$$
where $R_{\varepsilon,\bgl t,\bgl \xi}:=W_{\varepsilon,\bgl t,\bgl \xi}-L_g^{-1}\big(f_\varepsilon\big(W_{\varepsilon,\bgl t,\bgl \xi}\big)\big)$. Moreover, $\phi_{\varepsilon,\bgl t,\bgl \xi}$ is continuously differentiable with respect to $\bgl t$ and $\bgl \xi$.
\end{proposition}

We now give an estimate for $\big\|R_{\varepsilon,\bgl t,\bgl \xi}\big\|_{1,2}$.

\begin{proposition}\label{Pr2m}
Given   positive real numbers $a<b$ and $\eta$, there exists a positive constant $C=C\(a,b,\eta,k,n,M,g,h\)$ such that for $\varepsilon>0$ small, for any $\bgl t\in\[a,b\]^k$ and $\bgl \xi\in M^k$ with $d_g(\xi_i,\xi_j)\ge\eta$ for all $i\not=j$,  there holds
$$\big\|R_{\varepsilon,\bgl t,\bgl \xi}\big\|_{1,2}\le C\left\{\begin{aligned}&\varepsilon ^{ 5\over 2}&&\text{if }n= 5\\
& \varepsilon^2\left |\ln \varepsilon\right|^{\frac{2}{3}} &&\text{if }n=6\text{ and }\(M,g\)\text{ l.c.f.}\\
& \varepsilon^{ n+2\over 2(n-4)}&&\text{if }n\ge7\text{ and }\(M,g\)\text{ l.c.f.}
 \end{aligned}\right. $$
where $R_{\varepsilon,\bgl t,\bgl \xi}$ is as in Proposition~\ref{Pr1m}.
\end{proposition}

\proof
We argue exactly as in the proof of Proposition~\ref{Pr2}. We point out that in this case
\begin{multline*}
\big<L_g^{-1}\big(f_\varepsilon\big(W_{\varepsilon,\bgl t,\bgl \xi}\big)\big)-W_{\varepsilon,\bgl t,\bgl \xi},\phi\big>_{L_g}=\sum\limits_{i=1}^k
\big<L_g^{-1}\big(f_\varepsilon\big(W_{\varepsilon,  t_i,  \xi_i}\big)\big)-W_{\varepsilon,  t_i,  \xi_i},\phi\big>_{L_g}\\
+\int_M\bigg[ \bigg(\sum\limits_{i=1}^k W_{\varepsilon,t_i,\xi_i}\bigg)^{2^*-1}-\sum\limits_{i=1}^k   W_{\varepsilon,t_i,\xi_i}^{2^*-1} \bigg]\phi dv_g\,,
\end{multline*}
and then
\begin{multline*}
\big\|L_g^{-1}\big(f_\varepsilon\big(W_{\varepsilon,\bgl t,\bgl \xi}\big)\big)-W_{\varepsilon,\bgl t,\bgl \xi}\big\|_{1,2}=\operatorname{O}\bigg( \sum\limits_{i=1}^k \left\|L_g^{-1}\(f_\varepsilon\(W_{\varepsilon,  t_i,  \xi_i}\)\)-W_{\varepsilon,  t_i,  \xi_i}\right\|_{1,2}\bigg)\\
+\operatorname{O}\Bigg(\bigg\|\bigg(\sum\limits_{i=1}^k W_{\varepsilon,t_i,\xi_i}\bigg)^{2^*-1}-\sum\limits_{i=1}^k   W_{\varepsilon,t_i,\xi_i}^{2^*-1}\bigg\|_{\frac{2n}{n+2}}\Bigg).
\end{multline*}
The first $k$ terms are estimated in Proposition~\ref{Pr2}, and the last term can be estimated following the arguments used for \eqref{new1}.
\endproof

For  any $\bgl t \in (0,\infty)^k$ and $\bgl \xi\in M^k \setminus \Delta_k$ we define
$$\mathcal{J}_\varepsilon\(\bgl t,\bgl \xi\):=J_\varepsilon\big(W_{\varepsilon,\bgl t,\bgl \xi}+\phi_{\varepsilon,\bgl t,\bgl \xi}\big),$$
where $J_\varepsilon$ is defined in \eqref{Eq22} and $\phi_{\varepsilon,\bgl t,\bgl \xi}$ is given by Proposition~\ref{Pr1m}. As already observed for $k=1$, we can solve equation \eqref{Eq19} by searching critical points of $\mathcal{J}_\varepsilon$, and the asymptotic expansion given in Proposition~\ref{Pr3m} below is crucial.

\begin{proposition}\label{Pr3m}
There holds
$$\mathcal{J}_\varepsilon\(\bgl t,\bgl \xi\)=c_1+ \varepsilon^{n-2\over n-4} \widetilde E_k(\bgl{t}, \bgl{\xi}) +\operatorname{o}\big(\varepsilon^{n-2\over n-4}\big)$$
as $\varepsilon\to0$, uniformly with respect to $\bgl \xi$ in compact subsets of $ M^k\setminus \Delta_k$ and $\bgl t$ in compact subsets of $(0,\infty)^k$, where $\widetilde E_k$ is given by \eqref{eridk2} and $c_1>0$ depends only in $n$.
\end{proposition}

\proof
We argue exactly as in the proof of Proposition~\ref{Pr3}, taking into account Lemma~\ref{Lem1m} below and exploiting \eqref{Eq8}.
\endproof

\begin{lemma}\label{Lem1m}
Assume that either $n=5$ or [$n\ge6$ and $\(M,g\)$ l.c.f.]. We let $K_n$ be as in \eqref{Eq24}, $\delta_\varepsilon\(t\)$ be as in \eqref{Eq8}, and $A_{\xi_i}$ be as in \eqref{Eq3}. As $\varepsilon\to0$, the following expansion does hold
\begin{multline}\label{Lem1Eq1m}
J_\varepsilon\big(W_{\varepsilon,\bgl t,\bgl \xi}\big)=\frac{k}{n}K_n^{-n}+[n(n-2)]^{\frac{n-2}{2}} \bigg[ \frac{(n-1)\omega_n}{2^{n-1}(n-4)} \varepsilon \sum_{i=1}^k h(\xi_i) \delta_\varepsilon(t_i)^2-\frac{\beta_n^2}{2} \sum_{i=1}^k A_{\xi_i}\delta_\varepsilon(t_i)^{n-2}\\
-\frac{\beta_n^2}{2} \sum_{i,j=1 \atop i\not=j}^k \( \delta_\varepsilon(t_i) \delta_\varepsilon(t_j)\)^{\frac{n-2}{2}} G_g(\xi_i,\xi_j)\bigg] +\operatorname{o}\big(\varepsilon^{n-2\over n-4}\big),
\end{multline}
uniformly with respect to $\bgl \xi$ in compact subsets of $ M^k\setminus \Delta_k$ and $\bgl t$ in compact subsets of $(0,\infty)^k$.
\end{lemma}

\proof
We proceed exactly as in the proof of Lemma~\ref{Lem1}. We only point out that
\begin{multline}\label{new0}
J_\varepsilon\big(W_{\varepsilon,\bgl t,\bgl \xi}\big)=\sum\limits_{i=1}^kJ_\varepsilon\(W_{\varepsilon, t_i, \xi_i}\)+{1\over2}\sum\limits_{i,j=1\atop i\not=j}^k
\<W_{\varepsilon, t_i, \xi_i},W_{\varepsilon, t_j, \xi_j }\>_{L_g} \\
+{\varepsilon \over2}\sum\limits_{i,j=1\atop i\not=j}^k\int_M h W_{\varepsilon, t_i, \xi_i}\ W_{\varepsilon, t_j, \xi_j }dv_g-{1\over2^*}\int_M \bigg[\bigg(\sum\limits_{i =1 }^k W_{\varepsilon, t_i, \xi_i}\bigg) ^{2^*}-\sum\limits_{i =1 }^k W_{\varepsilon, t_i, \xi_i}^{2^*}\bigg]dv_g\,.
\end{multline}
We claim that if $i\not=j$, then
\begin{equation}\label{new2}
\<W_{\varepsilon, t_i, \xi_i},W_{\varepsilon, t_j, \xi_j }\>_{L_g}=  [n(n-2)]^{n-2 \over 2} \beta_n^2 \(\delta_\varepsilon\(t_i\)\delta_\varepsilon\(t_j\)\)^{n-2\over2}G_g(\xi_i,\xi_j)+\operatorname{o}\big(\varepsilon^{n-2\over n-4}\big).
\end{equation}
Indeed, by \eqref{Eq6}, \eqref{Eq8}, \eqref{Pr2Eq4}, \eqref{added1}, and Lemma~\ref{Lem4}, we get that
\begin{align*}
\<W_{\varepsilon, t_i, \xi_i},W_{\varepsilon, t_j, \xi_j }\>_{L_g}&=\int_M
\(L_gW_{\varepsilon, t_i, \xi_i}\)W_{\varepsilon, t_j, \xi_j }dv_{ {g} }\\
&\quad+\int_{\partial B_{\xi_i}\(r_0\)}\(\partial_{\nu_{\inward}}W_{\varepsilon,t_i,\xi_i}+\partial_{\nu_{\outward}}W_{\varepsilon,t_i,\xi_i} \)W_{\varepsilon, t_j, \xi_j } d\sigma_g\,\\
&=\beta_n\delta_\varepsilon(t_j)^{\frac{2-n}{2}}r_0^{n-2} U\(\delta_\varepsilon(t_j)^{-1}r_0\)\\
&\quad\times\int_{B_{\xi_i}\(r_0\)}L_{g_{\xi_i}}\big(G_{g_{\xi_i}}(\cdot,\xi_i)\widehat{W}_{\varepsilon,t_i,\xi_i}\big)\varLambda_{\xi_i}(x)^{-1}G_g(x,\xi_j)dv_{g_{\xi_i}} +\operatorname{o}\big(\varepsilon^{n-2\over n-4}\big)\allowdisplaybreaks\\
&=[n(n-2)]^{n-2\over 4}  \beta_n \delta_\varepsilon(t_j)^{n-2\over2}G_g(\xi_i,\xi_j)\\
&\quad\times\int_{B_0 (r_0)}
L_{g_{\xi_i}}\big(G_{g_{\xi_i}}(\cdot,\xi_i) \widehat{W}_{\varepsilon,t_i,\xi_i}\big) (\exp_{\xi_i} y) [1+\operatorname{o}(1)+\operatorname{O}(|y|)] dy+\operatorname{o}\big(\varepsilon^{n-2\over n-4}\big)\\
&=\frac{1}{2}[n(n-2)]^{n\over 2}\beta_n\omega_{n-1} \(\delta_\varepsilon\(t_i\)\delta_\varepsilon\(t_j\)\)^{n-2\over2}G_g (\xi_i,\xi_j)
I_{\frac{n+2}{2}}^{\frac{n-2}{2}} +\operatorname{o}\big(\varepsilon^{n-2\over n-4}\big)
\end{align*}
in view of $\varLambda_{\xi_i}\(\xi_i\)=1$, $L_g G_g(x,\xi_i)=0$, and
\begin{equation} \label{new1added}
W_{\varepsilon, t_i, \xi_i}=\beta_nG_g(x,\xi_i)\delta_\varepsilon(t_i)^{\frac{2-n}{2}}r_0^{n-2}U\(\delta_\varepsilon(t_i)^{-1}r_0\)
\end{equation}
for all $x\in M\setminus B_{\xi_i}(r_0)$. Since
\begin{equation}\label{new3}
I_{\frac{n+2}{2}}^{\frac{n-2}{2}}= \int_0^{+\infty}\frac{r^{\frac{n-2}{2}}}{(1+r)^{\frac{n+2}{2}}}dr= \int_0^{+\infty} (1-\frac{1}{r+1})^{\frac{n-2}{2}} \frac{dr}{(1+r)^2}= \int_0^1 (1-s)^{\frac{n-2}{2}}ds=\frac{2}{n}\,,\end{equation}
we deduce the validity of \eqref{new2}. We claim that
\begin{multline} \label{new1}
{1\over2^*}\int_M \bigg[\bigg(\sum\limits_{i =1 }^k W_{\varepsilon, t_i, \xi_i}\bigg) ^{2^*}-\sum\limits_{i =1 }^k W_{\varepsilon, t_i, \xi_i}^{2^*}\bigg]dv_{ {g} }\\
= [n(n-2)]^{\frac{n-2}{2}}
 \beta_n^2  \sum\limits_{i,j=1\atop i\not=j}^k\(\delta_\varepsilon\(t_i\)\delta_\varepsilon\(t_j\)\)^{n-2\over2}G_g (\xi_i,\xi_j) +\operatorname{o}\big(\varepsilon^{n-2\over n-4}\big).
\end{multline}
Indeed, by \eqref{Eq6}, \eqref{Eq8}, \eqref{new1added}, and Lemma~\ref{Lem4} , we deduce that
\begin{align*}
&\int_M \bigg[\bigg(\sum\limits_{i =1}^k W_{\varepsilon, t_i, \xi_i}\bigg)^{2^*} - \sum\limits_{i =1 }^k W_{\varepsilon, t_i, \xi_i}^{2^*} \bigg]dv_g\\
&\qquad= \sum\limits_{j=1}^k \int_{B_{\xi_j}(r_0)} \bigg[\bigg(\sum\limits_{i =1 }^k W_{\varepsilon, t_i, \xi_i}\bigg)^{2^*} - W_{\varepsilon, t_j, \xi_j}^{2^*} \bigg]dv_g+\operatorname{o}\big(\varepsilon^{n-2\over n-4}\big)\allowdisplaybreaks\\
&\qquad=2^* [n(n-2)]^{\frac{n-2}{4}} \beta_n\sum\limits_{i,j=1 \atop i \not=j}^k\delta_\varepsilon(t_i)^{\frac{n-2}{2}}(1+\operatorname{o}(1))\int_{B_{\xi_j}(r_0)} G_g(x,\xi_i)  W_{\varepsilon, t_j, \xi_j} ^{2^*-1} dv_g +\operatorname{o}\big(\varepsilon^{n-2\over n-4}\big)\\
&\qquad=\frac{2^*}{2} [n(n-2)]^{n\over2} \beta_n\omega_{n-1}  I_{\frac{n+2}{2}}^{\frac{n-2}{2}} \sum\limits_{i,j=1 \atop i\not=j}^k\(\delta_\varepsilon\(t_i\)\delta_\varepsilon\(t_j\)\)^{n-2\over2}G_g(\xi_j,\xi_i) +\operatorname{o}\big(\varepsilon^{n-2\over n-4}\big)
\end{align*}
in view of the estimate $(a+b)^{2^*}-a^{2^*}-2^*a^{2^*-1}b=\operatorname{O}(a^{2^*-2}b^2+b^{2^*})$ for all $a,b\geq 0$. Therefore, thanks to \eqref{new3}, we deduce the validity of \eqref{new1}. Moreover, by \eqref{Eq8} and \eqref{new1added}, we get that
\begin{multline}\label{new4}
\int_M h W_{\varepsilon, t_i, \xi_i}\ W_{\varepsilon, t_j, \xi_j }dv_g=\operatorname{O}\bigg(\delta_\varepsilon(t_j)^{\frac{n-2}{2}} \int_{B_{\xi_i}(r_0)} W_{\varepsilon, t_i, \xi_i}dv_g+\delta_\varepsilon(t_i)^{\frac{n-2}{2}} \int_{B_{\xi_j}(r_0)} W_{\varepsilon, t_j, \xi_j}dv_g \bigg)\\
+\operatorname{O}\(\delta_\varepsilon(t_i)^{\frac{n-2}{2}} \delta_\varepsilon(t_j)^{\frac{n-2}{2}}\)=\operatorname{O}\big(\varepsilon^{\frac{n-2}{n-4}}\big).
\end{multline}
Inserting \eqref{new2}, \eqref{new1}, and \eqref{new4} into \eqref{new0} and combining with Lemma~\ref{Lem1} we finally deduce the validity of \eqref{Lem1Eq1m}.
\endproof

\smallskip
\proof[Proof of Theorem~\ref{ThThThTh}]
The key point is to show that $\widetilde E_k$ attains its maximum value
$$m_k:=\sup\limits_{(0,\infty)^k \times (M^k\setminus \Delta_k)} \widetilde E_k$$
at interior points, i.e.
$$ \mathcal{\widetilde D}_k=\{(\bgl t,\bgl \xi)\in (0,\infty)^k \times (M^k\setminus \Delta_k): \widetilde E_k(\bgl t,\bgl \xi)=m_k \}$$
is a non-empty compact set. We claim that there exists $L>0$ so that
\begin{equation}\label{sup1}
\sup\limits_{(0,L]^l \times (M^l\setminus \Delta_l)}\widetilde E_l=m_l\quad\text{and}\quad\sup\limits_{((0,\infty)^l \setminus (0,L]^l) \times (M^l\setminus \Delta_l)}  \widetilde E_l<m_l
\end{equation}
for all $l=1,\dots,k$, and there holds
\begin{equation} \label{sup2}
0<m_1<m_2<\dots<m_k <+\infty\,.
\end{equation}
Indeed, since $\min\left\{A_\xi\,:\,\xi\in M\right\}>0$ and
$$\widetilde E_l(\bgl t,\bgl \xi)\leq\sum_{i=1}^l \( c_2  t_i^2 h(\xi_i)-c_3 t_i^{n-2} A_{\xi_i}\)$$
for all $(\bgl t,\bgl \xi)\in (0,\infty)^k \times (M^k\setminus \Delta_k)$, we have that $\widetilde E_l(\bgl t,\bgl \xi)\to -\infty$ uniformly as soon as $t_i \to +\infty$ for some $i=1,\dots,l$. Therefore, we can find some $L>0$ large so that \eqref{sup1} does hold, and $m_l<+\infty$. Since $\max_M h>0$, we can find $\xi_1 \in M$ with $h(\xi_1)>0$, and then $c_2 t_1^2 h(\xi_1)-c_3 t_1^{n-2} A_{\xi_1}>0$ for $t_1>0$ sufficiently small. It follows that $m_1>0$. To conclude the proof of \eqref{sup2}, observe that for $l\geq 2$, we have that
\begin{align}
\widetilde E_l(\bgl t,\bgl \xi)&=\widetilde E_{l-1}(t_1,\dots,t_{l-1}, \xi_1,\dots,\xi_{l-1})\nonumber\\
&\qquad+c_2 t_l^2 h(\xi_l)-c_3 t_l^{n-2} A_{\xi_l}- 2 c_3 t_l^{\frac{n-2}{2}} \sum_{i=1}^{l-1} t_i^{\frac{n-2}{2}}  G_g (\xi_i,\xi_l).\label{Ekeland0}
\end{align}
Let $(t_1^j,\dots,t_{l-1}^j, \xi_1^j,\dots,\xi_{l-1}^j)_{j\in\mathbb{N}}$ be a maximizing sequence for $m_{l-1}$. Up to a subsequence, we can assume that $\xi_i^j \to\overline\xi_i \in M$ as $j \to +\infty$ for all $i=1,\dots,l-1$. Now, we fix some $\xi_l \in \{h>0\}\setminus \{\overline \xi_1,\dots,\overline \xi_{l-1}\}$ and we choose $t_l>0$ small so that
\begin{align*}
&c_2 t_l^2 h(\xi_l)-c_3 t_l^{n-2} A_{\xi_l}- 2 c_3 t_l^{\frac{n-2}{2}} \sum_{i=1}^{l-1} t_i^{\frac{n-2}{2}} G_g (\overline \xi_i,\xi_l)\\
&\qquad\geq c_2 t_l^2 h(\xi_l)-c_3 t_l^{n-2} A_{\xi_l}- 2 c_3 L^{\frac{n-2}{2}} t_l^{\frac{n-2}{2}} \sum_{i=1}^{l-1}  G_g (\overline \xi_i,\xi_l)>0
\end{align*}
in view of $\frac{n-2}{2}>2$. Therefore, we get that
\begin{align*}
m_l&\geq\lim\limits_{j \to +\infty} \widetilde E_l(t_1^j,\dots,t_{l-1}^j,t^j_l, \xi_1^j,\dots,\xi_{l-1}^j,\xi^j_l)\\
&=m_{l-1}+c_2 t_l^2 h(\xi_l)-c_3 t_l^{n-2} A_{\xi_l}- 2 c_3 t_l^{\frac{n-2}{2}} \sum_{i=1}^{l-1} t_i^{\frac{n-2}{2}} G_g (\overline \xi_i,\xi_l)>m_{l-1}\,,
\end{align*}
and \eqref{sup2} is established.

\medskip
Now, we prove that $\mathcal{\widetilde D}_k \not= \emptyset$ and that $\mathcal{\widetilde D}_k$ is a compact set. By \eqref{sup1}, we can find a maximizing sequence $(\bgl t^j,\bgl \xi^j)_{j\in\mathbb{N}}$, $(\bgl t^j,\bgl \xi^j):=(t_1^j,\dots,t_k^j, \xi_1^j,\dots,\xi_k^j)$, for $\widetilde E_k$ so that $t_i^j \leq L$ for all $i=1,\dots,k$ and $j\in\mathbb{N}$. By \eqref{Ekeland0} with $l=k$, we get that if $t_k^j\to0$ as $j\to +\infty$, then $m_k\leq m_{k-1}$ in contradiction with \eqref{sup2}. Since the same argument applies for all the $t_i^j$'s, we get that there exists $\eta>0$ so that $t_i^j \geq \eta$ for all $i=1,\dots,k$ and $j\in\mathbb{N}$. By compactness of $[\eta,L]^k \times M^k$, up to a subsequence, we can assume that $(\bgl t^j,\bgl \xi^j) \to (\bgl t_0,\bgl \xi_0)\in[\eta,L]^k \times M^k$ as $j \to +\infty$.  Since $G_g(x,y) \to +\infty$ as $d_g(x,y)\to 0$, we get that $\bgl \xi_0 \in M^k \setminus \Delta_k$, and thus that $(\bgl t_0,\bgl \xi_0) \in \mathcal{\widetilde D}_k$. Hence, we get that $\mathcal{\widetilde D}_k \not= \emptyset$. As for the compactness of $\mathcal{\widetilde D}_k$, we let $(\bgl t^j,\bgl \xi^j)_{j\in\mathbb{N}}$ be a sequence in $\mathcal{\widetilde D}_k$, and by the same arguments as above, we deduce that $(\bgl t^j,\bgl \xi^j) \to (\bgl t_0,\bgl \xi_0)\in[\eta,L]^k \times (M^k\setminus \Delta_k)$ as $j \to +\infty$, and by continuity of $\widetilde E_k$, $(\bgl t_0,\bgl \xi_0)\in\mathcal{\widetilde D}_k$, which proves that $\mathcal{\widetilde D}_k$ is a compact set.

\medskip
To conclude the proof, let $\widetilde U$ be a compact neighborhood of $\mathcal{\widetilde D}_k$ in $(0,\infty)^k \times (M^k \setminus \Delta_k)$. Since by \eqref{sup1} $\widetilde E_k$ is not a constant function, we have that
\begin{equation}\label{final1}
\sup_{\partial \widetilde U}\widetilde E_k<\sup_{\mathcal{\widetilde D}_k} \widetilde E_k\,.
\end{equation}
By Proposition~\ref{Pr3m}, we get that
\begin{equation}\label{final2}
\mathcal{\widetilde J}_\varepsilon(\bgl t,\bgl \xi):=\varepsilon^{-{n-2\over n-4}}\[\mathcal{J}_\varepsilon \(\bgl t,\bgl \xi\)-c_1\]\longrightarrow\widetilde E_k(\bgl t,\bgl \xi)
\end{equation}
as $\varepsilon\to0$, uniformly with respect to $\bgl \xi$ in compact subsets of $ M^k \setminus \Delta_k$ and $\bgl t$ in compact subsets of $(0,\infty)^k$. It follows from \eqref{final1} and \eqref{final2} that for $\varepsilon>0$ small, we get that
$$\sup_{\partial \widetilde U}\mathcal{\widetilde J}_\varepsilon<\sup_{\mathcal{\widetilde D}_k} \mathcal{\widetilde J}_\varepsilon\,.$$
Then $\mathcal{\widetilde J}_\varepsilon$ achieves its maximum value in $\widetilde U$ at some interior point $(\bgl t_\varepsilon,\bgl \xi_\varepsilon)$ of $\widetilde U$, which is a critical point of $\mathcal{J}_\varepsilon$ in $\widetilde U$. As already observed, we then get that {$u_{k,\varepsilon}=W_{\varepsilon,\bgl t_\varepsilon, \bgl \xi_\varepsilon}+\phi_{\varepsilon,\bgl t_\varepsilon,\bgl \xi_\varepsilon}$} is a critical point of $J_\varepsilon$, and thus, by elliptic regularity, a classical solution of \eqref{Eq1}. Up to a subsequence and taking $\widetilde U$ smaller and smaller, we can assume that $(\bgl t_\varepsilon,\bgl \xi_\varepsilon)\to (\bgl t_0, \bgl \xi_0) \in \mathcal{\widetilde D}_k$ as $\varepsilon \to 0$. Arguing as in the proof of Theorem~\ref{Th}, since $\big\|\phi_{\varepsilon,\bgl t_\varepsilon, \bgl \xi_\varepsilon} \big\|_{1,2}\to0$, by the definition of $W_{\varepsilon,\bgl t_\varepsilon,\bgl \xi_\varepsilon}$, we get that { $u_{k,\varepsilon}^{2^*} \rightharpoonup K_n^{-n}\sum_{i=1}^k \delta_{\(\xi_0\)_i}$} in the measures sense as $\varepsilon \to 0$. Then the family {$(u_{k,\varepsilon})_\varepsilon$} blows up at the points $\(\xi_0\)_1,\dots,\(\xi_0\)_k$ as $\varepsilon \to 0$, where $(\bgl t_0,\bgl \xi_0)$ is so that $\widetilde E_k(\bgl t_0,\bgl \xi_0)=m_k$. { We also get that $\lim_{k\to+\infty}\limsup_{\varepsilon\to0}\left\|\nabla u_{k,\varepsilon}\right\|_{L^2\(M\)}=+\infty$ by the definition of $W_{\varepsilon,\bgl t_\varepsilon,\bgl \xi_\varepsilon}$.} This ends the proof of Theorem~\ref{ThThThTh}.
\endproof

\smallskip
\proof[Proof of Theorem~\ref{ThThThThTh}]
We assume that $n=6$ and $h\in C^1\(M\)$. { It is not difficult to show that the expansion \eqref{Lem1Eq1m} is $C^1$--uniform with respect to $\bgl t$ and $\bgl\xi$.} A straightforward adaptation of the $C^1$--estimates in Section~\ref{Sec6} below { then} yields to the $C^1$--uniformity of the expansion for $\mathcal{J}_\varepsilon$ in Proposition~\ref{Pr3m} for all integers $k\ge1$. For any $\bgl\xi:=(\xi_1,\dots,\xi_k)\in M^k\backslash\Delta_k$, let $A_{k,\bgl\xi}$ and $E_k\(\bgl\xi\)$ be as in \eqref{matrix1}--\eqref{matrix2}, and $\bgl\xi_0:=\(\(\xi_0\)_1,\dots,\(\xi_0\)_k\)$ be an isolated critical point of $E_k$ with non-trivial local degree so that $A_{k,\bgl\xi_0}^{-1}.H$\vspace{-3pt} { has positive coordinates}. Observe that with these notations, { we can write that}
$$\widetilde{E}_k\(\bgl{t},\bgl\xi\)=c_2\big<\bgl{T},H\big(\bgl\xi\big)\big>-c_3\big<\bgl{T}, A_{k,\bgl\xi}.\bgl{T}\big>,$$
where $\bgl{T}:=\(t_1^2,\dotsc,t_k^2\)$, $H\(\bgl\xi\):=\(h\(\xi_1\),\dotsc,h\(\xi_k\)\)$ and $\<\cdot,\cdot\>$ is the Euclidean scalar product. Arguing as in the proof of Theorem~\ref{ThThTh}, it suffices to find $\bgl{t}_0\(\bgl\xi_0\)\in\(0,\infty\)^k$ such that $\(\bgl{t}_0\(\bgl\xi_0\),\bgl\xi_0\)$ is a $C^1$--stable critical point of $\widetilde{E}_k$. { One then easily checks} that such a property is achieved when taking
$$\bgl{T}_0\(\bgl\xi_0\):=\frac{c_2}{2c_3}A_{k,\bgl\xi_0}^{-1}.H\big(\bgl\xi_0\big)\quad\text{with }\bgl{T}_0\(\bgl\xi_0\)=\(\(t_0\(\bgl\xi_0\)\)_1^2,\dotsc,\(t_0\(\bgl\xi_0\)\)_k^2\),$$
which is well defined since $A_{k,\bgl\xi_0}^{-1}.H\big(\bgl\xi_0\big)$\vspace{-3pt} { has positive coordinates}. This ends the proof of Theorem~\ref{ThThThThTh}.
\endproof

\section{First derivatives estimates}\label{Sec6}

This section is devoted to the end of the proof of Proposition~\ref{Pr3}. We assume that $h\in C^1\(M\)$ and we prove the $C^1$--uniformity of the expansion for $\mathcal{J}_\varepsilon$ in Proposition~\ref{Pr3}. Arguing as in the proof of Lemma~\ref{Lem1}, it is not difficult to show that \eqref{Lem1Eq1}--\eqref{Lem1Eq4} are $C^1$--uniform with respect to $\xi\in M$ and $t$ in compact subsets of $(0,\infty)$ as $\varepsilon\to0$. We only need to prove the $C^1$--uniformity of \eqref{Eq26}. We begin with proving the preliminary Lemmas~\ref{Lem2}--\ref{Lem3quater}. Throughout this section, we identify the tangent spaces $T_\xi M$ with $\mathbb{R}^n$ thanks to local, smooth, orthonormal frames, so that $\exp_\xi$ denotes the composition of the standard exponential map (with respect to $g_\xi$) with a linear isometry $\varUpsilon_\xi:\mathbb{R}^n\to T_\xi M$ which is smooth with respect to $\xi$. We denote by $\varOmega$ the domain in $M$ where the frame is defined. We use the notations
\begin{equation}\label{Eq25}
Z_{0,\varepsilon,t,\xi}:=Z_{\varepsilon,t,\xi}\quad\text{and}\quad Z_{i,\varepsilon,t,\xi}:=Z_{\varepsilon,t,\xi,e_i}
\end{equation}
for all $i=1,\dotsc,n$, where $e_i$ is the $i$--th vector in the canonical basis of $\mathbb{R}^n$, $Z_{\varepsilon,t,\xi}$ and $Z_{\varepsilon,t,\xi,e_i}$ are as in \eqref{Eq14}--\eqref{Eq15}. We let $W_{\varepsilon,t,\xi}$ be as in \eqref{Eq9}, $J_\varepsilon$ be as in \eqref{Eq22}, and $\mathcal{J}_\varepsilon$ be as in \eqref{Eq23}. All our estimates in this section are uniform with respect to $t\in\[a,b\]$, $\xi\in\varOmega$, and $\varepsilon\in\(0,\varepsilon_0\)$ for some fixed $\varepsilon_0>0$. In Lemma~\ref{Lem2} below, we approximate the first derivatives of $W_{\varepsilon,t,\xi}\(x\)$ with respect to $t$, $\xi$, and $x$.

\begin{lemma}\label{Lem2}
There hold
\begin{align}
\bigg\|\frac{d}{dt}W_{\varepsilon,t,\xi}-\frac{n^{\frac{n-2}{4}}\(n-2\)^{\frac{n+2}{4}}\delta'_\varepsilon\(t\)}{2\delta_\varepsilon\(t\)}Z_{0,\varepsilon,t,\xi}\bigg\|_{2^*}&=\operatorname{o}\(1\),\label{Lem2Eq1}\\
\bigg\|\frac{d}{d\eta_i}W_{\varepsilon,t,\exp_\xi\eta}\Big|_{\eta=0}-\frac{n^{\frac{n-2}{4}}\(n-2\)^{\frac{n+2}{4}}}{\delta_\varepsilon\(t\)}Z_{i,\varepsilon,t,\xi}\bigg\|_{2^*}&=\operatorname{o}\(1\)\label{Lem2Eq2}
\end{align}
as $\varepsilon\to0$ for all $i=1,\dotsc,n$, where $\delta_\varepsilon\(t\)$ is as in \eqref{Eq8}.
\end{lemma}

\proof
We begin with proving \eqref{Lem2Eq1}. For any $x\in M$, we find that
\begin{equation}\label{Lem2Eq4}
\frac{d}{dt}W_{\varepsilon,t,\xi}\(x\)=G_g\(x,\xi\)\frac{d}{dt}\widehat{W}_{\varepsilon,t,\xi}\(x\),
\end{equation}
where
\begin{align}
\frac{d}{dt}\widehat{W}_{\varepsilon,t,\xi}\(x\)&=\frac{n^{\frac{n-2}{4}}\(n-2\)^{\frac{n+2}{4}}\delta'_\varepsilon\(t\)}{2\delta_\varepsilon\(t\)^{\frac{n}{2}}}\beta_n\times\left\{\begin{aligned}
&d_{g_\xi}\(x,\xi\)^{n-2}V_0\(\delta_\varepsilon\(t\)^{-1}\exp_\xi^{-1}x\)&&\text{if }d_{g_\xi}\(x,\xi\)\le r_0\\
&r_0^{n-2}V_0\(\delta_\varepsilon\(t\)^{-1}r_0\)&&\text{if }d_{g_\xi}\(x,\xi\)>r_0
\end{aligned}\right.\nonumber\\
&=\frac{n^{\frac{n-2}{4}}\(n-2\)^{\frac{n+2}{4}}\delta'_\varepsilon\(t\)}{2\delta_\varepsilon\(t\)}{\widehat{Z}_{0,\varepsilon,t,\xi}\(x\)}+\operatorname{O}\(\delta_\varepsilon\(t\)^{\frac{n-4}{2}}\delta'_\varepsilon\(t\)\mathbf{1}_{M\backslash B_\xi\(r_0/2\)}\(x\)\).\label{Lem2Eq5}
\end{align}
\eqref{Lem2Eq1} follows from \eqref{Eq8}, \eqref{Lem2Eq4}, and \eqref{Lem2Eq5}. Now, we prove \eqref{Lem2Eq2}. For any $x\in B_\xi\(r_0\)$, we get that
\begin{equation}\label{Lem2Eq6}
\frac{d}{d\eta_i}W_{\varepsilon,t,\exp_\xi\eta}\(x\)\Big|_{\eta=0}=G_g\(x,\xi\)\frac{d}{d\eta_i}\widehat{W}_{\varepsilon,t,\exp_\xi\eta}\(x\)\Big|_{\eta=0}+\frac{d}{d\eta_i}G_g\(x,\exp_\xi\eta\)\Big|_{\eta=0}\widehat{W}_{\varepsilon,t,\xi}\(x\).
\end{equation}
Moreover, letting $y=\exp_\xi^{-1}x$ and using Lemma~\ref{Lem5} in appendix, we find that
\begin{equation}\label{Lem2Eq7}
\frac{d}{d\eta_i}\widehat{W}_{\varepsilon,t,\exp_\xi\eta}\(\exp_\xi y\)\Big|_{\eta=0}=-n^{\frac{n-2}{4}}\(n-2\)^{\frac{n+2}{4}}\beta_n\frac{\delta_\varepsilon\(t\)^{\frac{n+2}{2}}\left|y\right|^{n-4}y_i}{\(\delta_\varepsilon\(t\)^2+\left|y\right|^2\)^{\frac{n}{2}}}+\operatorname{O}\(\frac{\delta_\varepsilon\(t\)^{\frac{n+2}{2}}\left|y\right|^{n-1}}{\(\delta_\varepsilon\(t\)^2+\left|y\right|^2\)^{\frac{n}{2}}}\)
\end{equation}
and, using \eqref{Eq7}, we get that
\begin{align}\label{Lem2Eq8}
\frac{d}{d\eta_i}G_g\(x,\exp_\xi\eta\)\Big|_{\eta=0}=\varLambda_\xi\(x\)\frac{d}{d\eta_i}G_{g_{\exp_\xi\eta}}\(x,\exp_\xi\eta\)\Big|_{\eta=0}+\frac{d}{d\eta_i}\varLambda_{\exp_\xi\eta}\(x\)\Big|_{\eta=0}G_{g_\xi}\(x,\xi\).
\end{align}
Since we have chosen $\varLambda_\xi$ so that $\varLambda_\xi\(\xi\)=1$ and $\nabla\varLambda_\xi\(\xi\)=0$, we get that
\begin{equation}\label{Lem2Eq9}
\varLambda_\xi\(\exp_\xi y\)=1+\operatorname{O}\(\left|y\right|^2\)\quad\text{and}\quad\frac{d}{d\eta_i}\varLambda_{\exp_\xi\eta}\(\exp_\xi y\)\Big|_{\eta=0}=\operatorname{O}\(\left|y\right|\).
\end{equation}
By \eqref{Lem2Eq8}, \eqref{Lem2Eq9}, Lemmas~\ref{Lem4} and~\ref{Lem5}, we get that
\begin{align}
G_g\(\exp_\xi y, \xi\)&=\beta_n^{-1}\left|y\right|^{2-n}+\operatorname{O}\(\left|y\right|^{4-n}\),\label{Lem2Eq10}\\
\frac{d}{d\eta_i}G_g \(\exp_\xi y,\exp_\xi\eta\)\Big|_{\eta=0}&=\(n-2\)\beta_n^{-1}\left|y\right|^{-n}y_i+\operatorname{O}\(\left|y\right|^{3-n}\).\label{Lem2Eq11}
\end{align}
Moreover, using Lemma~\ref{Lem4}, we find that
\begin{equation}\label{Lem2Eq12}
Z_{i,\varepsilon,t,\xi}\(\exp_\xi y\)=\frac{\delta_\varepsilon\(t\)^{\frac{n}{2}}{\chi(|y|)}y_i}{\(\delta_\varepsilon\(t\)^2+\left|y\right|^2\)^{\frac{n}{2}}}+\operatorname{O}\(\frac{\delta_\varepsilon\(t\)^{\frac{n}{2}}\left|y\right|^3}{\(\delta_\varepsilon\(t\)^2+\left|y\right|^2\)^{\frac{n}{2}}}\).
\end{equation}
By \eqref{Lem2Eq6}--\eqref{Lem2Eq12}, we get that
\begin{equation}\label{Lem2Eq13}
\frac{d}{d\eta_i}W_{\varepsilon,t,\exp_\xi\eta}\(\exp_\xi y\)\Big|_{\eta=0}=\frac{n^{\frac{n-2}{4}}\(n-2\)^{\frac{n+2}{4}}}{\delta_\varepsilon\(t\)}Z_{i,\varepsilon,t,\xi}\(\exp_\xi y\)+\operatorname{O}\(\frac{\delta_\varepsilon\(t\)^{\frac{n-2}{2}}\left|y\right|}{\(\delta_\varepsilon\(t\)^2+\left|y\right|^2\)^{\frac{n-2}{2}}}\)
\end{equation}
uniformly with respect to $y\in B_0\(r_0\)$. It follows from \eqref{Lem2Eq13} that
\begin{multline}\label{Lem2Eq14}
\int_{B_\xi\(r_0\)}\bigg|\frac{d}{d\eta_i}W_{\varepsilon,t,\exp_\xi\eta}\(x\)\Big|_{\eta=0}-\frac{n^{\frac{n-2}{4}}\(n-2\)^{\frac{n+2}{4}}}{\delta_\varepsilon\(t\)}Z_{i,\varepsilon,t,\xi}\(x\)\bigg|^{2^*}dv_g\\
=\operatorname{O}\(\delta_\varepsilon\(t\)^n\int_0^{r_0}\frac{r^{2^*+n-1}dr}{\(\delta_\varepsilon\(t\)^2+r^2\)^n}\)=\operatorname{o}\(1\)
\end{multline}
as $\varepsilon\to0$. In $M\backslash\overline{B_\xi\(r_0\)}$, we find that
\begin{equation}\label{Lem2Eq15}
\frac{d}{d\eta_i}W_{\varepsilon,t,\exp_\xi\eta}\(x\)\Big|_{\eta=0}=\beta_n\delta_\varepsilon\(t\)^{\frac{2-n}{2}}r_0^{n-2}U\(\delta_\varepsilon\(t\)^{-1}r_0\)\frac{d}{d\eta_i}G_g\(x,\exp_\xi\eta\)\Big|_{\eta=0}=\operatorname{O}\(\delta_\varepsilon\(t\)^{\frac{n-2}{2}}\)
\end{equation}
uniformly with respect to $x\in M\backslash\overline{B_\xi\(r_0\)}$. Finally, \eqref{Lem2Eq2} follows from \eqref{Lem2Eq14} and \eqref{Lem2Eq15}. This ends the proof of Lemma~\ref{Lem2}.
\endproof

In Lemma~\ref{Lem3} below, we approximate the first derivatives of the energy of our test functions.

\begin{lemma}\label{Lem3}
There hold
\begin{align}
&\frac{d}{dt}J_\varepsilon\(W_{\varepsilon,t,\xi}\)=\frac{n^{\frac{n-2}{4}}\(n-2\)^{\frac{n+2}{4}}\delta'_\varepsilon\(t\)}{2\delta_\varepsilon\(t\)}DJ_\varepsilon\(W_{\varepsilon,t,\xi}\).Z_{0,\varepsilon,t,\xi}+\left\{\begin{aligned}&\operatorname{o}\(\delta_\varepsilon\(t\)^2\ln\delta_\varepsilon\(t\)\)&&\text{if }n=4\\&\operatorname{o}\(\varepsilon\delta_\varepsilon\(t\)^2\)&&\text{if }n\ge5\end{aligned}\right.\label{Lem3Eq1}\\
&\frac{d}{d\eta_i}J_\varepsilon\big(W_{\varepsilon,t,\exp_\xi\eta}\big)\Big|_{\eta=0}=\frac{n^{\frac{n-2}{4}}\(n-2\)^{\frac{n+2}{4}}}{\delta_\varepsilon\(t\)}DJ_\varepsilon\(W_{\varepsilon,t,\xi}\).Z_{i,\varepsilon,t,\xi}+\left\{\begin{aligned}&\operatorname{o}\(\delta_\varepsilon\(t\)^2\)&&\text{if }n=4\\&\operatorname{o}\(\varepsilon\delta_\varepsilon\(t\)^2\)&&\text{if }n\ge5\end{aligned}\right.\label{Lem3Eq2}
\end{align}
as $\varepsilon\to0$ for all $i=1,\dotsc,n$, where $\delta_\varepsilon\(t\)$ is as in \eqref{Eq8}.
\end{lemma}

\proof
We begin with proving \eqref{Lem3Eq1}. Integration by parts gives that
\begin{multline}\label{Lem3Eq3}
\frac{d}{dt}J_\varepsilon\(W_{\varepsilon,t,\xi}\)-\frac{n^{\frac{n-2}{4}}\(n-2\)^{\frac{n+2}{4}}\delta'_\varepsilon\(t\)}{2\delta_\varepsilon\(t\)}DJ_\varepsilon\(W_{\varepsilon,t,\xi}\).Z_{0,\varepsilon,t,\xi}\\
=\int_M\(L_gW_{\varepsilon,t,\xi}-f_\varepsilon\(W_{\varepsilon,t,\xi}\)\)\bigg(\frac{d}{dt}W_{\varepsilon,t,\xi}-\frac{n^{\frac{n-2}{4}}\(n-2\)^{\frac{n+2}{4}}\delta'_\varepsilon\(t\)}{2\delta_\varepsilon\(t\)}Z_{0,\varepsilon,t,\xi}\bigg)dv_g\\
{+\int_{\partial B_\xi\(r_0\)}\(\partial_{\nu_{\inward}}W_{\varepsilon,t,\xi}+\partial_{\nu_{\outward}}W_{\varepsilon,t,\xi}\)\frac{d}{dt}W_{\varepsilon,t,\xi}\,d\sigma_g\,. }
\end{multline}
By  \eqref{Pr2Eq5bis}, \eqref{Pr2Eq5}, \eqref{Pr2Eq8}, and \eqref{Lem2Eq5}, we get that
\begin{multline}\label{Lem3Eq4}
\(L_gW_{\varepsilon,t,\xi}\(x\)-f_\varepsilon\(W_{\varepsilon,t,\xi}\(x\)\)\)\bigg(\frac{d}{dt}W_{\varepsilon,t,\xi}\(x\)-\frac{n^{\frac{n-2}{4}}\(n-2\)^{\frac{n+2}{4}}\delta'_\varepsilon\(t\)}{2\delta_\varepsilon\(t\)}Z_{0,\varepsilon,t,\xi}\(x\)\bigg)\\
=\operatorname{O}\(\(\delta_\varepsilon\(t\)^{n-1}\delta'_\varepsilon\(t\)+\varepsilon\delta_\varepsilon\(t\)^{n-3}\delta'_\varepsilon\(t\)\)\mathbf{1}_{M\backslash B_\xi\(r_0/2\)}\(x\)\).
\end{multline}
\eqref{Lem3Eq1} follows from \eqref{Lem3Eq3} and \eqref{Lem3Eq4} in view of \eqref{Pr2Eq4} and \eqref{Lem2Eq5}. Now, we prove \eqref{Lem3Eq2}. Integration by parts gives that
\begin{multline}\label{Lem3Eq5}
\frac{d}{d\eta_i}J_\varepsilon\big(W_{\varepsilon,t,\exp_\xi\eta}\big)\Big|_{\eta=0}-\frac{n^{\frac{n-2}{4}}\(n-2\)^{\frac{n+2}{4}}}{\delta_\varepsilon\(t\)}DJ_\varepsilon\(W_{\varepsilon,t,\xi}\).Z_{i,\varepsilon,t,\xi}\\
=\int_M\(L_gW_{\varepsilon,t,\xi}-f_\varepsilon\(W_{\varepsilon,t,\xi}\)\)\bigg(\frac{d}{d\eta_i}W_{\varepsilon,t,\exp_\xi\eta}\Big|_{\eta=0}-\frac{n^{\frac{n-2}{4}}\(n-2\)^{\frac{n+2}{4}}}{\delta_\varepsilon\(t\)}Z_{i,\varepsilon,t,\xi}\bigg)dv_g\\
+\int_{\partial B_\xi\(r_0\)}\[\partial_{\nu_{\inward}}W_{\varepsilon,t,\xi} \(\frac{d}{d\eta_i}W_{\varepsilon,t,\exp_\xi\eta}\Big|_{\eta=0}\)_{\inward} +\partial_{\nu_{\outward}}W_{\varepsilon,t,\xi} \(\frac{d}{d\eta_i}W_{\varepsilon,t,\exp_\xi\eta}\Big|_{\eta=0}\)_{\outward}\] d\sigma_g\,,
\end{multline}
where
\begin{align*}
\(\frac{d}{d\eta_i}W_{\varepsilon,t,\exp_\xi\eta}\Big|_{\eta=0}\)_{\inward} \(x\)&=G_g\(x,\xi\)\frac{d}{d\eta_i}\widehat{W}_{\varepsilon,t,\exp_\xi\eta}\(x\)\Big|_{\eta=0}+\frac{d}{d\eta_i}G_g\(x,\exp_\xi\eta\)\Big|_{\eta=0}\widehat{W}_{\varepsilon,t,\xi}\(x\),\\
\(\frac{d}{d\eta_i}W_{\varepsilon,t,\exp_\xi\eta}\Big|_{\eta=0}\)_{\outward} \(x\)&=\frac{d}{d\eta_i}G_g\(x,\exp_\xi\eta\)\Big|_{\eta=0} \widehat{W}_{\varepsilon,t,\xi}\(x\)
\end{align*}
for all $x\in\partial B_\xi\(r_0\)$, in view of \eqref{Lem2Eq6} and \eqref{Lem2Eq15}. Regarding the second term in the right hand side of \eqref{Lem3Eq5}, on $\partial B_\xi\(r_0\)$, we find that
\begin{align}
&\partial_{\nu_{\inward}}W_{\varepsilon,t,\xi} \(\frac{d}{d\eta_i}W_{\varepsilon,t,\exp_\xi\eta}\Big|_{\eta=0}\)_{\inward}+\partial_{\nu_{\outward}}W_{\varepsilon,t,\xi}\(\frac{d}{d\eta_i}W_{\varepsilon,t,\exp_\xi\eta}\Big|_{\eta=0}\)_{\outward}\nonumber\\
&=\beta_nG_g\(\cdot,\xi\)\delta_\varepsilon\(t\)^{\frac{2-n}{2}}\frac{d}{dr}\(r^{n-2}U\(\delta_\varepsilon(t)^{-1}r\) \)\Big|_{r=r_0}\(\frac{d}{d\eta_i}W_{\varepsilon,t,\exp_\xi\eta}\Big|_{\eta=0}\)_{\inward} \nonumber \\
&+\beta_n \partial_{\nu_{\inward}}G_g\(\cdot,\xi\) \delta_\varepsilon\(t\)^{\frac{2-n}{2}}r_0^{n-2}U\(\delta_\varepsilon(t)^{-1}r_0\)  G_g\(x,\xi\)\frac{d}{d\eta_i}\widehat{W}_{\varepsilon,t,\exp_\xi\eta}\(x\)\Big|_{\eta=0}
\label{Lem3Eq6}
\end{align}
for all $x\in\partial B_\xi\(r_0\)$. By \eqref{Lem2Eq13}, we get that
\begin{equation}\label{Lem3Eq7}
\( \frac{d}{d\eta_i}W_{\varepsilon,t,\exp_\xi\eta}\Big|_{\eta=0}\)_{\inward}=\operatorname{O}\(\delta_\varepsilon\(t\)^{\frac{n-2}{2}}\).
\end{equation}
If follows from \eqref{Lem3Eq6}--\eqref{Lem3Eq7} that
\begin{equation}\label{Lem3Eq8}
\partial_{\nu_{\inward}}W_{\varepsilon,t,\xi} \(\frac{d}{d\eta_i}W_{\varepsilon,t,\exp_\xi\eta}\Big|_{\eta=0}\)_{\inward}+\partial_{\nu_{\outward}}W_{\varepsilon,t,\xi} \(\frac{d}{d\eta_i}W_{\varepsilon,t,\exp_\xi\eta}\Big|_{\eta=0}\)_{\outward}=\operatorname{O}\(\delta_\varepsilon\(t\)^n\)
\end{equation}
in view of \eqref{Pr2Eq4} and \eqref{Lem2Eq7}. Now, we estimate the first term in the right hand side of \eqref{Lem3Eq5}. In $B_\xi\(r_0\)$, using \eqref{Pr2Eq5}, \eqref{Pr2Eq8}, and \eqref{Lem2Eq13}, we find that
\begin{align}
&\int_{B_\xi\(r_0\)}\(L_gW_{\varepsilon,t,\xi}\(x\)-W_{\varepsilon,t,\xi}\(x\)^{2^*-1}\)\bigg(\frac{d}{d\eta_i}W_{\varepsilon,t,\exp_\xi\eta}\(x\)\Big|_{\eta=0}-\frac{n^{\frac{n-2}{4}}\(n-2\)^{\frac{n+2}{4}}}{\delta_\varepsilon\(t\)}Z_{i,\varepsilon,t,\xi}\(x\)\bigg)dv_g\nonumber\\
&\qquad=\delta_\varepsilon\(t\)^n\int_0^{r_0}\frac{r^{n-2}dr}{\(\delta_\varepsilon\(t\)^2+r^2\)^{n-1}}\times\left\{\begin{aligned}
&\operatorname{O}\(r^{n-2}\)&&\text{if }n=4,5\text{ or }\(M,g\)\text{ l.c.f.}\\
&\operatorname{O}\(r^4\ln r\)&&\text{if }n=6\text{ and }\(M,g\)\text{ non-l.c.f.}\\
&\operatorname{O}\(r^4\)&&\text{if }n\ge7\text{ and }\(M,g\)\text{ non-l.c.f.}
\end{aligned}\right.\nonumber\allowdisplaybreaks\\
&\qquad=\left\{\begin{aligned}
&\operatorname{O}\(\delta_\varepsilon\(t\)^{n-1}\)&&\text{if }n=4,5\text{ or }\(M,g\)\text{ l.c.f.}\\
&\operatorname{O}\(\delta_\varepsilon\(t\)^5\ln\delta_\varepsilon\(t\)\)&&\text{if }n=6\text{ and }\(M,g\)\text{ non-l.c.f.}\\
&\operatorname{O}\(\delta_\varepsilon\(t\)^5\)&&\text{if }n\ge7\text{ and }\(M,g\)\text{ non-l.c.f.}
\end{aligned}\right.\label{Lem3Eq9}
\end{align}
and, using \eqref{Lem2Eq13}, we find that
\begin{multline}\label{Lem3Eq10}
\int_{B_\xi\(r_0\)}hW_{\varepsilon,t,\xi}\(x\)\bigg(\frac{d}{d\eta_i}W_{\varepsilon,t,\exp_\xi\eta}\(x\)\Big|_{\eta=0}-\frac{n^{\frac{n-2}{4}}\(n-2\)^{\frac{n+2}{4}}}{\delta_\varepsilon\(t\)}Z_{i,\varepsilon,t,\xi}\(x\)\bigg)dv_g\\
=\operatorname{O}\(\delta_\varepsilon\(t\)^{n-2}\int_0^{r_0}\frac{r^ndr}{\(\delta_\varepsilon\(t\)^2+r^2\)^{n-2}}\)=\left\{\begin{aligned}
&\operatorname{O}\(\delta_\varepsilon\(t\)^2\)&&\text{if }n=4\\
&\operatorname{O}\(\delta_\varepsilon\(t\)^3\left|\ln\delta_\varepsilon\(t\)\right|\)&&\text{if }n=5\\
&\operatorname{O}\(\delta_\varepsilon\(t\)^3\)&&\text{if }n\ge6.
\end{aligned}\right.
\end{multline}
In $M\backslash\overline{B_\xi\(r_0\)}$, using \eqref{Pr2Eq5bis} and \eqref{Lem2Eq15}, we find that
\begin{equation}\label{Lem3Eq11}
\(L_gW_{\varepsilon,t,\xi}-f_\varepsilon\(W_{\varepsilon,t,\xi}\)\)\bigg(\frac{d}{d\eta_i}W_{\varepsilon,t,\exp_\xi\eta}\Big|_{\eta=0}-\frac{n^{\frac{n-2}{4}}\(n-2\)^{\frac{n+2}{4}}}{\delta_\varepsilon\(t\)}Z_{i,\varepsilon,t,\xi}\bigg)=\operatorname{O}\(\delta_\varepsilon\(t\)^n+\varepsilon\delta_\varepsilon\(t\)^{n-2}\)
\end{equation}
uniformly with respect to $x\in M\backslash\overline{B_\xi\(r_0\)}$. Finally, \eqref{Lem3Eq2} follows from \eqref{Lem3Eq5}--\eqref{Lem3Eq11}. This ends the proof of Lemma~\ref{Lem3}.
\endproof

Now, we prove the following error estimates.

\begin{lemma}\label{Lem3bis}
There exists a positive constant $C=C\(a,b,n,M,g,h\)$ such that
\begin{align}
&\left\|L_g^{-1}\(f'_\varepsilon(W_{\varepsilon,t,\xi})Z_{0,\varepsilon,t,\xi}\)-Z_{0,\varepsilon,t,\xi} \right\|_{1,2}\nonumber\\
&\qquad\le C\left\{\begin{aligned}
&\delta_\varepsilon(t)^{\frac{n-2}{2}}&&\text{if }4\le n\le9\text{ or }\(M,g\)\text{ l.c.f.}\\
&\delta_\varepsilon(t)^4&&\text{if }n\ge10\text{ and }\(M,g\)\text{ non-l.c.f.}
\end{aligned}\right.\label{Lem3bisEq1}\\
&\left\|L_g^{-1}\(f'_\varepsilon(W_{\varepsilon,t,\xi})Z_{i,\varepsilon,t,\xi}\)-Z_{i,\varepsilon,t,\xi}\right\|_{1,2}\nonumber\\
&\qquad\le C\left\{\begin{aligned}
&\delta_\varepsilon(t)^2\left|\ln\delta_\varepsilon(t)\right|^{\frac{3}{4}}&&\text{if }n=4\\
&\delta_\varepsilon(t)^{\frac{n}{2}}&&\text{if }n=5\text{ or (}n\ge6\text{ and }\(M,g\)\text{ l.c.f.)} \\
&\delta_\varepsilon(t)^3\left|\ln\delta_\varepsilon(t)\right|^{\frac{2}{3}}&&\text{if }n=6\text{ and }\(M,g\)\text{ non-l.c.f.}\\
&\delta_\varepsilon(t)^3&&\text{if }n\ge7\text{ and }\(M,g\)\text{ non-l.c.f.}
\end{aligned}\right.\label{Lem3bisEq2}
\end{align}
for all $i=1,\dotsc,n$, where $\delta_\varepsilon\(t\)$ is as in \eqref{Eq8}.
\end{lemma}

\proof
For any $i=0,\dotsc,n$ and $\phi\in H^2_1\(M\)$, an integration by parts gives
\begin{equation}\label{Lem3bisEq3}
\<L_g^{-1}\(f'_\varepsilon (W_{\varepsilon,t,\xi})Z_{i,\varepsilon,t,\xi}\)-Z_{i,\varepsilon,t,\xi},\phi\>_{L_g}=\int_M\(f'_\varepsilon\(W_{\varepsilon,t,\xi}\)Z_{i,\varepsilon,t,\xi}-L_gZ_{i,\varepsilon,t,\xi}\)\phi dv_g\,.
\end{equation}
By Sobolev's embedding $H_1^2(M)\hookrightarrow L^{2^*}(M)$, it follows from \eqref{Lem3bisEq3} that
\begin{equation}\label{Lem3bisEq4}
\left\|L_g^{-1}\(f'_\varepsilon(W_{\varepsilon,t,\xi})Z_{i,\varepsilon,t,\xi}\)-Z_{i,\varepsilon,t,\xi}\right\|_{1,2}=\operatorname{O}\(\left\|f'_\varepsilon\(W_{\varepsilon,t,\xi}\)Z_{i,\varepsilon,t,\xi}-L_gZ_{i,\varepsilon,t,\xi}\right\|_{\frac{2n}{n+2}}\).
\end{equation}
By conformal covariance \eqref{confinv} of $L_g$ and by \eqref{Eq7}, in $B_\xi\(r_0\)$, we can write that
\begin{multline}\label{Lem3bisEq5}
f'_\varepsilon\(W_{\varepsilon,t,\xi}\)Z_{i,\varepsilon,t,\xi}-L_gZ_{i,\varepsilon,t,\xi}\\
=\varLambda_\xi^{2^*-1}\big[\(2^*-1\)G_{g_\xi}\(\cdot,\xi\)^{2^*-1}\widehat{W}_{\varepsilon,t,\xi}^{2^*-2}\widehat{Z}_{i,\varepsilon,t,\xi}-L_{g_\xi}\big(G_{g_\xi}\(\cdot,\xi\)\widehat{Z}_{i,\varepsilon,t,\xi}\big)\big]-\varepsilon hZ_{i,\varepsilon,t,\xi}\,.
\end{multline}
Since $\widehat{Z}_{i,\varepsilon,t,\xi}\(\xi\)=0$ and $L_{g_\xi}G_{g_\xi}\(\cdot,\xi\)=\delta_\xi$, we get that
\begin{equation}\label{Lem3bisEq6}
L_{g_\xi}\big(G_{g_\xi}\(\cdot,\xi\)\widehat{Z}_{i,\varepsilon,t,\xi}\big)=G_{g_\xi}\(\cdot,\xi\)\Delta_{g_\xi}\widehat{Z}_{i,\varepsilon,t,\xi}-2\big<\nabla G_{g_\xi}\(\cdot,\xi\),\nabla\widehat{Z}_{i,\varepsilon,t,\xi}\big>_{g_\xi}\,.
\end{equation}
We begin with considering the case $i=0$. Since $\widehat{Z}_{0,\varepsilon,t,\xi}\circ\exp_\xi$ is radially symmetrical and $V_0$ is a solution to the equation $\Delta_{\Eucl}V_0=\(2^*-1\)U^{2^*-2}V_0$ in $\mathbb{R}^n$, writing $\Delta_{g_\xi}\widehat{Z}_{0,\varepsilon,t,\xi}\(\exp_\xi y\)$ in polar coordinates and using \eqref{Eq6}, we find that
\begin{eqnarray}
&&G_{g_\xi}(\exp_\xi y,\xi) \Delta_{g_\xi}\widehat{Z}_{0,\varepsilon,t,\xi}(\exp_\xi y)\nonumber\\
&&\quad=G_{g_\xi}(\exp_\xi y,\xi) \Delta_{\Eucl}(\widehat{Z}_{0,\varepsilon,t,\xi}\circ \exp_\xi)(y)+\operatorname{O}\(|y|^{N-n+1}|\nabla(\widehat{Z}_{0,\varepsilon,t,\xi}\circ\exp_\xi)(y)|\)\nonumber\allowdisplaybreaks\\
&&\quad= (2^*-1)\beta_n \delta_\varepsilon(t)^{-\frac{n+2}{2}}\chi |y|^{n-2} G_{g_\xi}(\exp_\xi y,\xi) U(\frac{|y|}{\delta_\varepsilon(t)})^{2^*-2}V_0(\frac{y}{\delta_\varepsilon(t)})\nonumber\\
&&\qquad-2(n-2)\beta_n \delta_\varepsilon(t)^{\frac{n+2}{2}} \chi |y|^{n-4} G_{g_\xi}(\exp_\xi y,\xi) \frac{(n+2)|y|^2-(n-2)\delta_\varepsilon(t)^2}{(\delta_\varepsilon(t)^2+|y|^2)^{\frac{n+2}{2}}}\nonumber\\
&&\qquad+\operatorname{O}\(\frac{\delta_\varepsilon(t)^{\frac{n+2}{2}}|y|^{N-2}}{(\delta_\varepsilon(t)^2+|y|^2)^{\frac{n}{2}}}
+\delta_\varepsilon\(t\)^{\frac{n-2}{2}}\)\label{Lem3bisEq7}
\end{eqnarray}
uniformly with respect to $y\in B_0\(r_0\)$. Moreover, since $\widehat{Z}_{0,\varepsilon,t,\xi}\circ\exp_\xi$ is radially symmetrical, we get that
\begin{eqnarray}\label{Lem3bisEq8}
&&\big<\nabla G_{g_\xi}\(\exp_\xi y,\xi\),\nabla\widehat{Z}_{0,\varepsilon,t,\xi}\(\exp_\xi y\)\big>_{g_\xi}=
\partial_r\[G_{g_\xi}\(\exp_\xi y,\xi\)\]\partial_r\big[\widehat{Z}_{0,\varepsilon,t,\xi}\circ\exp_\xi\big]\\
&&\qquad=\beta_n \delta_\varepsilon(t)^{\frac{n+2}{2}} \chi |y|^{n-3} \partial_r \[G_{g_\xi}(\exp_\xi y,\xi)\] \frac{(n+2)|y|^2-(n-2)\delta_\varepsilon(t)^2}{(\delta_\varepsilon(t)^2+|y|^2)^{\frac{n+2}{2}}}+\operatorname{O}\(\delta_\varepsilon(t)^{\frac{n-2}{2}}\).\nonumber
\end{eqnarray}
By \eqref{Lem3bisEq6}--\eqref{Lem3bisEq8} and Lemma~\ref{Lem4}, we get that
\begin{multline*}
\(2^*-1\)G_{g_\xi}\(\exp_\xi y,\xi\)^{2^*-1}\widehat{W}_{\varepsilon,t,\xi}\(\exp_\xi y\)^{2^*-2}\widehat{Z}_{0,\varepsilon,t,\xi}\(\exp_\xi y\)-L_{g_\xi}\big(G_{g_\xi}\(\cdot,\xi\)\widehat{Z}_{0,\varepsilon,t,\xi}\big)\(\exp_\xi y\)\\
=\operatorname{O}\(\delta_\varepsilon(t)^{\frac{n-2}{2}}\)+\frac{\delta_\varepsilon(t)^{\frac{n+2}{2}}}{\(\delta_\varepsilon(t)^2+\left|y\right|^2\)^{\frac{n}{2}}}\times\left\{\begin{aligned}
&\operatorname{O}\(|y|^{n-4}\)&&\text{if }n=4,5\text{ or }\(M,g\)\text{ l.c.f.}\\
&\operatorname{O}\(|y|^2\ln|y|\) &&\text{if }n=6\text{ and }\(M,g\)\text{ non-l.c.f.}\\
&\operatorname{O}\(|y|^2\)&&\text{if }n\ge7\text{ and }\(M,g\)\text{ non-l.c.f.},
\end{aligned}\right.
\end{multline*}
which, inserted into \eqref{Lem3bisEq4}--\eqref{Lem3bisEq5}, yields to the validity of \eqref{Lem3bisEq1} in view of \eqref{Eq8}. Now, we consider the case $i=1,\dotsc,n$. Using \eqref{Eq6}, we find that
\begin{multline}\label{Lem3bisEq13}
\Delta_{g_\xi}\widehat{Z}_{i,\varepsilon,t,\xi}\(\exp_\xi y\)=\Delta_{\Eucl}\big(\widehat{Z}_{i,\varepsilon,t,\xi}\circ\exp_\xi\big)\(y\)-\beta_n\chi\(\left|y\right|\)\partial_jg_\xi^{ij}\(\exp_\xi y\)\frac{\delta_\varepsilon\(t\)^{\frac{n}{2}}\left|y\right|^{n-2}}{\(\delta_\varepsilon\(t\)^2+\left|y\right|^2\)^{\frac{n}{2}}}\\
+\operatorname{O}\Bigg(\frac{\delta_\varepsilon\(t\)^{\frac{n}{2}}\left|y\right|^{N+n-3}}{\(\delta_\varepsilon\(t\)^2+\left|y\right|^2\)^{\frac{n}{2}}}+\left|y\right|^N \left|\nabla\(\chi\(\left|y\right|\)\frac{\delta_\varepsilon\(t\)^{\frac{n}{2}}\left|y\right|^{n-2}}{\(\delta_\varepsilon\(t\)^2+\left|y\right|^2\)^{\frac{n}{2}}}\)\right|\Bigg)
\end{multline}
uniformly with respect to $y\in B_0\(r_0\)$, where $\partial_jg_\xi^{ij}$ are the derivatives of the components of $g_\xi^{-1}$ in geodesic normal coordinates. Since $g_\xi$ defines conformal normal coordinates of order $N\ge3$, see Lee--Parker~\cite{LeePar}*{Theorem~5.1 and Lemma~5.5}, we get that the Ricci curvature $\Ric_{g_\xi}$ of $g_\xi$ vanishes at $\xi$, and thus
\begin{equation}\label{Lem3bisEq14}
\partial_jg_\xi^{ij}\(\exp_\xi y\)=\left\{\begin{aligned}&0&&\text{if }\(M,g\)\text{ l.c.f.}\\&-\frac{1}{3}\(\Ric_{g_\xi}\)_{ip}\(\xi\)y^p+\operatorname{O}\(|y|^2\)=\operatorname{O}\(|y|^2\)&&\text{if }\(M,g\)\text{ non-l.c.f.}\end{aligned}\right.
\end{equation}
Since $V_i$ is a solution of $\Delta_{\Eucl}V_i=(2^*-1)U^{2^*-2}V_i$ in $\mathbb{R}^n$, by \eqref{Lem3bisEq13}--\eqref{Lem3bisEq14}, we find that
\begin{eqnarray}
&&G_{g_\xi}(\exp_\xi y,\xi) \Delta_{g_\xi}\widehat{Z}_{i,\varepsilon,t,\xi}\(\exp_\xi y\)\nonumber\\
&&\quad=(2^*-1)\beta_n \delta_\varepsilon(t)^{-\frac{n+2}{2}} \chi |y|^{n-2}G_{g_\xi}(\exp_\xi y,\xi) U(\frac{|y|}{\delta_\varepsilon(t)})^{2^*-2}V_i(\frac{y}{\delta_\varepsilon(t)})\nonumber \\
&&\qquad+2(n-2)\beta_n \delta_\varepsilon(t)^{\frac{n}{2}} \chi |y|^{n-4}y_i G_{g_\xi}(\exp_\xi y,\xi) \frac{|y|^2-(n-1)\delta_\varepsilon(t)^2}{(\delta_\varepsilon(t)^2+|y|^2)^{\frac{n+2}{2}}}\nonumber\\
&&\qquad+\operatorname{O}\(\delta_\varepsilon(t)^{\frac{n}{2}}+\frac{\delta_\varepsilon\(t\)^{\frac{n}{2}}|y|^{\alpha}}{\(\delta_\varepsilon\(t\)^2+\left|y\right|^2\)^{\frac{n}{2}}}\)\label{Lem3bisEq15}
\end{eqnarray}
where $\alpha=N-2$ if $\(M,g\)$ is l.c.f.~and $\alpha=2$ is $\(M,g\)$ is not l.c.f. Moreover, we get that
\begin{eqnarray}
&&\big<\nabla G_{g_\xi}(\exp_\xi y,\xi),\nabla\widehat{Z}_{i,\varepsilon,t,\xi}(\exp_\xi y)\big>_{g_\xi} \nonumber \\
&&=y_i\partial_r\[G_{g_\xi}(\exp_\xi y,\xi)\]\partial_r\[ \frac{\beta_n
\chi \delta_\varepsilon(t)^{\frac{n}{2}}|y|^{n-2}}{(\delta_\varepsilon(t)^2+|y|^2)^{\frac{n}{2}}}\]+\partial_{y_i}\[G_{g_\xi}\(\exp_\xi y,\xi\)\]\frac{\beta_n\chi \delta_\varepsilon\(t\)^{\frac{n}{2}}\left|y\right|^{n-2}}{\(\delta_\varepsilon\(t\)^2+\left|y\right|^2\)^{\frac{n}{2}}}\nonumber\allowdisplaybreaks\\
&&=-\beta_n\delta_\varepsilon(t)^{\frac{n}{2}}  \chi(|y|) \partial_r\[G_{g_\xi}(\exp_\xi y,\xi)\] |y|^{n-3} y_i \frac{2|y|^2-(n-2)\delta_\varepsilon(t)^2}{(\delta_\varepsilon(t)^2+|y|^2)^{\frac{n+2}{2}}}\nonumber \\
&&+\beta_n\delta_\varepsilon(t)^{\frac{n}{2}}  \chi(|y|) \partial_{y_i}\[G_{g_\xi}(\exp_\xi y,\xi)\] \frac{\delta_\varepsilon(t)^{\frac{n}{2}}|y|^{n-2}}{(\delta_\varepsilon(t)^2+|y|^2)^{\frac{n}{2}}}+\operatorname{O}\(\delta_\varepsilon(t)^{\frac{n}{2}}\).\label{Lem3bisEq16}
\end{eqnarray}
By \eqref{Lem3bisEq6}, \eqref{Lem3bisEq15}--\eqref{Lem3bisEq16} and Lemma~\ref{Lem4}, we get that
\begin{multline*}
\(2^*-1\)G_{g_\xi}\(\exp_\xi y,\xi\)^{2^*-1}\widehat{W}_{\varepsilon,t,\xi}\(\exp_\xi y\)^{2^*-2}\widehat{Z}_{i,\varepsilon,t,\xi}\(\exp_\xi y\)-L_{g_\xi}\big(G_{g_\xi}\(\cdot,\xi\)\widehat{Z}_{i,\varepsilon,t,\xi}\big)\(\exp_\xi y\)\\
=\operatorname{O}\(\delta_\varepsilon(t)^{\frac{n}{2}}\)+\frac{\delta_\varepsilon(t)^{\frac{n}{2}}}{\(\delta_\varepsilon(t)^2+\left|y\right|^2\)^{\frac{n}{2}}}\times\left\{\begin{aligned}
&\operatorname{O}\(|y|^{n-3}\)&&\text{if }n=4,5\text{ or }\(M,g\)\text{ l.c.f.}\\
&\operatorname{O}\(|y|^2\)&&\text{if }n\ge6\text{ and }\(M,g\)\text{ non-l.c.f.}
\end{aligned}\right.
\end{multline*}
which, inserted into \eqref{Lem3bisEq4}--\eqref{Lem3bisEq5}, yields to the validity of \eqref{Lem3bisEq2} in view of \eqref{Eq8}. This ends the proof of Lemma~\ref{Lem3bis}.
\endproof

By Proposition~\ref{Pr1}, for $\varepsilon$ small, for any $t\in\[a,b\]$ and $\xi\in M$, there exist $\lambda_{\varepsilon,t,\xi}\in\mathbb{R}$ and $\omega_{\varepsilon,t,\xi}\in T_\xi M$ such that
\begin{equation}\label{Eq27}
DJ_\varepsilon\(W_{\varepsilon,t,\xi}+\phi_{\varepsilon,t,\xi}\)=\<\lambda_{\varepsilon,t,\xi}Z_{\varepsilon,t,\xi}+Z_{\varepsilon,t,\xi,\omega_{\varepsilon,t,\xi}},\cdot\>_{L_g},
\end{equation}
where $Z_{\varepsilon,t,\xi}$ and $Z_{\varepsilon,t,\xi,\omega_{\varepsilon,t,\xi}}$ are as in \eqref{Eq13}. We let $Z_{0,\varepsilon,t,\xi},\dotsc,Z_{n,\varepsilon,t,\xi}\in\mathbb{R}$ be as in \eqref{Eq25}. We let $\lambda_{0,\varepsilon,t,\xi},\dotsc,\lambda_{n,\varepsilon,t,\xi}\in\mathbb{R}$ be such that
\begin{equation}\label{Eq28}
\lambda_{0,\varepsilon,t,\xi}:=\lambda_{\varepsilon,t,\xi}\quad\text{and}\quad \sum_{i=1}^n\lambda_{i,\varepsilon,t,\xi}e_i:=\omega_{\varepsilon,t,\xi}\,,
\end{equation}
where $\lambda_{\varepsilon,t,\xi}$ and $\omega_{\varepsilon,t,\xi}$ are as in \eqref{Eq27}. We estimate the $\lambda_{i,\varepsilon,t,\xi}$'s in Lemma~\ref{Lem3ter} below.

\begin{lemma}\label{Lem3ter}
For any $i=0,\dotsc,n$, in case $n\ge5$, there holds
\begin{equation}\label{Lem3terEq1}
\lambda_{i,\varepsilon,t,\xi}=\frac{DJ_\varepsilon\(W_{\varepsilon,t,\xi}\).Z_{i,\varepsilon,t,\xi}}{\left\|\nabla V_i\right\|_2^2}+\left\{\begin{aligned}
&\operatorname{o}\(\varepsilon\delta_\varepsilon(t)^2\)&&\text{if }i=0\\
&\operatorname{o}\(\varepsilon\delta_\varepsilon(t)^3\)&&\text{if }i=1,\dots,n
\end{aligned}\right.
\end{equation}
and in case $n=4$, there holds
\begin{equation}\label{Lem3terEq1bis}
\lambda_{i,\varepsilon,t,\xi}=\frac{DJ_\varepsilon\(W_{\varepsilon,t,\xi}\).Z_{0,\varepsilon,t,\xi}}{\left\|\nabla V_0\right\|_2^2}+\left\{\begin{aligned}
&\operatorname{o}\(\delta_\varepsilon(t)^2\)&&\text{if }i=0\\
&\operatorname{O}\(\varepsilon^2\delta_\varepsilon(t)^2\)&&\text{if }i=1,\dots,n
\end{aligned}\right.
\end{equation}
as $\varepsilon\to0$, where $\delta_\varepsilon\(t\)$ is as in \eqref{Eq8}. In particular, there holds
\begin{equation}\label{Lem3terEq1ter}
\lambda_{i,\varepsilon,t,\xi}=\left\{\begin{aligned}
&\operatorname{O}\(\delta_\varepsilon(t)^2 \ln |\delta_\varepsilon(t)| \)&&\text{if }i=0\hbox{ and }n=4\\
&\operatorname{O}\(\delta_\varepsilon(t)^2 \)&&\text{if }i=1,\dots,n\hbox{ and }n=4\\
&\operatorname{O}\(\varepsilon\delta_\varepsilon(t)^2\)&&\text{if }n \geq 5.
\end{aligned}\right.
\end{equation}
\end{lemma}

\proof
For any $i=0,\dotsc,n$, by \eqref{Eq27}--\eqref{Eq28}, we get that
\begin{equation}\label{Lem3terEq2}
DJ_\varepsilon\(W_{\varepsilon,t,\xi}+\phi_{\varepsilon,t,\xi}\).Z_{i,\varepsilon,t,\xi}=\sum_{j=0}^n\lambda_{j,\varepsilon,t,\xi}\<Z_{i,\varepsilon,t,\xi},Z_{j,\varepsilon,t,\xi}\>_{L_g}.
\end{equation}
For any $i,j=0,\dotsc,n$, we find that
\begin{equation}\label{Lem3terEq3}
\<Z_{i,\varepsilon,t,\xi},Z_{j,\varepsilon,t,\xi}\>_{L_g}=\left\|\nabla V_i\right\|_2^2\delta_{ij}+\operatorname{o}\(\delta_\varepsilon\(t\)\)
\end{equation}
as $\varepsilon\to0$, where the $\delta_{ij}$'s are the Kronecker symbols. It follows from \eqref{Lem3terEq2}--\eqref{Lem3terEq3} that
\begin{equation}\label{Lem3terEq4}
DJ_\varepsilon\(W_{\varepsilon,t,\xi}+\phi_{\varepsilon,t,\xi}\).Z_{i,\varepsilon,t,\xi}=\lambda_{i,\varepsilon,t,\xi}\left\|\nabla V_i\right\|_2^2+\operatorname{o}\bigg(\delta_\varepsilon\(t\)\sum_{j=0}^n\left|\lambda_{j,\varepsilon,t,\xi}\right|\bigg)
\end{equation}
as $\varepsilon\to0$. Independently, we get that
\begin{multline}\label{Lem3terEq5}
DJ_\varepsilon\(W_{\varepsilon,t,\xi}+\phi_{\varepsilon,t,\xi}\).Z_{i,\varepsilon,t,\xi}=DJ_\varepsilon\(W_{\varepsilon,t,\xi}\).Z_{i,\varepsilon,t,\xi}+\<Z_{i,\varepsilon,t,\xi}-L_g^{-1}\(f'_\varepsilon\(W_{\varepsilon,t,\xi}\)Z_{i,\varepsilon,t,\xi}\),\phi_{\varepsilon,t,\xi}\>_{L_g}\\
-\int_M\(f_0\(W_{\varepsilon,t,\xi}+\phi_{\varepsilon,t,\xi}\)-f_0\(W_{\varepsilon,t,\xi}\)-f'_0\(W_{\varepsilon,t,\xi}\)\phi_{\varepsilon,t,\xi}\)Z_{i,\varepsilon,t,\xi}dv_g\,.
\end{multline}
By Cauchy--Schwarz inequality, we get that
\begin{equation}\label{Lem3terEq6}
\left|\<Z_{i,\varepsilon,t,\xi}-L_g^{-1}\(f'_\varepsilon\(W_{\varepsilon,t,\xi}\)Z_{i,\varepsilon,t,\xi}\),\phi_{\varepsilon,t,\xi}\>_{L_g}\right|\le\left\|Z_{i,\varepsilon,t,\xi}-L_g^{-1}\(f'_\varepsilon\(W_{\varepsilon,t,\xi}\)Z_{i,\varepsilon,t,\xi}\)\right\|_{L_g}\left\|\phi_{\varepsilon,t,\xi}\right\|_{L_g}\,.
\end{equation}
By the Mean Value Theorem and H\"older's inequality, we get that
\begin{align}
&\int_M\(f_0\(W_{\varepsilon,t,\xi}+\phi_{\varepsilon,t,\xi}\)-f_0\(W_{\varepsilon,t,\xi}\)-f'_0\(W_{\varepsilon,t,\xi}\)\phi_{\varepsilon,t,\xi}\)Z_{i,\varepsilon,t,\xi}dv_g\nonumber\\
&\qquad=\left\{\begin{aligned}
& \operatorname{O}\big(\big(\left\|W_{\varepsilon,t,\xi}\right\|^{2^*-3}_{2^*}+\left\|\phi_{\varepsilon,t,\xi}\right\|^{2^*-3}_{2^*}\big)\left\|\phi_{\varepsilon,t,\xi}\right\|^2_{2^*}\left\|Z_{i,\varepsilon,t,\xi}\right\|_{2^*}\big)&&\text{if }n=4,5\\
&\operatorname{O}\big(\left\|W_{\varepsilon,t,\xi}^{2^*-3}Z_{i,\varepsilon,t,\xi}\right\|_{2^*}\left\|\phi_{\varepsilon,t,\xi}\right\|^2_{2^*}\big)&&\text{if }n\ge6.
\end{aligned}\right.\nonumber\\
&\qquad=\operatorname{O}\big(\left\|\phi_{\varepsilon,t,\xi}\right\|^2_{2^*}\big).\label{Lem3terEq7}
\end{align}
By \eqref{Eq8}, \eqref{Lem3terEq5}--\eqref{Lem3terEq7}, Propositions~\ref{Pr1},~\ref{Pr2}, and Lemma~\ref{Lem3bis}, we get that in case $n\ge5$, there holds
\begin{equation}\label{Lem3terEq8}
DJ_\varepsilon\(W_{\varepsilon,t,\xi}+\phi_{\varepsilon,t,\xi}\).Z_{0,\varepsilon,t,\xi}=DJ_\varepsilon\(W_{\varepsilon,t,\xi}\).Z_{i,\varepsilon,t,\xi}+\left\{\begin{aligned}&\operatorname{o}\(\varepsilon\delta_\varepsilon(t)^2\)&&\text{if }i=0\\&\operatorname{o}\(\varepsilon\delta_\varepsilon(t)^3\)&&\text{if }i=1,\dotsc,n\end{aligned}\right.
\end{equation}
and in case $n=4$, there holds
\begin{equation}\label{Lem3terEq9}
DJ_\varepsilon\(W_{\varepsilon,t,\xi}+\phi_{\varepsilon,t,\xi}\).Z_{0,\varepsilon,t,\xi}=DJ_\varepsilon\(W_{\varepsilon,t,\xi}\).Z_{i,\varepsilon,t,\xi}+\left\{\begin{aligned}&\operatorname{o}\(\delta_\varepsilon(t)^2\)&&\text{if }i=0\\&\operatorname{O}\(\varepsilon^2\delta_\varepsilon(t)^2\)&&\text{if }i=1,\dotsc,n\end{aligned}\right.
\end{equation}
as $\varepsilon\to0$. Finally, \eqref{Lem3terEq1} and \eqref{Lem3terEq1bis} follow from \eqref{Lem3terEq4}, \eqref{Lem3terEq8}, and \eqref{Lem3terEq9}. The estimate \eqref{Lem3terEq1ter} follows from \eqref{Lem3terEq1}--\eqref{Lem3terEq1bis} and the validity of \eqref{Lem1Eq1}--\eqref{Lem1Eq4} in a $C^1$--uniform way with respect to $\xi\in M$ and $t$ in compact subsets of $(0,\infty)$ as $\varepsilon\to0$. This ends the proof of Lemma~\ref{Lem3ter}.
\endproof
We then prove the following result.

\begin{lemma}\label{Lem3quater}
For any $i=1,\dotsc,n$, there hold
\begin{align}
\frac{d}{d\eta_i}\mathcal{J}_\varepsilon\(t,\exp_\xi\eta\)\Big|_{\eta=0}&=\frac{n^{\frac{n-2}{4}}\(n-2\)^{\frac{n+2}{4}}}{\delta_\varepsilon\(t\)}\left\|\nabla V_i\right\|_2^2\lambda_{i,\varepsilon,t,\xi}+\operatorname{o}\bigg(\sum_{j=0}^n\left|\lambda_{j,\varepsilon,t,\xi}\right|\bigg)\label{Pr3Eq9}\\
\frac{d}{dt}\mathcal{J}_\varepsilon\(t,\xi\)&=\frac{n^{\frac{n-2}{4}}\(n-2\)^{\frac{n+2}{4}}\delta'_\varepsilon\(t\)}{2\delta_\varepsilon\(t\)}\left\|\nabla V_0\right\|_2^2\lambda_{0,\varepsilon,t,\xi}+\operatorname{o}\bigg(\sum_{j=0}^n\left|\lambda_{j,\varepsilon,t,\xi}\right|\bigg)\label{Pr3Eq10}
\end{align}
as $\varepsilon\to0$, where $\delta_\varepsilon\(t\)$ is as in \eqref{Eq8}.
\end{lemma}

\proof
For any $i=1,\dotsc,n$, by \eqref{Eq27}--\eqref{Eq28}, we get that
\begin{equation}\label{Pr3Eq3}
\frac{d}{d\eta_i}\mathcal{J}_\varepsilon\(t,\exp_\xi\eta\)\Big|_{\eta=0}=\sum_{j=0}^n\lambda_{j,\varepsilon,t,\xi}\Big<Z_{j,\varepsilon,t,\xi},\frac{d}{d\eta_i}\big(W_{\varepsilon,t,\exp_\xi\eta}+\phi_{\varepsilon,t,\exp_\xi\eta}\big)\Big|_{\eta=0}\Big>_{L_g}.
\end{equation}
For any $i=1,\dotsc,n$ and $j=0,\dotsc,n$, an integration by parts gives that
\begin{multline}\label{Pr3Eq4}
\Big<Z_{j,\varepsilon,t,\xi},\frac{d}{d\eta_i}W_{\varepsilon,t,\exp_\xi\eta}\Big|_{\eta=0}\Big>_{L_g}=\int_ML_gZ_{j,\varepsilon,t,\xi}\bigg(\frac{d}{d\eta_i}W_{\varepsilon,t,\exp_\xi\eta}\Big|_{\eta=0}-\frac{n^{\frac{n-2}{4}}\(n-2\)^{\frac{n+2}{4}}}{\delta_\varepsilon\(t\)}Z_{i,\varepsilon,t,\xi}\bigg)dv_g\\
+\frac{n^{\frac{n-2}{4}}\(n-2\)^{\frac{n+2}{4}}}{\delta_\varepsilon\(t\)}\<Z_{j,\varepsilon,t,\xi},Z_{i,\varepsilon,t,\xi}\>_{L_g}.
\end{multline}
By H\"older's inequality and \eqref{Lem2Eq2}, we get that
\begin{multline}\label{Pr3Eq5}
\int_ML_gZ_{j,\varepsilon,t,\xi}\bigg(\frac{d}{d\eta_i}W_{\varepsilon,t,\exp_\xi\eta}\Big|_{\eta=0}-\frac{n^{\frac{n-2}{4}}\(n-2\)^{\frac{n+2}{4}}}{\delta_\varepsilon\(t\)}Z_{i,\varepsilon,t,\xi}\bigg)dv_g\le\left\|L_gZ_{j,\varepsilon,t,\xi}\right\|_{\frac{2n}{n+2}}\\
\times\bigg\|\frac{d}{d\eta_i}W_{\varepsilon,t,\exp_\xi\eta}\Big|_{\eta=0}-\frac{n^{\frac{n-2}{4}}\(n-2\)^{\frac{n+2}{4}}}{\delta_\varepsilon\(t\)}Z_{i,\varepsilon,t,\xi}\bigg\|_{2^*}=\operatorname{o}\(\left\|L_gZ_{j,\varepsilon,t,\xi}\right\|_{\frac{2n}{n+2}}\)=\operatorname{o}\(1\)
\end{multline}
as $\varepsilon\to0$. It follows from \eqref{Lem3terEq3} and \eqref{Pr3Eq4}--\eqref{Pr3Eq5} that
\begin{equation}\label{Pr3Eq6}
\Big<Z_{j,\varepsilon,t,\xi},\frac{d}{d\eta_i}W_{\varepsilon,t,\exp_\xi\eta}\Big|_{\eta=0}\Big>_{L_g}=\frac{n^{\frac{n-2}{4}}\(n-2\)^{\frac{n+2}{4}}}{\delta_\varepsilon\(t\)}\left\|\nabla V_i\right\|_2^2\delta_{ij}+\operatorname{o}\(1\)
\end{equation}
as $\varepsilon\to0$, where the $\delta_{ij}$'s are the Kronecker symbols. Since $\phi_{\varepsilon,t,\xi}$ belongs to $K^\perp_{\varepsilon,t,\xi}$, differentiating the equation $\<Z_{j,\varepsilon,t,\xi},\phi_{\varepsilon,t,\xi}\>_{L_g}=0$, we find that
\begin{equation}\label{Pr3Eq7}
\Big<Z_{j,\varepsilon,t,\xi},\frac{d}{d\eta_i}\phi_{\varepsilon,t,exp_\xi\eta}\Big|_{\eta=0}\Big>_{L_g}=-\Big<\frac{d}{d\eta_i}Z_{j,\varepsilon,t,exp_\xi\eta}\Big|_{\eta=0},\phi_{\varepsilon,t,\xi}\Big>_{L_g}.
\end{equation}
By \eqref{Pr3Eq7}, Cauchy--Schwarz inequality, Propositions~\ref{Pr1} and~\ref{Pr2}, we get that
\begin{align}
\Big|\Big<Z_{j,\varepsilon,t,\xi},\frac{d}{d\eta_i}\phi_{\varepsilon,t,exp_\xi\eta}\Big|_{\eta=0} \Big>_{L_g}\Big|&\le\Big\|\frac{d}{d\eta_i}Z_{j,\varepsilon,t,\exp_\xi\eta}\Big|_{\eta=0}\Big\|_{L_g}\left\|\phi_{\varepsilon,t,\xi}\right\|_{L_g}\nonumber\\
&=\operatorname{O}\big(\delta_\varepsilon\(t\)^{-1}\left\|\phi_{\varepsilon,t,\xi}\right\|_{1,2}\big)=\operatorname{o}\(1\)\label{Pr3Eq8}
\end{align}
as $\varepsilon\to0$. \eqref{Pr3Eq9} follows from \eqref{Pr3Eq3}, \eqref{Pr3Eq6}, and \eqref{Pr3Eq8}. \eqref{Pr3Eq10} follows from similar arguments by using \eqref{Lem2Eq1}. This ends the proof of Lemma~\ref{Lem3quater}.
\endproof

\proof[Proof of the $C^1$--uniformity of \eqref{Eq26} with respect to $t$]
The result follows directly from Lemmas~\ref{Lem3},~\ref{Lem3ter}, and~\ref{Lem3quater}.
\endproof

\smallskip
\proof[Proof of the $C^1$--uniformity of \eqref{Eq26} with respect to $\xi$]
In case $n\ge5$, the result follows directly from Lemmas~\ref{Lem3},~\ref{Lem3ter}, and~\ref{Lem3quater}. The case $n=4$ is trickier. In this case, the estimate \eqref{Lem3terEq1bis} is not sufficient to get such a direct proof as in higher dimensions. We prove the result in case $n=4$ in what follows. For any $i=1,\dotsc,n$ and $x\in M$, we define
$$Y_{i,\varepsilon,t,\xi}\(x\):=\chi\(d_{g_\xi}\(x,\xi\)\)\varLambda_\xi\(x\)\partial_{y_i}\big[\widetilde{W}_{\varepsilon,t,\xi}\circ\exp_\xi+\widetilde{\phi}_{\varepsilon,t,\xi}\circ\exp_\xi\big]\(\exp^{-1}_\xi x\)\,,$$
where $\widetilde{W}_{\varepsilon,t,\xi}:=W_{\varepsilon,t,\xi}/\varLambda_\xi$ and $\widetilde{\phi}_{\varepsilon,t,\xi}:=\phi_{\varepsilon,t,\xi}/\varLambda_\xi$.
We claim that
\begin{equation}\label{Pr3Eq12}
DJ_\varepsilon\(W_{\varepsilon,t,\xi}+\phi_{\varepsilon,t,\xi}\).Y_{i,\varepsilon,t,\xi}=-\frac{n^{\frac{n-2}{4}}\(n-2\)^{\frac{n+2}{4}}}{\delta_\varepsilon\(t\)}\left\|\nabla V_i\right\|_2^2\lambda_{i,\varepsilon,t,\xi}+\operatorname{o}\bigg(\sum_{j=0}^n\left|\lambda_{j,\varepsilon,t,\xi}\right|\bigg).
\end{equation}
It follows from \eqref{Eq27}--\eqref{Eq28} and \eqref{Lem3terEq3} that in order to prove \eqref{Pr3Eq12}, it suffices to prove that for any $i=1,\dotsc,n$ and $j=0,\dotsc,n$, there holds
\begin{equation}\label{Pr3Eq13}
\Big<Z_{j,\varepsilon,t,\xi},Y_{i,\varepsilon,t,\xi}+\frac{n^{\frac{n-2}{4}}\(n-2\)^{\frac{n+2}{4}}}{\delta_\varepsilon\(t\)}Z_{i,\varepsilon,t,\xi}\Big>_{L_g}=\operatorname{o}\(1\)
\end{equation}
as $\varepsilon\to0$. Integrating by parts and applying H\"older's inequality, we get that
\begin{multline}\label{Pr3Eq14}
\Big<Z_{j,\varepsilon,t,\xi},Y_{i,\varepsilon,t,\xi}+\frac{n^{\frac{n-2}{4}}\(n-2\)^{\frac{n+2}{4}}}{\delta_\varepsilon\(t\)}Z_{i,\varepsilon,t,\xi}\Big>_{L_g}=\int_ML_gZ_{j,\varepsilon,t,\xi}\Big(Y_{i,\varepsilon,t,\xi}+\frac{n^{\frac{n-2}{4}}\(n-2\)^{\frac{n+2}{4}}}{\delta_\varepsilon\(t\)}Z_{i,\varepsilon,t,\xi}\Big)dv_g\\
\le\left\|L_gZ_{j,\varepsilon,t,\xi}\right\|_{\frac{2n}{n+2}}\Big\|\chi\(d_{g_\xi}\(\cdot,\xi\)\)\varLambda_\xi\partial_{y_i}\big[\widetilde{W}_{\varepsilon,t,\xi}\circ\exp_\xi\big]\circ\exp^{-1}_\xi+\frac{n^{\frac{n-2}{4}}\(n-2\)^{\frac{n+2}{4}}}{\delta_\varepsilon\(t\)}Z_{i,\varepsilon,t,\xi}\Big\|_{2^*}\\
+\left\|L_gZ_{j,\varepsilon,t,\xi}\right\|_2\big\|\chi\(d_{g_\xi}\(\cdot,\xi\)\)\varLambda_\xi\nabla\widetilde{\phi}_{\varepsilon,t,\xi}\big\|_2\,.
\end{multline}
Using similar arguments as in the proof of Lemma~\ref{Lem2}, we find that
\begin{equation}\label{Pr3Eq15}
\Big\|\chi\(d_{g_\xi}\(\cdot,\xi\)\)\varLambda_\xi\partial_{y_i}\big[\widetilde{W}_{\varepsilon,t,\xi}\circ\exp_\xi\big]\circ\exp^{-1}_\xi+\frac{n^{\frac{n-2}{4}}\(n-2\)^{\frac{n+2}{4}}}{\delta_\varepsilon\(t\)}Z_{i,\varepsilon,t,\xi}\Big\|_{2^*}=\operatorname{o}\(1\)
\end{equation}
as $\varepsilon\to0$. Rough estimates give that
\begin{equation}\label{Pr3Eq16}
\left\|L_gZ_{j,\varepsilon,t,\xi}\right\|_{\frac{2n}{n+2}}=\operatorname{O}\(1\)\quad\text{and}\quad\left\|L_gZ_{j,\varepsilon,t,\xi}\right\|_2=\operatorname{O}\(\delta_\varepsilon\(t\)^{-1}\).
\end{equation}
Moreover, by \eqref{Eq8}, and Propositions~\ref{Pr1} and~\ref{Pr2}, we get that in dimension $n=4$, there holds
\begin{equation}\label{Pr3Eq17}
\big\|\chi\(d_{g_\xi}\(\cdot,\xi\)\)\varLambda_\xi\nabla\widetilde{\phi}_{\varepsilon,t,\xi}\big\|_2=\operatorname{O}\big(\big\|\widetilde{\phi}_{\varepsilon,t,\xi}\big\|_{1,2}\big)=\operatorname{O}(\left\|\phi_{\varepsilon,t,\xi}\right\|_{1,2})=\operatorname{O}\big(\left\|\phi_{\varepsilon,t,\xi}\right\|_{L_g}\big)=\operatorname{O}\(\varepsilon\delta_\varepsilon\(t\)\).
\end{equation}
\eqref{Pr3Eq13} follows from \eqref{Pr3Eq14}--\eqref{Pr3Eq17}, and thus we get \eqref{Pr3Eq12}. It follows from \eqref{Lem2Eq1}, \eqref{Pr3Eq12}, and Lemmas~\ref{Lem3ter} and~\ref{Lem3quater} that in dimension $n=4$, there holds
\begin{equation}\label{Pr3Eq18}
\frac{d}{d\eta_i}\mathcal{J}_\varepsilon\(t,\exp_\xi\eta\)\Big|_{\eta=0}=-DJ_\varepsilon\(W_{\varepsilon,t,\xi}+\phi_{\varepsilon,t,\xi}\).Y_{i,\varepsilon,t,\xi}+\operatorname{o}\(\delta_\varepsilon\(t\)^2 \ln |\delta_\varepsilon(t)| \)
\end{equation}
as $\varepsilon\to0$. By the conformal change of metric $g_\xi=\varLambda_\xi^{2^*-2}g$, we get that
\begin{multline}\label{Pr3Eq19}
DJ_\varepsilon\(W_{\varepsilon,t,\xi}+\phi_{\varepsilon,t,\xi}\).Y_{i,\varepsilon,t,\xi}=\int_M\Big(\big<\nabla\big(\widetilde{W}_{\varepsilon,t,\xi}+\widetilde{\phi}_{\varepsilon,t,\xi}\big),\nabla \widetilde{Y}_{i,\varepsilon,t,\xi}\big>_{g_\xi}\\
+\(\alpha_n\Scal_{g_\xi}+\varepsilon\varLambda_\xi^{2-2^*}h\)\big(\widetilde{W}_{\varepsilon,t,\xi}+\widetilde{\phi}_{\varepsilon,t,\xi}\big)\widetilde{Y}_{i,\varepsilon,t,\xi}-\big(\widetilde{W}_{\varepsilon,t,\xi}+\widetilde{\phi}_{\varepsilon,t,\xi}\big)_+^{2^*-1}\widetilde{Y}_{i,\varepsilon,t,\xi}\Big)dv_{g_\xi}\,,
\end{multline}
where $\widetilde{Y}_{i,\varepsilon,t,\xi}:=Y_{i,\varepsilon,t,\xi}/\varLambda_\xi$. Since $\widetilde{Y}_{i,\varepsilon,t,\xi}\equiv0$ in $M\backslash B_\xi\(r_0\)$, letting $y=\exp_\xi^{-1}x$ and integrating by parts \eqref{Pr3Eq19}, we find that
\begin{equation}\label{Pr3Eq20}
DJ_\varepsilon\(W_{\varepsilon,t,\xi}+\phi_{\varepsilon,t,\xi}\).Y_{i,\varepsilon,t,\xi}=I_{1,i,\varepsilon,t,\xi}+I_{2,i,\varepsilon,t,\xi}+I_{3,i,\varepsilon,t,\xi}+I_{4,i,\varepsilon,t,\xi}\,,
\end{equation}
where
\begin{align*}
I_{1,i,\varepsilon,t,\xi}&:=DJ_\varepsilon\(W_{\varepsilon,t,\xi}\).\(\chi\(d_{g_\xi}\(\cdot,\xi\)\)\varLambda_\xi\partial_{y_i}\big[\widetilde{W}_{\varepsilon,t,\xi}\circ\exp_\xi\big]\circ\exp^{-1}_\xi\),\allowdisplaybreaks\\
I_{2,i,\varepsilon,t,\xi}&:=\int_{B_0\(r_0\)}\Big\{ g_\xi^{pq}\partial_{y_p}\chi \partial_{y_q}\big[\widetilde{\phi}_{\varepsilon,t,\xi}\circ\exp_\xi\big]\partial_{y_i}\big[\widetilde{W}_{\varepsilon,t,\xi}\circ\exp_\xi\big]\left|\exp_\xi^*g_\xi\right|  \\
&\quad+g_\xi^{pq}\partial_{y_p} \chi \partial_{y_q}\big[\widetilde{W}_{\varepsilon,t,\xi}\circ\exp_\xi\big]\partial_{y_i}\big[\widetilde{\phi}_{\varepsilon,t,\xi}\circ\exp_\xi\big]\left|\exp_\xi^*g_\xi\right|\\
&\quad-\chi \partial_{y_i}\big[g_\xi^{pq}\big]\partial_{y_p}\big[\widetilde{W}_{\varepsilon,t,\xi}\circ\exp_\xi\big]\partial_{y_q}\big[\widetilde{\phi}_{\varepsilon,t,\xi}\circ\exp_\xi\big]\left|\exp_\xi^*g_\xi\right|\\
&\quad-\chi \big(\widetilde{W}_{\varepsilon,t,\xi}\circ\exp_\xi\big)\big(\widetilde{\phi}_{\varepsilon,t,\xi}\circ\exp_\xi\big)\partial_{y_i}\big[\(\alpha_n\Scal_{g_\xi}+\varepsilon\varLambda_\xi^{2-2^*}h\)\circ\exp_\xi\big]\left|\exp_\xi^*g_\xi\right|\\
&\quad-\Big(\big<\nabla\widetilde{W}_{\varepsilon,t,\xi},\nabla\widetilde{\phi}_{\varepsilon,t,\xi}\big>_{g_\xi}+\(\alpha_n\Scal_{g_\xi}+\varepsilon\varLambda_\xi^{2-2^*}h\)\widetilde{W}_{\varepsilon,t,\xi}\widetilde{\phi}_{\varepsilon,t,\xi}-\widetilde{W}_{\varepsilon,t,\xi}^{2^*-1}\widetilde{\phi}_{\varepsilon,t,\xi}\Big)\circ\exp_\xi\\
&\qquad \times\partial_{y_i}\big[\chi \left|\exp_\xi^*g_\xi\right|\big]\Big\}dy\,,\allowdisplaybreaks\\
I_{3,i,\varepsilon,t,\xi}&:=\int_{B_0\(r_0\)}\Big\{ g_\xi^{pq}\partial_{y_p} \chi \partial_{y_q}\big[\widetilde{\phi}_{\varepsilon,t,\xi}\circ\exp_\xi\big]\partial_{y_i}\big[\widetilde{\phi}_{\varepsilon,t,\xi}\circ\exp_\xi\big]\left|\exp_\xi^*g_\xi\right|\\
&\quad-\frac{1}{2}\chi \partial_{y_i}\big[g_\xi^{pq}\big]\partial_{y_p}\big[\widetilde{\phi}_{\varepsilon,t,\xi}\circ\exp_\xi\big]\partial_{y_q}\big[\widetilde{\phi}_{\varepsilon,t,\xi}\circ\exp_\xi\big]\left|\exp_\xi^*g_\xi\right|\\
&\quad-\frac{1}{2}\chi \big(\widetilde{\phi}_{\varepsilon,t,\xi}\circ\exp_\xi\big)^2\partial_{y_i}\big[\(\alpha_n\Scal_{g_\xi}+\varepsilon\varLambda_\xi^{2-2^*}h\)\circ\exp_\xi\big]\left|\exp_\xi^*g_\xi\right|\\
&\quad  -\frac{1}{2}\Big(\big|\nabla\widetilde{\phi}_{\varepsilon,t,\xi}\big|_{g_\xi}^2+\(\alpha_n\Scal_{g_\xi}+\varepsilon\varLambda_\xi^{2-2^*}h\)\widetilde{\phi}_{\varepsilon,t,\xi}^2\Big)\circ\exp_\xi  \partial_{y_i}\big[\chi \left|\exp_\xi^*g_\xi\right|\big]\Big\}  dy\,,\allowdisplaybreaks\\
I_{4,i,\varepsilon,t,\xi}&:=\frac{1}{2^*}\int_{B_0\(r_0\)}\Big\{ \Big( \big(\widetilde{W}_{\varepsilon,t,\xi}+\widetilde{\phi}_{\varepsilon,t,\xi}\big)_+^{2^*}-\widetilde{W}_{\varepsilon,t,\xi}^{2^*}-2^*\widetilde{W}_{\varepsilon,t,\xi}^{2^*-1}\widetilde{\phi}_{\varepsilon,t,\xi}\Big)\circ\exp_\xi \\
&\qquad  \times\partial_{y_i}\big[\chi \left|\exp_\xi^*g_\xi\right|\big] \Big\} dy\,,
\end{align*}
where $g_\xi^{pq}$ are the components of $g_\xi^{-1}$ in geodesic normal coordinates. Using similar arguments as in the proof of Lemma~\ref{Lem3}, we find that
\begin{equation}\label{Pr3Eq21}
I_{1,i,\varepsilon,t,\xi}=-\frac{d}{d\eta_i}J_\varepsilon\(W_{\varepsilon,t,\xi}\)\Big|_{\eta=0}+\operatorname{o}\(\delta_\varepsilon\(t\)^2\)
\end{equation}
as $\varepsilon\to0$. Now, we estimate $I_{2,i,\varepsilon,t,\xi}$. Since $g_\xi$ defines conformal normal coordinates of order $N\ge5$, see Lee--Parker~\cite{LeePar}*{Theorem~5.1}, we get that
\begin{align} \label{Pr3Eq22}
\partial_{y_i}\big[\Scal_{g_\xi}\circ\exp_\xi\big]\(y\)=\operatorname{O}\(\left|y\right|\),\qquad \partial_{y_i}\big[\left|\exp_\xi^*g_\xi\right|\big]\(y\)=\operatorname{O}\big(\left|y\right|^{N-1}\big).
\end{align}
Since $\partial_{y_i}\big[g_\xi^{pq}\big]\(y\)=\operatorname{O}\(\left|y\right|\)$, by using H\"older's inequality and \eqref{Pr3Eq22}, we find that
\begin{multline}\label{Pr3Eq25}
I_{2,i,\varepsilon,t,\xi}=\operatorname{O}\Big(\big\|d_{g_\xi}\(\cdot,\xi\)\nabla\widetilde{W}_{\varepsilon,t,\xi}\big\|_2\big\|\nabla\widetilde{\phi}_{\varepsilon,t,\xi}\big\|_2+\Big(\big\|d_{g_\xi}\(\cdot,\xi\)\widetilde{W}_{\varepsilon,t,\xi}\big\|_{\frac{2n}{n+2}}\\
+\varepsilon\big\|\widetilde{W}_{\varepsilon,t,\xi}\big\|_{\frac{2n}{n+2}}+\big\|d_{g_\xi}\(\cdot,\xi\)^{N-1}\widetilde{W}_{\varepsilon,t,\xi}^{2^*-1}\big\|_{\frac{2n}{n+2}}\Big)\big\|\widetilde{\phi}_{\varepsilon,t,\xi}\big\|_{2^*}\Big).
\end{multline}
In dimension $n=4$, rough estimates give that
\begin{multline}\label{Pr3Eq26}
\big\|d_{g_\xi}\(\cdot,\xi\)\nabla\widetilde{W}_{\varepsilon,t,\xi}\big\|_2=\operatorname{O}\big(\delta_\varepsilon\(t\)\sqrt{\left|\ln\delta_\varepsilon\(t\)\right|}\big),\quad\big\|d_{g_\xi}\(\cdot,\xi\)\widetilde{W}_{\varepsilon,t,\xi}\big\|_{\frac{2n}{n+2}}=\operatorname{O}\big(\delta_\varepsilon\(t\)\big),\\
\big\|\widetilde{W}_{\varepsilon,t,\xi}\big\|_{\frac{2n}{n+2}}=\operatorname{O}\big(\delta_\varepsilon\(t\)\big),\quad\text{and}\quad\big\|d_{g_\xi}\(\cdot,\xi\)^{N-1}\widetilde{W}_{\varepsilon,t,\xi}^{2^*-1}\big\|_{\frac{2n}{n+2}}=\operatorname{O}\big(\delta_\varepsilon\(t\)^3\big).
\end{multline}
It follows from \eqref{Pr3Eq25}--\eqref{Pr3Eq26} and from Sobolev's embedding $ H_1^2(M) \hookrightarrow L^{2^*}(M)$ that
\begin{equation}\label{Pr3Eq27}
I_{2,i,\varepsilon,t,\xi}=\operatorname{O}\Big(\delta_\varepsilon\(t\)\sqrt{\left|\ln\delta_\varepsilon\(t\)\right|}\big\|\widetilde{\phi}_{\varepsilon,t,\xi}\big\|_{1,2}\Big).
\end{equation}
It remains to estimate $I_{3,i,\varepsilon,t,\xi}$ and $I_{4,i,\varepsilon,t,\xi}$. Clearly, we get that
\begin{equation}\label{Pr3Eq28}
I_{3,i,\varepsilon,t,\xi}=\operatorname{O}\Big(\big\|\widetilde{\phi}_{\varepsilon,t,\xi}\big\|_{1,2}^2\Big).
\end{equation}
Similarly as in \eqref{Pr3Step1Eq3}, we get that
\begin{equation}\label{Pr3Eq29}
I_{4,i,\varepsilon,t,\xi}=\operatorname{O}\Big(\big\|\widetilde{\phi}_{\varepsilon,t,\xi}\big\|_{2^*}^2+\big\|\widetilde{\phi}_{\varepsilon,t,\xi}\big\|_{2^*}^{2^*}\Big).
\end{equation}
Since $\left\|\phi_{\varepsilon,t,\xi}\right\|_{1,2}=\operatorname{O}\big(\left\|\phi_{\varepsilon,t,\xi}\right\|_{L_g}\big)=\operatorname{O}\(\varepsilon\delta_\varepsilon\(t\)\)$, by \eqref{Pr3Eq27}--\eqref{Pr3Eq29}, and by Sobolev's embedding $H_1^2(M) \hookrightarrow L^{2^*}(M)$, we get that
\begin{equation}\label{Pr3Eq30}
I_{2,i,\varepsilon,t,\xi}+I_{3,i,\varepsilon,t,\xi}+I_{4,i,\varepsilon,t,\xi}=\operatorname{O}\(\sqrt\varepsilon\delta_\varepsilon\(t\)^2\)=\operatorname{o}\(\delta_\varepsilon\(t\)^2\)
\end{equation}
as $\varepsilon\to0$. Finally, it follows from \eqref{Pr3Eq18}, \eqref{Pr3Eq20}, \eqref{Pr3Eq21}, and \eqref{Pr3Eq30} that \eqref{Eq26} is $C^1$--uniform with respect to $\xi$.
\endproof

{\section*{Appendix}
\setcounter{theorem}{0}
\setcounter{equation}{0}
\renewcommand{\thesection}{A}

In this appendix, first, we state the following result which is due to Lee--Parker~\cite{LeePar}.

\begin{lemma}\label{Lem4}
For any $\xi\in M$, let $g_\xi$ be as in Section~\ref{Sec2} and $G_{g_\xi}$ be the Green's function of the conformal Laplacian with respect to $g_\xi$. Then the following expansions do hold:
\begin{enumerate}
\renewcommand{\labelenumi}{(\roman{enumi})}
\item when $n=4,5$ or $\(M,g\)$ is l.c.f.,
$$G_{g_\xi}\(\exp_\xi y,\xi\)=\beta_n^{-1}\left|y\right|^{2-n}+A_\xi+\operatorname{O}\(\left|y\right|\),$$
\item when $n=6$,
$$G_{g_\xi}\(\exp_\xi y,\xi\)=\beta_6^{-1}\(\left|y\right|^{-4}-\frac{1}{1440}\big|\Weyl_g\(\xi\)\big|_g^2\ln\left|y\right|\)+\operatorname{O}\(\left|y\right|\),$$
\item when $n\ge7$,
\begin{align*}
G_{g_\xi}\(\exp_\xi y,\xi\)&=\beta_n^{-1}\left|y\right|^{2-n}\bigg(1+\frac{\alpha_n\left|y\right|^2}{12\(n-4\)}\bigg(\frac{1}{12\(n-6\)}\big|\Weyl_g\(\xi\)\big|_g^2\left|y\right|^2\nonumber\\
&\qquad\qquad-\partial_{y_iy_j}\[\Scal_{g_\xi}\circ\exp_\xi\]\(0\)y_iy_j\bigg)\bigg)+\operatorname{O}\(\left|y\right|^{7-n}\)
\end{align*}
$C^1$--uniformly with respect to $y$ and $\xi$, where $\Weyl_g$ is the Weyl curvature tensor with respect to $g$, and for any $\xi\in M$, $\Scal_{g_\xi}$ is the scalar curvature with respect to $g_\xi$.
\end{enumerate}
\end{lemma}

\proof
We refer to Lee--Parker~\cite{LeePar}*{Lemma~6.4}. The only point which is not discussed in~\cite{LeePar} is the $C^1$--uniformity of the expansions with respect to $\xi$. This point follows from standard arguments.
\endproof

Finally, let us estimate the first derivatives of the geodesic distance in the following result.
\begin{lemma}\label{Lem5}
For any $\xi\in M$, let $g_\xi$ be as in Section~\ref{Sec2} and $d_{g_\xi}$ be the geodesic distance with respect to $g_\xi$. Then there holds
\begin{equation}\label{Lem5Eq1}
\nabla\big(d_{g_{\exp_\xi\eta}}\(\exp_\xi y,\exp_\xi\eta\)\big)\big|_{\eta=0}=-\frac{y}{\left|y\right|}+\operatorname{O}\(\left|y\right|^2\)
\end{equation}
uniformly with respect to $\xi\in M$ and $y\in T_\xi M\cong\mathbb{R}^n$, $y\ne0$, $\left|y\right|\ll1$.
\end{lemma}

\proof
By compactness of $M$ and since $g_\xi$ is smooth with respect to $\xi$, we get that there exists a positive real number $r_0$ such that $r_0<i_{g_\xi}$ for all $\xi\in M$, where $i_{g_\xi}$ is the injectivity radius of the manifold $\(M,g_\xi\)$. For any $\xi\in M$, we let $B_0\(r_0\)$ be the ball in $T_\xi M\cong\mathbb{R}^n$ of center 0 and radius $r_0$. For any $\xi\in M$ and $y\in B_0\(r_0\)\backslash\left\{0\right\}$, we write that
\begin{equation}\label{Lem5Eq2}
\nabla\big(d_{g_{\exp_\xi\eta}}\(\exp_\xi y,\exp_\xi\eta\)\big)\big|_{\eta=0}=\nabla\big(d_{g_{\exp_\xi\eta}}\(\exp_\xi y,\xi\)\big)\big|_{\eta=0}+\nabla\big(d_{g_\xi}\(\exp_\xi y,\exp_\xi\eta\)\big)\big|_{\eta=0}\,.
\end{equation}
We begin with estimating the first term in the right hand side of \eqref{Lem5Eq2}. Since we have chosen $\varLambda_\xi$ so that $\varLambda_\xi\(\xi\)=1$ and $\nabla\varLambda_\xi\(\xi\)=0$, we get that
$$\varLambda_{\exp_\xi\eta}\(\exp_\xi y\)=\varLambda_\xi\(\exp_\xi y\)\(1+\operatorname{O}\(\left|y\right|^2+\left|\eta\right|^2\)\)$$
uniformly with respect to $\xi\in M$ and $y,\eta\in B_0\(r_0\)$. It follows that
$$d_{g_{\exp_\xi\eta}}\(\exp_\xi y,\xi\)=d_{g_\xi}\(\exp_\xi y,\xi\)\(1+\operatorname{O}\(\left|y\right|^2+\left|\eta\right|^2\)\)=\left|y\right|+\operatorname{O}\(\left|y\right|^3+\left|y\right|\left|\eta\right|^2\),$$
and thus we get that
\begin{equation}\label{Lem5Eq3}
\nabla\big(d_{g_{\exp_\xi\eta}}\(\exp_\xi y,\xi\)\big)\big|_{\eta=0}=\operatorname{O}\(\left|y\right|^2\)
\end{equation}
uniformly with respect to $\xi\in M$ and $y\in B_0\(r_0\)$. Now, we estimate the second term in the right hand side of \eqref{Lem5Eq2}. For any $\xi\in M$ and $\eta\in T_\xi M$, we let $T_{\xi,\exp_\xi\eta} M$ and $\exp_{\xi,\exp_\xi\eta}$ be the respective tangent space and exponential map at $\exp_\xi\eta$ with respect to $g_\xi$. For any $\eta\in B_0\(r_0\)$, we identify $T_{\xi,\exp_\xi\eta}M$ with $\mathbb{R}^n$ thanks to a local orthonormal frame, parallel at $\xi$. For any $i=1,\dotsc,n$, $\xi\in M$, and $y\in B_0\(r_0\)\backslash\left\{0\right\}$, we get that
\begin{equation}\label{Lem5Eq4}
\frac{d}{d\eta_i}\big(d_{g_\xi}\(\exp_\xi y,\exp_\xi\eta\)\big)\Big|_{\eta=0}=\sum_{j=1}^n\frac{y_j}{\left|y\right|}\frac{d}{d\eta_i}\big(\exp_{\xi,\exp_\xi\eta}^{-1}\exp_\xi y\big)_j\Big|_{\eta=0}\,.
\end{equation}
We claim that for any $i,j=1,\dotsc,n$, there holds
\begin{equation}\label{Lem5Eq5}
\frac{d}{d\eta_i}\big(\exp_{\xi,\exp_\xi\eta}^{-1}\exp_\xi y\big)_j\Big|_{\eta=0}=-\delta_{ij}+\operatorname{O}\(\left|y\right|^2\)
\end{equation}
uniformly with respect to $\xi\in M$ and $y\in B_0\(r_0\)$, where the $\delta_{ij}$'s are the Kronecker symbols. We prove this claim. For any $j=1,\dotsc,n$, $\xi\in M$ and $\eta,y\in B_0\(r_0\)$, we define
$$\mathcal{E}_{j,\xi}\(y,\eta\):=\big(\exp_{\xi,\exp_\xi\eta}^{-1}\exp_\xi y\big)_j\,.$$
Clearly, $\mathcal{E}_{j,\xi}\(y,\eta\)$ is smooth with respect to $\xi$, $y$, and $\eta$. In order to prove the Taylor expansion \eqref{Lem5Eq5}, we compute the first and second order derivatives of $\mathcal{E}_{j,\xi}\(y,\eta\)$ with respect to $y$ and $\eta$. Since the frame is parallel at $\xi$, we get that $\mathcal{E}_{j,\xi}\(0,\eta\)=-\eta_j$ for all $\eta\in B_0\(r_0\)$. Differentiating this equation gives that
\begin{equation}\label{Lem5Eq6}
\partial_{\eta_i}\mathcal{E}_{j,\xi}\(0,0\)=-\delta_{ij}\quad\text{and}\quad\partial_{\eta_i\eta_k}\mathcal{E}_{j,\xi}\(0,0\)=0
\end{equation}
for all $i,j,k=1,\dotsc,n$. We also remark that $\mathcal{E}_{j,\xi}\(y,0\)=y_j$ for all $y\in B_0\(r_0\)$, and thus we get that
\begin{equation}\label{Lem5Eq7}
\partial_{y_i}\mathcal{E}_{j,\xi}\(0,0\)=\delta_{ij}\quad\text{and}\quad\partial_{y_iy_k}\mathcal{E}_{j,\xi}\(0,0\)=0
\end{equation}
for all $i,j,k=1,\dotsc,n$. As a third equation, we find $\mathcal{E}_{j,\xi}\(y,y\)=0$ for all $y\in B_0\(r_0\)$. Differentiating this equation and using \eqref{Lem5Eq6} and \eqref{Lem5Eq7}, we find that
\begin{equation}\label{Lem5Eq8}
\partial_{\eta_iy_k}\mathcal{E}_{j,\xi}\(0,0\)=-\frac{1}{2}\(\partial_{\eta_i\eta_k}\mathcal{E}_{j,\xi}\(0,0\)+\partial_{y_iy_k}\mathcal{E}_{j,\xi}\(0,0\)\)=0
\end{equation}
for all $i,j,k=1,\dotsc,n$. \eqref{Lem5Eq5} follows from \eqref{Lem5Eq6} and \eqref{Lem5Eq8}. Finally, \eqref{Lem5Eq1} follows from \eqref{Lem5Eq2}--\eqref{Lem5Eq5}. This ends the proof of Lemma~\ref{Lem5}.
\endproof}

\end{document}